\documentclass[10pt]{article}
\usepackage[leqno]{amsmath}   
\usepackage{amsthm}   
\usepackage{amssymb}  
\usepackage{amscd}   
\usepackage{cdextra}

\newcommand{\free}{\operatorname{free}}
\newcommand{\Col}{\operatorname{Col}}
\newcommand{\DR}{\mathrm{DR}}
\renewcommand{\det}{\operatorname{det}}

\newcommand{\cris}{{\mathrm{cris}}}
\newcommand{\hm}{H_{\Mh}}
\newcommand{\Ind}{\operatorname{Ind}}

\newcommand{\loc}{\mathrm{loc}}
\newcommand{\gl}{\mathrm{gl}}

\newcommand{\zetatilde}{\tilde{\zeta}}
\renewcommand{\Bbb}{\mathbb}  

\newcommand{\Q}{{\Bbb{Q}}}  
\newcommand{\R}{{\Bbb{R}}}  
\newcommand{\C}{{\Bbb{C}}}  
\newcommand{\Z}{{\Bbb{Z}}}  

\newcommand{\F}{{\Bbb{F}}}  

\renewcommand{\bar}{\overline}

\newcommand{\tensor}{\otimes}   
\newcommand{\isom}{\cong}       
\newcommand{\ohne}{\smallsetminus}
\newcommand{\Hom}{\operatorname{Hom}}

\newcommand{\coker}{\operatorname{Coker}}
\newcommand{\prolim}{\varprojlim}
\newcommand{\id}{\operatorname{id}}

\newcommand{\bew}{\begin{proof}}
\newcommand{\bewende}{\end{proof}}
\newtheorem{lemma}{Lemma}[subsection]
\newtheorem{prop}[lemma]{Proposition}
\newtheorem{thm}[lemma]{Theorem}
\newtheorem{defn}[lemma]{Definition}
\newtheorem{cor}[lemma]{Corollary}
\newcommand{\bem}{\noindent{\bf Remark:\ }}
\newcommand{\rem}{\noindent{\bf Remark:\ }}
\newcommand{\example}{\noindent{\bf Example:\ }}
\newtheorem{conj}[lemma]{Conjecture}

\newtheorem{alemma}{Lemma}[section]
\newtheorem{aprop}[alemma]{Proposition}
\newtheorem{acor}[alemma]{Corollary}

\newcommand{\spec}{\operatorname{Spec}} 

\newcommand{\cone}{\operatorname{Cone}}

\newcommand{\vol}{\operatorname{vol}}

\newcommand{\cycl}{\operatorname{cycl}}
\newcommand{\tors}{\operatorname{tors}}

\newcommand{\mal}{\times}

\newcommand{\Hyp}{\operatorname{Hyp}}
\newcommand{\ind}{\operatorname{ind}}

\newcommand{\rk}{\operatorname{rk}}

\newcommand{\af}{\frak{a}}

\newcommand{\Tr}{\operatorname{Tr}}

\newcommand{\car}{\operatorname{char}}

\newcommand{\Fr}{\operatorname{Fr}}

\newcommand{\Qbar}{\overline{\Q}}

\newcommand{\ctilde}{\tilde{c}}

\newcommand{\pf}{\frak{p}}

\newcommand{\Hh}{{\bf H}}

\newcommand{\Rh}{{\bf R}}
\newcommand{\Oh}{{\cal O}}
\newcommand{\Dh}{{\cal D}}
\newcommand{\Mh}{{\cal M}}

\newcommand{\Lh}{{\cal L}}

\renewcommand{\tilde}{\widetilde}

\newcommand{\tr}{\operatorname{tr}} 

\newcommand{\Norm}{\operatorname{N}}

\renewcommand{\Bbb}{\mathbb}  
\newcommand{\Lbb}{{\Bbb{L}}}
\renewcommand{\epsilon}{\varepsilon}
\renewcommand{\rho}{\varrho}

\newtheorem{mainconjecture}[lemma]{Main Conjecture}

\newcommand{\cores}{\operatorname{cores}}

\newcommand{\Gal}{\operatorname{Gal}}

\begin{document}

\begin{center}
{\bf\Large
Bloch-Kato Conjecture and Main Conjecture of Iwasawa Theory for
Dirichlet Characters
}\\[2ex]
revised version, February 2002\\
by Annette Huber and Guido Kings\\[3ex]
\end{center}
\begin{tabular}{ll}
Univ. Leipzig&Univ. Regensburg\\
Math. Institut &NWF-I Mathematik\\
PF 920&93040 Regensburg\\
04009 Leipzig&Germany \\
Germany&guido.kings@mathematik.uni-regensburg.de\\
huber@mathematik.uni-leipzig.de &
\end{tabular}

\section*{Introduction}
The Tamagawa number conjecture proposed by Bloch and Kato in
\cite{Bloch-Kato} describes the ``special values'' of $L$-functions
of motives in terms of cohomological data. Special
value means here the leading coefficient of the Taylor series 
of an $L$-function at an integral point.

The main conjecture of Iwasawa theory describes a $p$-adic $L$-function
in terms of the structure of Galois-modules, more precisely modules
for the Iwasawa algebra. We prove:

\noindent{\bf Theorem \ref{mcthm}, \ref{mc2true}, \ref{maintheorem}:}
{\it 
Let $\chi$ be a Dirichlet character. Then the Main Conjecture of
Iwasawa theory holds for $\chi$ and $p\neq 2$ and the Bloch-Kato
conjecture holds up ot powers of $2$ for $L(\chi,r)$ with $r\in\Z$.
}

Already in the simplest case of the Tamagawa number conjecture -- the
Riemann $\zeta$-function (see \cite{Bloch-Kato} Theorem 6.1) -- the
Main Conjecture due to  Mazur and Wiles
is used by Bloch and Kato in an essential way. It was clear to the
experts that  the Main conjecture for a Dirichlet character and 
a prime $p$ should be used in order to 
prove the $p$-part of the
Bloch-Kato conjecture (see e.g.  the list in Fontaine's
Bourbaki talk \cite{FontBou} {\S} 10 and \cite{Perrin-Riou-Asterisque}). 
There are two versions of the Main Conjecture depending
on parity.
Both are needed in order to deduce the Bloch-Kato conjecture in general. 
The first version involving the $p$-adic $L$-function 
 was proved by Mazur and Wiles.
 Under the condition that $p$ does not divide
$\Phi(N)$, the second version involving cyclotomic units follows 
 from the Mazur-Wiles case. Under the same restriction, it was also proved directly 
by Rubin
using Kolyvagin's Euler systems.
The missing cases need to be addressed.

 
What we want to advocate strongly is
the insight, due to Kato in \cite{galaxy} and \cite{Kato Main conj} and
Perrin-Riou \cite{Perrin-Riou-Asterisque}, that the 
Bloch-Kato conjecture and the Main Conjecture
are  two incarnations  of  the same mathematical content.
Knowing the $p$-part of the Bloch-Kato conjecture for all fields in the 
cyclotomic tower is equivalent to the Main Conjecture formulated 
appropriately. Our approach is to take Kato's viewpoint seriously and 
to prove the Main Conjecture without any restriction on the prime $p$
from the Bloch-Kato conjecture. 
As there is one basic case in which the Bloch-Kato conjecture is known --
the analytic class number formula for $\zeta_F(0)$ where $F$ is
a number field -- we can start a bootstrapping process which in the
end proves both conjectures for all Dirichlet characters.

The elements of the proof are not really new, but basically  
copy  Rubin's proof using Euler systems. 
What is new is the right formulation in which 
all problems with bad primes disappear. What is the problem? If $p$
divides $\Phi(N)$, then the projectors to $\chi$-eigenparts are
not $\Z_p$-integral. It is not possible to decompose for example the
class group of $\Q(\mu_N)$ into $\chi$-eigenparts.  Following Kato's
idea, we formulate the Main conjecture as it is dictated by the Bloch-Kato
conjecture. It turns out that the ``zeta elements'' of Kato
lead to Euler systems, an observation already due to Kato.
We then follow Rubin's proof using the machinery of Euler systems developed
by Rubin and independently by Perrin-Riou and Kato.
This general machinery leads to a divisibility statement for the
characteristic power series of the modules involved. To get equality,
 the classical class number trick is replaced by the same trick
with the Bloch-Kato conjecture for the Dedekind-$\zeta$-function. 

The use of the Bloch-Kato formulation has two important 
advantages over the more classical approach to Iwasawa theory. On the
one hand it explains clearly why, depending on parity, cyclotomic units or the $p$-adic $L$-function
appear in the formulation of the Main Conjecture. The reason is precisely
that they interpret the $L$-value of the complex $L$-function in the
``fundamental line''. On the other hand the Bloch-Kato conjecture
is ``rational'', so that a decomposition into characters is possible, i.e.,
one can formulate the conjecture for the motives attached to Dirichlet 
characters.
 The isogeny invariance of the formulation allows to choose 
the lattices in the Galois representations according to our needs. Thus we
do not have to decompose a fixed lattice into eigenspaces of the 
characters. This is an advantage over the classical situation and
allows to resolve the technical issues connected with the decomposition
into eigenspaces.

The proof uses very few ingredients and a  minimum of  computations.
As indicated earlier, we do not use 
any results on the Main Conjecture from the
literature. We see this also as an advantage, as the literature is not
free of mistakes and the computations with cyclotomic units are hard to 
follow. 
Here is the list of key ingredients:
\begin{enumerate}
\item standard results for Galois cohomology like Poitou-Tate duality;
\item the analytic class number formula for Dedekind-$\zeta$-functions;
\item the explicit reciprocity law for $\Z(r)$ with $r\geq 2$ over unramified cyclotomic fields as proved by Kato \cite{galaxy} or Perrin-Riou \cite{PerrinEuler};
\item the theory of Euler systems due to Kolyvagin, Rubin, Kato and Perrin-Riou;
\item a comparison result on the image of cyclotomic elements in Deligne-cohomology and $p$-adic cohomology known as Bloch-Kato Conjecture 6.2 from \cite{Bloch-Kato} (proved in \cite{Hu-Wi}
following Beilinson and Deligne; a second proof is
given in \cite{Hu-Ki}).
\end{enumerate}

Shortly before finishing the first version of this paper, we learned at
the Obernai conference that Burns and Greither had been working
independently on nearly the same problem. They  prove 
the equivariant Tamagawa number conjecture for the Tate-motive $\Q(r)$ and
$r\leq 0$ over
abelian number fields. This also allows to deduce Tamagawa number
formulas for the $L$-functions of Dirichlet characters at
negative integers.  Their proof
is quite different, for example they use what was previously
known on the Main conjecture and difficult computations of $\mu$-invariants.

The Bloch-Kato conjecture for all Dirichlet characters immediately implies
the Bloch-Kato conjecture for Dedekind-$\zeta$-functions of abelian number
fields. For negative values this amounts to the cohomological
Lichtenbaum conjecture, see theorem \ref{lichtenbaumconj}.
The reduction of the conjecture at positive values to the Lichtenbaum 
conjecture was shown independently by
Benois and Nguyen Quang Do (\cite{BeNg}). 
Kolster, Nguyen
Quang Do and Fleckinger consider the latter in \cite{Kolster}. Together with
 the
correction of Euler factors worked out by Kolster and Nguyen Quang Do in
\cite{KoNg}, they prove the cohomological Lichtenbaum conjecture for abelian number fields up
to an explicit set of bad primes. The problem are the bad primes in the
above mentioned second version of the Main Conjecture. None of the
references given in \cite{Kolster} 5.2, nor \cite{BelNg} 3.2 (quoted in the \cite{BeNg} A.2.3)
seems to address this point correctly.

\noindent{\bf Overview:}
The text is organized as follows: In section 1 we state the Bloch-Kato
conjecture for all Artin motives. This of course includes the
case of motives attached to Dirichlet characters.

Section 2 starts with the proof of the Bloch-Kato conjecture by establishing
compatibility with the functional equation and the special cases
$r=0$ and $r=1$ for Dedekind-$\zeta$-functions. In section 3 facts
about cyclotomic elements and their relation with the Bloch-Kato conjecture
are assembled.

Most of section 4 is independent from the previous results. The Main
conjecture is stated and proved. Finally, the Bloch-Kato conjecture
is proved from the Main conjecture in section 5.
The last section 6 establishes a few facts about certain Iwasawa modules.
These results are needed in the discussion of the Main conjecture.

Appendix
A proves a well-known lemma on $p$-adic cohomology of local fields  but 
where we did not find a reference. 
Appendix B establishes compatibility of the conjecture 
with the functional equation.
This would also be a consequence of an unpublished result of Kato \cite{epsilon}.

\noindent{\bf Acknowledgments:} We are thankful to D. Blasius for
suggesting that this would be a worthwhile project. 
We would also like to thank
D.~Benois and P.~Colmez, who answered our questions on explicit reciprocity laws, and P. Schneider
for help with $p$-adic cohomology. 
We appreciate very much the comments by
D.~Benois, Th.~Nguyen Quang Do and the referees on the first version of the
paper.

\section{The Bloch-Kato conjecture for Artin motives}
Our first aim is to present a formulation of the Bloch-Kato conjecture
for Artin motives and their $L$-values at integral points. We follow
Fontaine's approach in \cite{FontBou}. 

\subsection{Artin motives, realizations and regulators}\label{artinmotives}
We work over the base field $\Q$ and with coefficients in some
number field $E$. Let $G_{\Q}$ be the absolute Galois group of $\Q$ 
and $\Oh$ be the ring of integers of $E$. An 
{\em Artin motive
over $\Q$ with coefficients in $E$} is a direct summand of the motive
of a number field in the category
of all Grothendieck motives with coefficients in $E$. They form a well-defined
abelian category. The {\em dual motive}
of $V$ is denoted $V^\lor$. 

\example
Let $F$ be a number field. Then $h^0(F)$ itself is an Artin motive. 
It is self-dual.

\example 
A Dirichlet character is  a homomorphism 
\[
\chi:(\Z/N\Z)^{\mal}\to \C^{\mal}.
\]
The conductor $f$ is the smallest number such that $\chi$ factors 
through $(\Z/f\Z)^{\mal}$. The character is {\em primitive} if $f=N$.
It is convenient to extend $\chi$ to $\Z/N\Z$ by $\chi(a)=0$ if 
$a\in \Z/N\Z\ohne (\Z/N\Z)^{\mal}$. Consider the isomorphism
\[
rec:\Gal(\Q(\mu_N)/\Q)\xrightarrow{\isom} (\Z/N\Z)^{\mal}
\]
which maps the geometric Frobenius $\Fr_p$ at $p\nmid N$ to $p$. Via this isomorphism
we view $\chi$ as character
\[
\chi:\Gal(\Q(\mu_N)/\Q)\to \C^{\mal}.
\]
Let $E$ be the number field generated by all values of Dirichlet characters
of $(\Z/N\Z)^{\mal}$. Explicitly, it is the field of $\Phi(N)$-th roots
of unity, where $\Phi(N)=\#(\Z/N\Z)^{\mal}$. 
The character $\chi$ defines a rank one
Artin motive $V(\chi)$ with coefficients in $E$ as follows:
 $V(\chi)$ is the image of the projector $p_{\chi^{-1}}$ ({\em sic!})
on $h^0(\Q(\mu_N))$ 
\begin{align*}
p_{\chi^{-1}}:h^0(\Q(\mu_N))&\to h^0(\Q(\mu_N)) \\
\alpha&\mapsto\frac{1}{\Phi(N)} \sum_{\sigma\in G}\chi(\sigma)\sigma^* \alpha.\nonumber
\end{align*}
where we use the identification $G:=\Gal(\Q(\mu_N)/\Q)\isom (\Z/N\Z)^{\mal}$.
We call these {\em Dirichlet motives}.
In fact,
\[ h^0(\Q(\mu_N))=\bigoplus_{\chi}V(\chi) \]
where the sum is over all characters of conductor dividing $N$. If $F$ is
an abelian number field, its motive is a direct summand of the motive of
$h^0(\Q(\mu_N))$ for some $N$ and hence (after extension of coefficients
to big enough $E$) isomorphic to a direct sum of Dirichlet motives.

Let $V$ be an Artin motive. $V$ has {\em $p$-adic realizations} $V_p$ for
all primes $p$. They are  $E_p=E\tensor \Q_p$-modules with a
continuous operation of $G_\Q$. 

\example
 The $p$-adic realization of $V(\chi)$ is defined as
\[ V_p(\chi)=p_{\chi^{-1}}H^0_{et}(\spec F\tensor \Qbar,\Q_p)\]
 where the
projector is taken with respect to the action of the Galois group of $F$ 
over $\Q$. It is a $G_\Q$-module via the action on $\Qbar$ over $\Q$.
This also factors through $G(F/\Q)$ but gives the contragredient
representation. Hence
$V_p(\chi)$ is the
rank one $E_p$-module with operation of $G_\Q$ via $\chi$.

Let $V_B$ be the {\em Betti-realization} and $V_{\DR}$ the {\em de Rham realization}
of $V$. They are finite dimensional
$E$-vector spaces linked by an $E\tensor_\Q\C$-linear isomorphism
$V_B\tensor_\Q\C\isom V_{\DR}\tensor_\Q\C$. 

\example For $V=h^0(F)$ with $E=\Q$, the Betti-realization $V_B$ is given by the maps
from the set of embeddings $F\to \C$ to $\Q$. The operation of $G(F/\Q)$
is via $(g f)(\sigma)=f(\sigma g)$. 
It has a natural self-dual lattice $T_B$ given by the functions with values in
$\Z$. A natural basis is given by the functions $\delta_\tau$ with
$\delta_\tau(\tau)=1$ and $\delta_\tau(\sigma)=0$ for $\sigma\neq \tau$.
The de Rham realization is
$V_\DR=F\tensor_\Q E$.

{\em Motivic cohomology} of a number field
is given by Adams eigenspaces of $K$-theory of its ring of integers. Taking
direct summands this also defines motivic cohomology with coefficients
in $V$. The only non-vanishing motivic cohomology groups are
$\hm^1(\spec \Z,V(r))$ for $r\geq 1$ and $\hm^0(\spec\Z,V)$. 
For $r\geq 1$ there is the {\em Beilinson regulator map}
\[
r_\infty:\hm^1(\spec\Z,V(r))\to H^1_\Dh(\R,V_\R(r))
\]
with values in real Deligne cohomology. We use the identification
\[ H^1_\Dh(\R,V_\R(r))\isom V_{\DR}(r)_\R/V_B(r)^+_\R\isom \left(V_B(r-1)_{\R}\right)^+ \]
where $+$ denotes the invariants under complex conjugation and 
where the second isomorphism is induced by $\C=\R\oplus \R(-1)$.
$r_\infty\tensor\R$ is an isomorphism for $r>1$ by \cite{Borel} and \cite{Ra}.
For $r=0$ the cycle class map to singular cohomology induces
\[ z:\hm^0(\spec\Z,V)\to \left(V_B\tensor\R\right)^+=H^1_\Dh(\R,V_\R(1))\ .\]

Let $S$ be a finite set of rational primes and $j:\spec\Q\to \spec\Z[1/Sp]$
the natural inclusion.
We abbreviate
$H^i(\Z[1/Sp],V_p(r)):=H^i_{et}(\spec\Z[1/Sp],j_*V_p(r))$
where $j_*$ is the direct image in the {\'e}tale topology (sic, {\em not}
the derived direct image). Note that for $V=h^0(F)$, we have the
equality $H^i(\Z[1/Sp],V_p(r))=H^i(\Oh_F[1/Sp],\Q_p(r))$ as
the direct image of $\Q_p$ under the map $\spec \Oh_F[1/Sp]\to \spec
\Z[1/Sp]$ is $j_*V_p$.

For all primes $p$, there are also {\em $p$-adic regulator maps}
\[
 r_p:\hm^1(\spec\Z,V(r))\to H^1(\Z[1/Sp],V_p(r))
\]
induced by
\[
r_p:\hm^1(\Oh_F,\Q(r))\to H^1(\Oh_F[1/Sp],\Q_p(r)).
\]
The map $r_p\tensor\Q_p$ is again an isomorphism for $r> 1$ see \cite{Soulepadic} Theorem 1. We will need a refined
version of this regulator. 
Let $v\neq p$ be a rational prime. We denote by $I_v$ its inertia group
and by $\F_v$ the residue class field. One defines 
\[ H^1_f(\Q_v,V_p(r))=H^1(\F_v,V_p(r)^{I_v}) \ .\]
For $v=p$ we put
\[ D_{\cris}(V_p(r))=\left(B_{\cris}\tensor V_p(r)\right)^{G_{\Q_p}}\ ,\
\tan(V_p(r))=\left((B_{\DR}/Fil^0B_{\DR})\tensor V_p(r)\right)^{G_{\Q_p}}\ .
\]
The latter vanishes for $r\leq 0$ and is naturally isomorphic to $V_{\DR}(r)\tensor\Q_p$ for $r\geq 1$.

We consider
the fundamental short exact sequence of Bloch-Kato and Fontaine
\[ 0\to \Q_p\to B_{\cris}\xrightarrow{(\phi-1,\pi)} B_{\cris}\oplus B_{\DR}/Fil^0 B_{\DR}\to 0, \]
where $\phi$ is the Frobenius on $ B_{\cris}$ as in \cite{FontBou} 3.2
p. 210 footnote and $\pi$ is the 
composition $B_{\cris}\to B_{\DR}\to B_{\DR}/Fil^0 B_{\DR}$. 
This sequence induces after tensoring with $V_p(r)$ and taking cohomology 
\begin{multline*}
0\to H^0(\Q_p,V_p(r))\to  D_{\cris}(V_p(r))\to
 D_{\cris}(V_p(r))\oplus \tan(V_p(r))\to H^1(\Q_p,V_p(r))
\end{multline*}
The space $H^1_f(\Q_p,V_p(r))$ is defined as the image of the last map. Under local 
duality, it is dual to $H^1_{/f}(\Q_p,V_p^{\lor}(1-r))$ where
$H^1_{/f}=H^1/H^1_f$. Moreover, the $p$-adic exponential 
\[\exp_p:\tan(V_p(r))\to H^1(\Q_p,V_p(r))\]
is also induced from the same sequence. By definition it induces an isomorphism
of $\tan(V_p(r))$ with $H^1_f(\Q_p,V_p(r))$ if $r\neq 0$.
The local cohomology groups $H^1_f(\Q_p,V_p(r))$ are used to
define
\[ H^1_f(\Z[1/Sp],V_p(r))=\left\{ a\in H^1(\Z[1/Sp],V_p(r))\mid
                     a\in H^1_f(\Q_v,V_p(r)) \text{ for all $v\in Sp$}\right\} \]
Now, the $p$-adic regulator factors through
\[ r_p:\hm^1(\Z,V(r))\to H^1_f(\Z[1/Sp],V_p(r)) \]
and $r_p\tensor\Q_p$ is in fact always an isomorphism. (This follows
from the previous isomorphism for $r>1$, compatibility of $H^1_f$ with local duality \cite{Bloch-Kato} 3.8 and the explicit shape of $\exp_p$ for $r=1$, \cite{Bloch-Kato} p.358).

\subsection{The Bloch-Kato conjecture for Artin motives}\label{formulation}
Attached to an Artin motive $V$, there is an  $E\tensor \C$-valued
$L$-function
\[ L(V,s)=\prod_{v}\frac{1}{P_v(V, v^{-s})} \]
with $v$ a rational prime and
\[ P_v(V,t)=\det_{E_l}(1- \Fr_v t\mid V_l^{I_v} )\]
the characteristic polynomial of the geometric Frobenius $\Fr_v$ operating on 
the $l$-adic realization  ($I_v$ is the inertia group at $v$)
for any auxiliary prime $l\neq v$. It is in fact an element of the
polynomial ring $E[t]$.
For $V=h^0(F)$ and $E=\Q$ this gives the Dedekind-$\zeta$-function of $F$. 
Let $V(\chi)$ be a Dirichlet motive with coefficients in $E$ big enough. 
Let $N$ be the conductor of $\chi$ and consider $\chi$ 
 as a map
$\chi:\Gal(\Q(\mu_N)/\Q)\to E^*$. We identify 
\[ E\tensor\C\isom\bigoplus_{\tau:E\to \C}\C .\]
Then 
\[ L(V(\chi),s)= \big( L(\tau\circ\chi,s) \big)_\tau \]
where 
\[L(\tau\circ\chi,s) =
\sum_{n\geq 1} \frac{\tau\circ\chi(n)}{n^s}\ .
\]

\begin{defn}
 We write
$ L(V,r)^*\in (E\tensor \C)^*$ for the leading coefficient of the Taylor expansion of the $E\tensor\C$-valued function
$L(V,s)$  at $r$.
\end{defn}

\begin{defn}[\cite{FontBou} {\S} 1]\label{fundline}
Let $V$ be an Artin motive. The {\em fundamental line} $\Delta_f(V(r))$ is
the one-dimensional $E$-vector space
\begin{align*}
&\det_E \hm^0(\Z,V(r))\tensor \det_E^{-1}\hm^1(\Z,V^{\lor}(1-r))\tensor
\det_E^{-1}V_B(r)^+ &\text{if $r\leq 0$}\\
&\det_E \hm^0(\Z,V^{\lor}(1-r))\tensor \det_E^{-1}\hm^1(\Z,V(r))\tensor
\det_E^{-1}V_B(r)^+\tensor \det_E V_{\DR}(r) &\text{if $r\geq 1$}
\end{align*} 
\end{defn}

\begin{prop}[\cite{FontBou} 6.10]\label{firstiso}
There is a natural isomorphism
\[ \Delta_f(V(r))\otimes_{\Q}\R\isom E\tensor_\Q\R =:E_\infty\]
induced by the short exact sequence
\[
0\to \hm^0(\Z,V(r))_\R \to V_B(r)^+_\R\to \hm^1(\Z,V^{\lor}(1-r))^{\lor}_\R \to 0 
\]
in the case $r\leq 0$
and by the short exact sequence
\[
0\to  \hm^1(\Z,V(r))_\R \to V_B(r-1)^+_\R\to \hm^0(\Z,V^{\lor}(1-r))^\lor_\R\to 0
\]
together with the isomorphism $V_B(r-1)_\R^+\isom V_{\DR}\tensor\R/V_B(r)_\R^+$
in the case $r\geq 1$.
\end{prop}
\bew The maps are the cycle maps $z$ and the Beilinson regulator $r_\infty$
respectively their duals.
\bewende

Let $S$ be a finite set of primes. We are going
to define an isomorphism of $E_p:=E\otimes_{\Q}\Q_p$-modules
\[ \Delta_f(V(r))\tensor_{\Q}\Q_p\isom \det_{E_p}R\Gamma_c(\Z[1/Sp],V_p(r))\ , \]
where $R\Gamma_c(\Z[1/Sp],V_p(r))$ is defined in \cite{FontBou}, i.e., it
sits in an exact triangle
\[
R\Gamma_c(\Z[1/Sp],V_p(r))\to R\Gamma(\Z[1/Sp],V_p(r))\to V_p(r)^+\oplus
\bigoplus_{v|pS} R\Gamma(\Q_v,V_p(r)).
\]

As in \cite{FontBou} {\S} 4, we define
 $R\Gamma_{f}(\Q_v,V_p(r))\subset\tau_{\leq 1} R\Gamma(\Q_v,V_p(r))$ to be the
subcomplex with $H^0_f(\Q_v,V_p(r))=H^0(\Q_v,V_p(r))$ and
$H^1_f(\Q_v,V_p(r))$ as in section \ref{artinmotives}.
If $v\neq p$, we have a quasi isomorphism 
\[ R\Gamma_f(\Q_v,V_p(r))\isom[V_p(r)^{I_v}\xrightarrow{\Fr_v-1}V_p(r)^{I_v}]. \]
If $v=p$, recall that $D_{\cris}(V_p(r))=\left(B_{\cris}\tensor V_p(r)\right)^{G_{\Q_p}}$.
On this space  the Frobenius $\phi$ of $B_{\cris}$ still acts. One has a 
quasi isomorphism
\[ 
R\Gamma_f(\Q_p,V_p(r))\isom [D_{\cris}(V_p(r))\xrightarrow{(\phi-1,\pi)}
D_{\cris}(V_p(r))\oplus \tan({V_p(r)})],
\]
where $\pi$ is the canonical projection.
The above quasi isomorphisms induce by the theory of determinants
isomorphisms
\[
{\det}_{E_p}R\Gamma_f(\Q_v,V_p(r))\isom{\det}_{E_p}[V_p(r)^{I_v}\xrightarrow{\Fr_v-1}V_p(r)^{I_v}],
\]
for $v\neq p$ and
\[
{\det}_{E_p}R\Gamma_f(\Q_p,V_p(r))\isom{\det}_{E_p}[D_{\cris}(V_p(r))\xrightarrow{(\phi-1,\pi)}
D_{\cris}(V_p(r))\oplus \tan({V_p(r)})],
\]
for $v=p$. The determinants of the latter complexes are equal 
to ${\det}_{E_p}V_p(r)^{I_v}\otimes {\det}^{-1}_{E_p}V_p(r)^{I_v}$ and
 ${\det}_{E_p}D_{\cris}(V_p(r))\otimes {\det}^{-1}_{E_p}D_{\cris}(V_p(r))\otimes {\det}^{-1}_{E_p}\tan(V_p(r))$ respectively.
As these are tensor products of a rank one $E_p$-module and its dual,
these are identified with $E_p$ 
and ${\det}^{-1}_{E_p}\tan_{V_p(r)}$ respectively. We record this in a definition:

\begin{defn}\label{localident}
We identify
\[ \alpha:\det_{E_p}R\Gamma_f(\Q_v,V_p(r))\xrightarrow{\isom}\begin{cases}
E_p&\text{$v\neq p$ or $r\leq 0$}\\
\det_{E_p}^{-1}\tan(V_p(r))&\text{$v=p,r\geq 1$}
\end{cases}\]
where the isomorphisms are the ones  discussed above.
\end{defn}

 One has to be very careful in working with these isomorphisms.
Assume that $\phi-1$ is an isomorphism. This implies that the complex
 $[D_{\cris}(V_p(r))\xrightarrow{(\phi-1,\pi)}
D_{\cris}(V_p(r))\oplus \tan({V_p(r)})]$ is
quasi isomorphic to $\tan({V_p(r)})$. Hence the theory of determinants gives an
isomorphism
$\beta:[D_\cris(V_p(r))\xrightarrow{\phi-1}D_\cris(V_p(r))\oplus\tan({V_p(r)})]
\isom \tan({V_p(r)})$.
\begin{prop}\label{localfactor}
Under the above conditions
\[ \alpha=P_p(V,p^{-r}) \beta \]
where $P_p$ is as before the characteristic polynomial of $\Fr_p$.
\end{prop}
\bew As $\alpha$ and $\beta$ are the identity on $\det \tan(V_p(r))$ we may
as well assume it to be $E_p$.
Let $d_1,\dots,d_n$ be a basis of $D=D_\cris(V_p(r))$. Then
$d=d_1\land\dots \land d_n$ is a basis of $\det D$. By definition
$\alpha(d\tensor d^{-1})=1$
 and $\beta(d\tensor (\det(\phi-1) d)^{-1})=1$. This implies
$\beta=\det(\phi-1)^{-1}\alpha$. 
By the normalization of $\phi$ (\cite{FontBou} 3.2, footnote p. 210),
$\det(\phi-1)$ is the characteristic polynomial of $\Fr_p$ (up to sign).
\bewende

\rem
The reason for using
this normalization lies in the following fact: the formulation of 
Fontaine and Perrin-Riou is concerned with the full $L$-function. As we 
will see later, the determinant of the (cyclotomic) elements in 
Galois cohomology $H^1(\Z[1/pS], V_p(r))$
is related to the $L$-function with the Euler factors at $pS$ removed. The 
difference between these two $L$-functions are the Euler factors at $pS$,
which are exactly introduced by the above normalization of the determinant
of $R\Gamma_f(\Q_v,V_p(r))$.

Again as in \cite{FontBou} {\S} 4 we put:
\begin{align*}
R\Gamma_{/f}(\Q_v,V_p(r))&=\cone\left[ R\Gamma_f(\Q_v,V_p(r))
                                \to R\Gamma(\Q_v,V_p(r))\right] \\
R\Gamma_f(\Z[1/Sp],V_p(r))&=\cone\left[ R\Gamma(\Z[1/Sp],V_p(r))
\to R\Gamma_{/f}(\Q_p,V_p(r))\right][-1]
\end{align*}

\begin{lemma}[\cite{FontBou} 4.4]\label{firstlemma}
Let $V$ be an Artin motive.
Let $S$ be a finite set of primes.
\begin{enumerate}
\item 
For $r>1$, the complex $R\Gamma_f(\Z[1/Sp],V_p(r))$ is concentrated in
degree $1$ and its first cohomology is isomorphic to
$H^1(\Z[1/Sp],V_p(r))$. 
\item For $r=1$, the complex is concentrated in degrees $1$ and $3$. There
is a distinguished triangle
\[ H^1_f(\Z[1/Sp,V_p(1))[-1]\to R\Gamma_f(\Z[1/Sp],V_p(1))\to  H^0(\Z[1/Sp],V^{\lor}_p)^\lor[-3]\ .\]
\item
For $r<0$, the complex  $R\Gamma_f(\Z[1/Sp],V_p(r))$ is concentrated in
degree $2$ and its second cohomology is isomorphic to
$H^1(\Z[1/Sp],V^{\lor}_p(1-r))^{\lor}$.
\item For $r=0$, the complex $R\Gamma_f(\Z[1/Sp],V_p(r))$ is concentrated in
degrees $0$ and $2$. There is a distinguished triangle
\[  H^0(\Z[1/Sp],V_p) \to R\Gamma_f(\Z[1/Sp],V_p)\to H^1_f(\Z[1/Sp],V^{\lor}_p(1))^{\lor}[-2] \]   
\end{enumerate}
\end{lemma}
\bew See \cite{FontBou} p. 215-216. The maps are either inclusions or
their duals via the localization sequence. \bewende

Note finally, that there is a natural distinguished triangle (see \cite{FontBou} p. 215)
\[ R\Gamma_c(\Z[1/Sp],V_p(r)) \to R\Gamma_f(\Z[1/Sp],V_p(r))\to V_p(r)^+\oplus\bigoplus_{v\in Sp} R\Gamma_f(\Q_v,V_p(r)) \]
\begin{defn}\label{globalident}
We identify
\begin{align*} 
 \det_{E_p} R\Gamma_c&(\Z[1/Sp],V_p(r))\\
&\isom \det_{E_p} R\Gamma_f(\Z[1/Sp],V_p(r))\tensor \det^{-1}_{E_p}V_p(r)^+\tensor\det_{E_p}^{-1}\bigoplus_{v\in Sp}R\Gamma_f(\Q_v,V_p(r)) \\
&\isom \det_{E_p} R\Gamma_f(\Z[1/Sp],V_p(r))\tensor \det^{-1}_{E_p}V_p(r)^+
\tensor \det_{E_p} \tan({V_p(r)})\\
&\isom \Delta_f(V(r))\tensor_{\Q_p}E_p
\end{align*}
where the second isomorphism uses \ref{localident} and the last is by term by term
comparison via the $p$-adic regulators (and their duals) together with the 
isomorphism
of $\tan({V_p(r)})$ with $V_{\DR}(r)\tensor\Q_p$ in the case $r\geq 1$.
\end{defn}

We can now formulate the Tamagawa number conjecture of Bloch and Kato. The
precise formulation is taken from Fontaine's Bourbaki talk \cite{FontBou} who also shows equivalence
to the original statement. In his formulation a set of primes $S$ appears, but
see \ref{indep}  below.

\begin{conj}[Bloch-Kato]\label{conj}
Let $V$ be an Artin motive over $\Q$ with coefficients in a number field $E$
and $r\in \Z$.  Denote by $L(V,r)^*$ the leading coefficient of the Taylor expansion of the $L$-function
of $V$ at $r$. Let $\delta\in \Delta_f(V(r))\tensor_{\Q}\R$ be such that 
$L(V,r)^*\delta$ is mapped to $1$ under
the isomorphism of \ref{firstiso}
\[ \Delta_f(V(r))\tensor_{\Q}\R\isom E\tensor_{\Q}\R \ .\]
Then, $\delta\in \Delta_f(V(r))$ and 
for all primes $p$, the image of $\delta$  under the
identification of \ref{globalident} 
\[ \Delta_f(V(r))\tensor\Q_p\isom\det_{E_p}R\Gamma_c(\Z[1/p],V_p(r))
\]
is a generator of
$\det_{\Oh_p}R\Gamma_c(\Z[1/p],T_p(r))$ where $\Oh_p=\Oh\otimes\Z_p$ and $T_p$ is any $G_\Q$-stable 
$\Oh_p$-lattice in $V_p$.
\end{conj}
\bem The element $\delta\in \Delta_f(V(r))$ is uniquely determined up to units
in $\Oh=\Oh_E$ by the integrality condition. Hence the above can be
read as a formula for $L(V,r)^*$. We are going to make the explicit
translation at least in a special case in section \ref{lichtenbaum}.

\begin{prop}\label{indep}
\begin{description}
\item[a)] The conjecture is well-posed, i.e., independent of the choice of lattice $T_p$.
\item[b)] It is equivalent to the conjecture formulated by Fontaine in \cite{FontBou} and hence to the original conjecture in \cite{Bloch-Kato}.
\end{description}
\end{prop}
\bew
a) Let $T_p'\subset T_p$ be
a sublattice. The quotient is a finite group, say $Q$.
The assertion now follows from the computation of Galois-cohomology
with torsion coefficients and the fact that the Euler characteristic
of $R\Gamma_c( \Z[1/p], Q)$ is $1$, see e.g. \cite{Milne} I{\S} 5 Theorem 5.1.

b) Fontaine's formulation of the conjecture coincides with the 
above but he replaces 
$R\Gamma_c(\Z[1/p],T_p(r))$ by $R\Gamma_c(\Z[1/Sp],T_p(r))$ where
$S$ is any set of primes including the primes of bad reduction of $V$.
What we have to show is independence of $S$ (including all primes of bad
reduction or not.)
Hence let $v$ be a prime different from $p$ which is not in $S$. 
We only have to remark that
there is a localization sequence
\[ R\Gamma_c(\Z[1/Svp],T_p(r))\to R\Gamma_c(\Z[1/Sp],T_p(r))\to
  R\Gamma(\F_v,T_p(r)^{I_v}) \]
whether $v$ is a good reduction prime or not (see e.g. \cite{Milne} II prop. 2.3 d)). Then
\[ R\Gamma(\F_v,V_p(r)^{I_v})=R\Gamma_f(\Q_v,V_p(r)) \]
and the isomorphisms in \ref{globalident} are compatible with enlarging $S$.

The equivalence with original conjecture of Bloch and Kato was shown
by Fontaine in \cite{FontBou}.
\bewende

It is often useful to reformulate the $p$-adic condition in terms of global
duality, which gives for $p\neq 2$
\[
\det R\Gamma_c(\Z[1/p],V_p(r))\tensor \det V_p(r)^+\isom
\det R\Gamma(\Z[1/p],V_p^{\lor}(1-r)).
\]
\begin{prop}\label{ohnec}
Let $T_p$ be a Galois stable lattice in $V_p$ and $p\neq 2$. Then the integral structures
on the left resp. right hand side of the last formula given by
\[
 \det R\Gamma_c(\Z[1/p],T_p(r))\tensor\det T_p(r)^+\mbox{ resp. } 
\det  R\Gamma(\Z[1/p],T_p^{\lor}(1-r))
\]
agree.
\end{prop}
\rem For $p=2$ the complex $R\Gamma(\Z[1/p],T_p^{\lor}(1-r))$ is not perfect and one
has to consider $\tau_{\leq 2}R\Gamma(\Z[1/p],T_p^{\lor}(1-r))$ instead. Moreover $T_p(r)^+=H^0(\R, T_p(r))$
has to be substituted by the zeroeth Tate-cohomology of $T_p(r)$.\\
 
\bew Let $^{*}$ denote the Pontryagin dual $\Hom(\cdot,E_p/\Oh_p)$.
Global duality implies
\[ \det R\Gamma_c(\Z[1/p],T_p(r))\tensor\det T_p(r)^+\isom \det^{-1}R\Gamma(\Z[1/p],T_p^{*}(1-r))^{*} \]
under the above rational isomorphism. Hence we have to show that the integral 
structures $\det  R\Gamma(\Z[1/p],T_p^{\lor}(1-r))$ and 
$\det^{-1} R\Gamma(\Z[1/p],T_p^{*}(1-r))^{*}$ in the $E_p$-module 
$\det R\Gamma(\Z[1/p],V_p^{\lor}(1-r))$ coincide.
For a bounded complex $A$, we have
\begin{multline*} 
R\Hom(A,\Oh_p)=R\Hom(A,R\Hom(E_p/\Oh_p,E_p/\Oh_p))\\
  =R\Hom(A\tensor^{\Bbb{L}}E_p/\Oh_p,E_p/\Oh_p)=(A\tensor^{\Bbb{L}}E_p/\Oh_p)^* 
\end{multline*}
Moreover,
\begin{multline*}
R\Gamma(\Z[1/p],T_p^{*}(1-r))= \\
R\Gamma(\Z[1/p],T_p^{\lor}(1-r)\tensor^{\Lbb}_{\Oh_p}
E_p/\Oh_p)=R\Gamma(\Z[1/p],T_p^{\lor}(1-r))\tensor^{\Lbb}_{\Oh_p}
E_p/\Oh_p
\end{multline*}
because $R\Gamma(\Z[1/p],T_p^{\lor}(1-r)\tensor^{\Lbb}_{\Oh_p} E_p)=R\Gamma(\Z[1/p],T_p^{\lor}(1-r))\tensor^{\Lbb}_{\Oh_p} E_p$. Together
 this implies
\begin{align*} 
\det^{-1}R\Gamma(&\Z[1/p],T_p^{*}(1-r))^{*}
=
\det^{-1}\left(R\Gamma(\Z[1/p],T_p^{\lor}(1-r))
\tensor^{\Lbb}_{\Oh_p} E_p/\Oh_p\right)^{*}\\
&=\det^{-1}R\Hom(R\Gamma(\Z[1/p],T_p^{\lor}(1-r)),\Oh_p)
=\det R\Gamma(\Z[1/p],T_p^{\lor}(1-r))
\end{align*}
\bewende

\rem
Kato formulates versions of the Tamagawa number conjecture for values at
positive integers in \cite{Kato Main conj} and for values at negative
integers in \cite{galaxy}. 
The above proposition implies that his conjectures
are equivalent to the ones formulated by Fontaine. Note that in
the formulation of the Tamagawa number conjecture for the full
$L$-function the integral structure of the fundamental line defined
using $\det R\Gamma(\Z[1/p],T_p^\lor(1-r))$ still involves the
identification \ref{localident} of the local factor.

\subsection{Main Theorem}
We now state our main theorem on the Bloch-Kato conjecture and sum up the
known cases.

\begin{thm}\label{maintheorem}
The Bloch-Kato conjecture \ref{conj} holds up to powers of $2$ for all $r\in \Z$ and  all abelian Artin motives,
i.e.,
for Artin motives which are direct sums of direct summands of $h^0(F)$'s where $F$ is an
abelian number field. In particular, it holds for all Dirichlet motives.
\end{thm}

\begin{enumerate}
\item
Some parts of the proof certainly work for the prime $2$ but we did not check whether there are serious problems or not.
\item The case $V=h^0(\Q)$ already appears in the original paper
by Bloch and Kato \cite{Bloch-Kato} Theorem 6.1 (conjecture 6.2 has
meanwhile been proved in \cite{Hu-Wi}.)
 In \cite{Kato Main conj} {\S} 6 Kato proves the $p$-part of the equivariant
conjecture for $\Q(\mu_{p^n})^+$ for $r\geq 2$ even.
\item The Beilinson conjecture for Dirichlet motives, i.e.
the first part of conjecture \ref{conj},  was already proved 
by Beilinson in \cite{Be} (with some corrections in \cite{Ne} and \cite{Es}). 
\item
Fontaine \cite{FontBou} {\S} 10 mentions the 
case of Dirichlet motives $V(\chi)$ but assumes that $p^2$ resp. $p$ does divide $N$ or $\phi(N)$ (where $N$ is the conductor of $\chi$) depending on parity conditions.  
\item
The Bloch-Kato conjecture for $V=h^0(F)$ with coefficients in $\Q$ and
$r\leq 0$ is easily seen to be equivalent to the cohomological Lichtenbaum
conjecture (see \ref{lichtenbaumconj} below).  The case $V=h^0(F)$ where $F$ is 
totally real and $r\leq 0$ is odd  
 is a direct consequence of the main conjecture
due to Wiles. In \cite{Kolster} the cohomological Lichtenbaum conjecture is treated
for abelian number fields. Together with
 the
correction of Euler factors worked out by Kolster and Nguyen Quang Do in
\cite{KoNg}, they prove the cohomological Lichtenbaum conjecture for abelian number fields up
to an explicit set of bad primes. 
It follows (for all primes $p\neq 2$) as a corollary of our main theorem \ref{maintheorem}.
\item In recent work \cite{BeNg}, Benois and Nguyen show how to reduce the Bloch-Kato
conjecture for $h^0(F)$ and $r\geq 1$ to the case of $r\leq 0$ by showing
compatibility under the functional equation. See also section \ref{fcteq} below.
\item Recall that the Beilinson conjecture (which
determines the $L$-value up to rational factor) is known for $h^0(F)$ where
$F$ is any number field. However, it is not known for all Artin motives!

\item
There is also an equivariant version of the Bloch-Kato conjecture,
formulated by Kato in the abelian case and by Burns-Flach in general (see also
conjecture \ref{equivconj} below).
A proof of the equivariant conjecture
for Tate motives $\Q(r)$ and $r\leq 0$
 over abelian number fields and $p\neq 2$  is given by Burns and 
Greither in \cite{BuGr}. 
Our theorem \ref{maintheorem} is equivalent to the equivariant conjecture
with respect to the maximal order $\widetilde{\Z[G]}$ in $\Q[G]$
and hence for $r\leq 0$ a consequence of the result of Burns and Greither.  
Using a key observation by Burns and Greither,
it is also possible to deduce the  equivariant case from  
theorem \ref{maintheorem}.
\end{enumerate}

\subsection{Relation to Lichtenbaum conjecture}\label{lichtenbaum}
The above theorem \ref{maintheorem} can be restated to give a formula for the $L$-value.
For simplicity of notation, we do this in the case $E=\Q$, $V=h^0(F)$ where
$F$ is a number field and for $r\leq -1$.

\begin{thm}[Cohomological Lichtenbaum conjecture]\label{lichtenbaumconj}
Let $F$ be a number field, $k\geq 2$.  
Let $R_k(F)$ be the Beilinson regulator of $F$, i.e., the covolume of
the image of $K_{2k-1}(F)$ in $h^0(F)^{\lor}(k-1)_\R^+$
under the Beilinson regulator map.
Then the Bloch-Kato conjecture for $h^0(F)$ and $1-k$ implies that
up to sign and up to powers of $2$,
\[ \zeta_F(1-k)^*=R_k(F)\prod_{p}\frac{\# H^2(\Oh_F[1/p],\Z_p(k))}
                                  {\# H^1(\Oh_F[1/p],\Z_p(k))_{\mathrm{tors}}} \ .\]
In particular, the cohomological Lichtenbaum conjecture holds for abelian
number fields $F$.
\end{thm}

\bew Let us show that conjecture  \ref{conj} with $r=1-k\leq -1$ 
implies this formula for the $L$-value. 
Let 
$T_B^{\lor}(-r)^+$ be the natural integral structure 
in $h^0(F)^{\lor}(-r)_\R^+$ and 
$\Omega\subset\hm^1(F,\Q(1-r))$ the quotient of $K_{2k-1}(F)$ by the 
torsion subgroup.
In this case 
\[ \Delta_f(V(r))=\det^{-1}\hm^1(F,\Q(1-r))\tensor\det^{-1}
\left(h^0(F)_B(r)\right)^+.\] 
In the fundamental line one has the 
lattice $\det^{-1}\Omega\tensor \det^{-1} T_B(r)^+$. By definition 
\[ \det\Omega=R_{1-r} \det T_B^{\lor}(-r)^+=R_{1-r}\det^{-1}T_B(r)^+ \]
where $R_{1-r}$ is the Beilinson regulator of $F$.
Recall that $\delta$ is the element of the fundamental line corresponding
to $1/\zeta_F(r)^*$ under its identification with $\R$. In terms of the above lattice
\[ \delta\Z=\frac{R_{1-r}}{\zeta_F(r)^*} \det^{-1}\Omega\tensor \det^{-1}T_B(r)^+ \]
Then conjecture \ref{conj} implies that $\delta$ 
is mapped to  a generator
of $\det R\Gamma_c(\Z[1/p],T_p(r))$ under the isomorphism
$\Delta_f(V(r))\tensor \Q_p\isom \det R\Gamma_c(\Z[1/p],V_p(r))$. 
In our case 
\[
 R\Gamma_f(\Z[1/p],V_p(r))\isom H^1(\Z[1/p],V^{\lor}_p(1-r))^{\lor}[-2]=
R\Gamma(\Z[1/p],V^{\lor}_p(1-r))^{\lor}[-3] 
\]
The isomorphism of the fundamental line was defined by
\begin{multline*}
\Delta_f(V(r))\tensor \Q_p\xrightarrow{\det r_p^\lor\tensor \id} 
\det H^1(\Z[1/p],V^{\lor}_p(1-r))^{\lor}\otimes \det^{-1} V_p(r)^+\\
\xrightarrow{\id\tensor\alpha^{\lor}} 
\det H^1(\Z[1/p],V^{\lor}_p(1-r))^{\lor}\otimes \det^{-1} V_p(r)^+
\tensor \det^{-1} R\Gamma_f(\Q_p,V_p(r))
\end{multline*}
with $\alpha$ as in \ref{localident}. By lemma \ref{localfactor} this
isomorphism differs from the one induced by $\det r_p^\lor\tensor \id$ alone by
the Euler factor $P_p(V,p^{-r})=\prod_{v\mid p}(1-N(v)^{k-1})$.
This is a $p$-unit as $k\geq 2$.

Using proposition \ref{ohnec} this implies that $r_p(\delta)$ is
mapped to
a generator of 
\[ \det R\Gamma(\Z[1/p],T_p^{\lor}(1-r)) \tensor\det^{-1}T_p(r)^+\]
or equivalently,
\begin{align*}
\frac{R_{1-r}}{\zeta_F(r)^*} \Z_p&=\det r_p(\Omega\otimes\Z_p)\tensor
\det R\Gamma(\Z[1/p],T_p^{\lor}(1-r))\\
&=\det r_p(\Omega\otimes\Z_p) \tensor\det^{-1}H^1(\Z[1/p],T_p^{\lor}(1-r))
\tensor \det H^2(\Z[1/p],T_p^{\lor}(1-r))\\
&=\det^{-1}(H^1(\Z[1/p],T_p^{\lor}(1-r))/r_p(\Omega\tensor\Z_p))
\tensor \det H^2(\Z[1/p],T_p^{\lor}(1-r))
\end{align*}
Write again $k=1-r$.
 We have
$H^i(\Z[1/p],T_p^\lor(k))=H^i(\Oh_F[1/p],\Z_p(k))$.
Soul{\'e} has proved that the map $K_{2k-1}(F)\otimes\Z_p\to H^1(\Oh_F[1/p],\Z_p(r))$
is surjective and the ranks agree by the computation of Borel. Hence
$r_p(\Omega\tensor\Z_p)$ is identified with the free part of
$H^1(\Oh_F[1/p],\Z_p(k))$, i.e., the quotient is
$H^1(\Oh_F[1/p],\Z_p(k))_{\tors}$.
 Note finally that a finite $\Z_p$-module of order $h$ has determinant $h^{-1}\Z_p$.
\bewende

Many cases were known before, see the discussion after theorem \ref{maintheorem}.
Lichtenbaum's original conjecture \cite{LC} involves the Borel regulator rather
than the Beilinson regulator and $K$-groups rather than Galois cohomology.
The missing ingredient is the
Quillen-Lichtenbaum conjecture
\[ K_{2r-i}(F)\tensor\Z_p\isom H^i(\Oh_F[1/p],\Z_p(r)) \]
for $i=1,2$ and $r\geq 2$.

\subsection{The equivariant conjecture}
There is also an equivariant version of the Tamagawa number conjecture.
It was first introduced by Kato for abelian groups in \cite{Kato Main conj}.
The case of non-abelian groups is treated by Burns and Flach in
\cite{BuFl}. We are going to treat the abelian case. 

As before let
$M$ be an Artin motive over $\Q$ with coefficients in $E$. Let
$K$ be an abelian extension of $\Q$ with Galois group $G=\Gal(K/\Q)$.
Let $X(K/\Q)$ be the set of $\C$-valued characters of $G$.
We assume for simplicity that $E$ contains all values of elements
of $G$. Note that we have
\[ E[G] =\bigoplus_{\omega\in X(K/\Q)} E(\omega) \]
where $G$ operates on $E(\omega)$ via $\omega^{-1}$. For every Artin motive $M$ let 
$M(\omega):=M\otimes V(\omega)$.
\begin{defn}\label{Lequidefn}
The equivariant $L$-value
\[ L(K/\Q,M,r)^*\in E_\infty[G]^* \]
has $\omega$-component $L(M(\omega),r)^*$. 
The equivariant fundamental line is defined as
\[ \Delta_f(K/\Q,M(r))=\bigoplus_{\omega\in X(K/\Q)} \Delta_f(M(\omega)(r)) \]
\end{defn}

\rem 
\begin{enumerate}\item $L(K/ \Q,M,r)^*$ is an element of $E_\infty[G]=E\tensor\R[G]$ rather
than $E\tensor\C[G]$ as the element is invariant under
complex conjugation on coefficients.
\item Let $M$ be an Artin motive with coefficients in $E$. Then 
$L(K/ \Q,M,r)^*$ is indeed an element in $E_\infty[G]$, even though
the $\omega$-components have coefficients in some bigger field
$E'$. This is seen by checking invariance of $L(K/ \Q,M,r)^*$
under $\Gal(E'/E)$.
\item If $K'/K$ is an extension such that $K'/ \Q$ is still
abelian, then $L(K'/\Q,M,r)^*$ is mapped to $L(K/ \Q,M,r)$ under
the norm map. 
\end{enumerate}

The absolute comparison morphisms for all $\omega$ together give 
natural isomorphisms
\begin{align*}
 \Delta_f(K/\Q,M(r))\tensor \R &\isom E_\infty[G] \\
\Delta_f(K/\Q,M(r))\tensor \Q_p &\isom \det_{E_p[G]}R\Gamma_c(\Z[1/p],\bigoplus_\omega M_p(\omega)(r))
\end{align*}
The isormophism $\bigoplus_\omega M_p(\omega)\isom M_p\tensor E_p[G]$ induces
\[ R\Gamma_c(\Z[1/p],\bigoplus_\omega M_p(\omega)(r))=R\Gamma_c(\Oh_K[1/p],M_p(r)) \]
Now we can formulate the equivariant conjecture:

\begin{conj}[Kato, Burns-Flach]\label{equivconj}
Let $M$ be an Artin motive over $\Q$ with coefficients in $E$ and $r\in \Z$.
Let $K/\Q$ be an abelian extension with Galois group $G$. Let 
$\delta(K/\Q)\in \Delta_f(K/\Q,M(r))\otimes\R$
be such that $L(K/\Q,M,r)^*\delta(K/\Q)$ is mapped to $1$ under the
isomorphism with $E_\infty[G]$. Then $\delta (K/\Q)\in \Delta_f(K/\Q,M(r))$,
 and for all primes the image of $\delta(K/ \Q)$ under 
the isomorphism with $\det_{E_p[G]}R\Gamma_c(\Oh_K[1/p],M_p(r))$ is a generator
of $\det_{\Oh_p[G]}R\Gamma_c(\Oh_K[1/p],T_p(r))$ where $T_p$ is any
$G_\Q$-stable lattice of $M_p$.
\end{conj}

Clearly, this is identical with conjecture \ref{conj} if $K=\Q$. 
Note
that the $\omega$-component of $\delta(K/\Q)$ is nothing but
$\delta$ for the motive $M(\omega)(r)$. Also the $\omega$-component of
the image of $\delta(K/\Q)$ in
$\det_{E_p[G]}R\Gamma(\Oh_K[1/p],M_p(r))$ is given by the image
of $\delta$ (for $M(\omega)(r))$ in 
$\det_{E_p}R\Gamma(\Z[1/p],M_p(\omega)(r))$.

\section{First steps in the proof}
\subsection{Overview}\label{secoverview}
We start with a few reductions.
As extension of coefficients is faithfully flat, it suffices to prove
the conjecture after extension of coefficients. The conjecture is
also compatible with direct sums of motives. As any abelian Artin motive
decomposes into a direct sum of Dirichlet motives after a suitable extension
of coefficients, it suffices to consider the case of Dirichlet motives $V(\chi)$. The proof is organized as follows:

\begin{enumerate}
\item We prove the Bloch-Kato conjecture directly in the case of the
full motive $h^0(F)$ where $F$ is a number field and $r=0$, see section
\ref{secclassnumber}. By the compatibility with the functional equation
  (which can be checked directly in this case) the conjecture also holds
for $h^0(F)$  and $r=1$.
\item We then use Euler system methods to establish a divisibility statement
for Iwasawa modules in the case $\chi(-1)=(-1)^{r}$, see \ref{strongdivisibility}.
By the class number trick --- with the class number formula replaced by the Bloch-Kato conjecture for $F=\Q(\mu_N)$ and $r=1$ --- we prove the
Main Conjecture of Iwasawa theory from it, see \ref{equality}.
\item Using Kato's explicit reciprocity law, we prove the Bloch-Kato conjecture
for $r\geq 1$ and $\chi(-1)=(-1)^{r}$ from the Main Conjecture, 
see \ref{Bloch-Katorge1}. 
The necessary computation was already used in the ``class number trick'' in
the previous step.
In fact, the argument proves the $p$-part of
the equivariant Bloch-Kato conjecture for the same $r$ and $\chi$ and the
cyclotomic $\Z/p^n$-extension over $\Q$. 

\item Using the precise understanding of the regulators of cyclotomic elements
we prove the Bloch-Kato conjecture
for $r\leq 0$ and $\chi(-1)=(-1)^{r}$ from the Main Conjecture, 
see \ref{Bloch-Katorle0}. 
However, 
the argument does not work for $r=0$ and $\chi(p)=1$, the case of ``trivial zeroes''.
The trivial zeroes will be treated in in the last step.
\item By the compatibility of the Bloch-Kato conjecture under the functional
equation, we deduce the Bloch-Kato conjecture in the remaining cases for $r\neq 0,1$, see \ref{Bloch-Katoung1} and \ref{Bloch-Katoung2}.  This
works even equivariantly and hence shows the $p$-part of the equivariant
conjecture for $r>1$, $\chi(-1)=(-1)^{r-1}$ and the cyclotomic $\Z/p^n$-extension of $\Q$.
\item From the equivariant  Bloch-Kato conjecture in the previous step,
we can deduce the second version of the Main conjecture, see \ref{mc2}.
Conversely, the new version allows to prove the Bloch-Kato conjecture
for $r=0$ and $r=1$ unless there are trivial zeroes.
\item The last exceptional case $r=0$ and $\chi(p)=1$, i.e., the case of 
trivial zeroes, follows again by
the functional equation see \ref{Bloch-Katofcteq}.
\end{enumerate}


\subsection{Functional equation}\label{fcteq}
We first study the compatibility of the equivariant Bloch-Kato conjecture under the
functional equation. The corresponding statements can already be found
in \cite{Motiveband} in the absolute case ($K=\Q$) and in \cite{BuFl} even
for the case of non commutative coefficients.
 
Let $K/ \Q$ be abelian and $M$ an Artin motive with coefficients in $E$. As before let $X(K/ \Q)$
be the set of $\C$-valued characters of $G:=\Gal(K/ \Q)$. For $\omega\in X(K/ \Q)$
define $M(\omega):=M\otimes V(\omega)$. Consider for each $M(\omega)$ and $r\geq 1$ the $E$-module of rank one
\[\Delta_{\loc}(M(\omega)(r))=\det^{-1}M_B(\omega)(r)^+\tensor\det M_{\DR}(\omega)(r)\tensor\det M_B(\omega)^{\lor}(1-r)^+ \]
After tensoring with $\R$ we identify it with $E\otimes\R$ 
using the isomorphism
\[ M_{\DR}(\omega)(r)_\R/M_B(\omega)(r)^+_\R \isom M_B(\omega)(r-1)^+_\R \]
Let $\epsilon(\omega)(r)$ be the element of $\Delta_{\loc}(M(\omega)(r))$
corresponding to $\frac{L(M(\omega)^{\lor},1-r)^*}{L(M(\omega),r)^*}$ under this isomorphism.
Define
\[
\Delta_{\loc}(K/ \Q,M(r)):=\bigoplus_{\omega\in X(K/ \Q)}\Delta_{\loc}(M(\omega)(r))
\]
(this is a module of rank $1$ over $E[G]=\bigoplus_{\omega}E(\omega)$)
and let $\epsilon (K/ \Q, M(r)):=\left(\epsilon(\omega)(r)\right)_{\omega}\in \Delta_{\loc}(K/ \Q,M(r))$.
On the other hand, each  $\Delta_{\loc}(M(\omega)(r))\tensor E_p$ is isomorphic to
\[ \det^{-1}R\Gamma(\Q_p,M_p(\omega)(r))\tensor\det^{-1} M_p(\omega)(r)^+\tensor\det M_p(\omega)^{\lor}(1-r)^+ \]
via the $p$-adic exponential and using \ref{localident}. Thus we get
\begin{align*}
\Delta_{\loc}(K/ \Q,M(r))\otimes E_p\isom&\det^{-1}_{E_p[G]}R\Gamma(\Q_p\otimes K,M_p(r))\tensor
\det^{-1}_{E_p[G]}\bigoplus_{\omega\in X(K/ \Q)} M_p(\omega)(r)^+\\&\tensor\det_{E_p[G]} \bigoplus_{\omega\in X(K/ \Q)}M_p(\omega)^{\lor}(1-r)^+,
\end{align*}
where we have used as before that $R\Gamma(\Q_p,\bigoplus_{\omega}M_p(\omega)(r))\isom R\Gamma(\Q_p\otimes K,M_p(r))$.
Observe that for any Galois-stable lattice $T_p$ in $M_p$ the $\Oh_p[G]$-module 
$T_p[G]:=T_p\otimes \Oh_p[G]$ is
a lattice in $\bigoplus_{\omega\in X(K/ \Q)} M(\omega)_p$. We also write $M_p[G]:=M_p\otimes E_p[G]$ and similarly for $M_B$ and $M_{\DR}$.
\begin{conj}\label{localconj}
Let $M$ be an Artin motive, $r\geq 1$, $p$ a prime.
Then the image of $\epsilon(K/ \Q, M(r))$ under the above isomorphism is a generator of
\[  \det_{\Oh_p[G]}^{-1}R\Gamma(\Q_p\otimes K,T_p(r))\tensor\det_{\Oh_p[G]}^{-1} T_p[G](r)^+\tensor\det_{\Oh_p[G]}( T_p[G])^{\lor}(1-r)^+\]
where $T_p$ is any Galois-stable lattice in $M_p$.
\end{conj}

\begin{thm}[Kato]\label{Katofcteq}
The above conjecture holds for $K/ \Q$ abelian, $M$ an abelian Artin motive and $r >1$, $p\neq 2$.
\end{thm}
\bew 
This was shown for the motive $\Q$ and $K=\Q$ by Bloch and Kato
in \cite{Bloch-Kato}. Perrin-Riou \cite{PR}  deduced it in the case that $p$ is
unramified in $M$. Benois and Nguyen Quang Do in \cite{BeNg} show it for
the full motive $h^0(F)$ of an abelian number field. Finally, there
is unpublished work of Kato \cite{epsilon} which settles the conjecture for
$K/ \Q$  abelian and $M=\Q$, thus for all abelian Artin motives. 
In appendix \ref{appfcteq}  we give a proof in the special case
for  $K=\Q(\mu_{p^n})$, $M$ arbitrary (abelian), $r>1$ and 
$p\neq 2$ (\ref{localconjequiv}). 
\bewende

\rem The case $r=1$, $K=\Q$ follows at
the very end of our arguments from the Bloch-Kato conjecture at $r=0$ and
$r=1$.

\begin{prop} \label{functionaleq}
Let $p\neq 2$, $K/ \Q$ abelian, $M$ be an Artin motive and $r\geq 1$. Then two of the following 
assertions imply the third.
\begin{enumerate}
\item
The $p$-part of the equivariant Bloch-Kato conjecture for $K/ \Q$, $M$ and $r$.
\item
The $p$-part of the equivariant Bloch-Kato conjecture for $K/ \Q$, $M^{\lor}$ and $1-r$.
\item Conjecture \ref{localconj} for $K/ \Q$, $M$, $r$ and $p$.
\end{enumerate}
\end{prop}
\bew 
We give the argument for the convenience of the reader. 
By definition
\begin{multline*} 
\Delta_f(K/ \Q,M(r))\tensor\Delta_f(K/ \Q,M^{\lor}(1-r))^{\lor}=\\
\det_{E[G]}^{-1}M_B[G](r)^+\tensor\det_{E[G]} M_{\DR}[G](r)\tensor\det_{E[G]}(M_B[G])^{\lor}(1-r)^+
\end{multline*}
The identification of the left hand side with $E[G]\otimes\R$ induced by \ref{firstiso} is
compatible with the above identification of the right hand side with $E[G]\otimes\R$.
The element $\delta(K/ \Q,M(r))\tensor\delta^{-1}(K/ \Q,M^{\lor},1-r)$ is nothing but $\epsilon(K/ \Q,M(r))$.

On the other hand, the tensor product of the fundamental lines
with $E_p$ in \ref{globalident} can be simplified using global duality (proposition  \ref{ohnec} over $K$ instead of
$\Q$)
and the localization sequence:
\begin{align*}
& \Delta_f(K/ \Q,M(r))\tensor\Delta_f(K/ \Q,M^{\lor}(1-r))^{\lor}\tensor E_p\\
&\isom
 \det_{\Oh[G]} R\Gamma_c(\Oh_K[1/p],M_p(r))\tensor \det_{\Oh[G]}^{-1} R\Gamma(\Oh_K[1/p],M_p(r))
\tensor
  \det_{\Oh[G]} (M_p[G]) ^{\lor}(1-r)^+\\
  &\isom \det_{\Oh[G]}^{-1} R\Gamma(\Q_p\otimes K,M_p(r))\tensor
  \det_{\Oh[G]}^{-1} M_p[G](r)^+\tensor
  \det_{\Oh[G]} (M_p[G])^{\lor}(1-r)^+.
\end{align*}
This is the second term appearing in the local conjecture. The integral
structures are compatible with this simplification by
proposition \ref{ohnec} (over $K$).
\bewende

\subsection{The class number formula case}\label{secclassnumber}
We check the first case of the Bloch-Kato conjecture by hand. This
is an old result, mentioned e.g. in \cite{FontBou} {\S} 8, but we did not
find a proof in the literature. As it is the crucial case to which all
the other cases are finally reduced, we discuss it in some detail.
\begin{prop}\label{classnumber}
The Bloch-Kato conjecture holds in the case $V=h^0(F)$, $r=0$ and $r=1$
up to powers of $2$.
\end{prop}
The rest of this section is devoted to the proof of this proposition.
The arguments are parallel to the ones used in deducing the Lichtenbaum
conjecture. 

\bew We start with the case $V=h^0(F)$ and $r=0$.
We write $\hm^i(\Oh_F,\Q(j))=\hm^i(\Z,V(j))$.
The fundamental line in this case is 
\[
\Delta_f(V):=\det \hm^0(\Oh_F, \Q)\otimes
\det^{-1}\hm^1(\Oh_F,\Q(1))\otimes\det^{-1}V_B^+.
\]
Let $T_B$ be the natural lattice in $V_B$ (see section \ref{artinmotives}).
 Recall
that $\hm^0(\Oh_F,\Q)=\Q$ and consider the splitting $s_\infty$ of
$(V_B^{\lor})_\R^+\to \hm^0(\Oh_F,\Q)_\R$ which maps $1$ to the dual of $\delta_\sigma$ for a fixed $\sigma$.
 Then $r_{\infty}\oplus s_\infty:\hm^1(\Oh_F,\Q(1))_\R\oplus \hm^0(\Oh_F,\Q)_\R\to (V_B^{\lor})_\R^+$ is an isomorphism and we
consider the lattice 
\[
\det\Z\otimes\det^{-1}\Oh_F^{*,\free}\otimes\det^{-1}T_B^+
\]
in $\Delta_f(V)$, where $\Oh_F^{*,\free}$ is the quotient of $\Oh_F^*$ by the group of 
roots of unity in $F$. 
Let $R_{\infty} (F):=
\vol\left((T_B^{\lor}\otimes\R)^+/(r_{\infty}(\Oh_F^{*,\free})\oplus s_\infty(\Z))\right)$,
then 
\[
\det r_{\infty}(\Oh_F^{*,\free})\otimes \det s_\infty (\Z)=R_{\infty} (F)\det (T_B^{\lor})^+
=R_{\infty} (F)\det^{-1}T_B^+.
\]
This allows to identify the element $\delta$ which maps to $1/\zeta_F(0)^*$ under the
isomorphism $\Delta_f(V)\otimes\R\isom \R$ as
\[
\delta\Z=\frac{R_{\infty} (F)}{\zeta_F(0)^*}\det\Z\otimes\det^{-1}\Oh_F^{*,\free}\otimes\det^{-1}T_B^+.
\]
The analytic class number formula gives
\[
\zeta_F(0)^*=-\frac{R_{\infty} (F)h}{\#\mu(F)},
\]
where $h=\# Cl(\Oh_F)$ is the class number. Hence
\[
\delta\Z= \det Cl(\Oh_F)\otimes\det\Z\otimes\det^{-1}\Oh_F^{*}\otimes\det^{-1}T_B^+.
\]
Now we define following Bloch and Kato 
$H^1_f(\Q_p, T_p):=\iota^{-1}(H^1_f(\Q_p, V_p))$, where 
$\iota:H^1(\Q_p, T_p)\to H^1(\Q_p, V_p)$. 
We put $H^1_{/f}:=H^1/H^1_f$. 
\begin{prop}[Appendix A]\label{flocalization} 
a) The Poitou-Tate localization sequence induces an exact sequence
\begin{multline*}
0\to \Oh_F^*\otimes\Z_p\to H^1(\Oh_F[1/p],\Z_p(1))\to 
H^1_{/f}(F\otimes\Q_p,\Z_p(1)) 
\to Cl(\Oh_F)\otimes\Z_p\to \\
\to H^2(\Oh_F[1/p],\Z_p(1))\to H^2(F\otimes\Q_p,\Z_p(1))\to
H^0(\Oh_F[1/p],\Q_p/\Z_p)^*\to 0.
\end{multline*}

b) Under the composition
\[ 
\det^{-1}R\Gamma_{/f}(F\tensor\Q_p,\Q_p(1))
\isom \det R\Gamma_{f}(F\tensor\Q_p,\Q_p)
\xrightarrow{\alpha} E_p
\]
of local duality with the isomorphism of \ref{localident} 
the lattice $\det^{-1}R\Gamma_{/f}(F\tensor\Q_p,\Z_p(1))$ is identified
with $\Z_p$. 
\end{prop} 
\bew See \ref{anhangPT} and \ref{anhanglocallattices}.\bewende
Let us show the Bloch-Kato 
conjecture for $h^0(F)$ and $r=0$: 
the exact sequence in part a) of the proposition implies 
\begin{multline*}
\det R\Gamma(\Z[1/p], T_p^{\lor} (1) )\otimes\det^{-1} R\Gamma_{/f}(\Q_p,T_p^{\lor}(1))\isom \\
\isom
\det^{-1}(\Oh_F^*\otimes\Z_p)\otimes\det (Cl(\Oh_F)\otimes\Z_p)\otimes\det^{-1}H^0(\Z[1/p], (T_p^{\lor})^*)^*.
\end{multline*}
If we tensor both sides with $\det^{-1}T_p^+$, use proposition \ref{ohnec}
and the fact that $1\in \Q\isom \hm^0(\Z, V)$ maps to 
$1\in \Z_p\isom H^0(\Z[1/p], (T_p^{\lor})^*)^*$ we get 
\[
\det R\Gamma_c(\Z[1/p], T_p)\otimes\det^{-1} R\Gamma_{/f}(\Q_p,T_p^{\lor}(1))\isom \delta\Z_p.
\]
The isomorphism $\Delta_f(V)\tensor\Q_p\isom \det R\Gamma_c(\Z[1/p], V_p)$
used in the formulation of the Bloch-Kato conjecture 
is the same as the composition 
\begin{multline*}
 \Delta_f(V)\tensor\Q_p\isom \det R\Gamma_c(\Z[1/p], V_p)
\tensor \det^{-1}R\Gamma_{/f}(\Q_p,V_p^{\lor}(1))\\
\isom \det R\Gamma_c(\Z[1/p], V_p)
\tensor \det R\Gamma_{f}(\Q_p,V_p)
\xrightarrow{\id\tensor\alpha}\det R\Gamma_c(\Z[1/p], V_p)
\end{multline*}
Hence part b) of the proposition implies that $\delta\Z_p\isom
\det R\Gamma_c(\Z[1/p], T_p)$ as required.

Now we consider the case $r=1$.
We deduce this case from the case $r=0$ by checking the compatibility
conjecture \ref{localconj} for $K=\Q$ in this special case. We keep the
same lattices as before. Let moreover
\[ T_\DR=\Oh_F \]
By the functional equation
\[ \frac{\zeta_F(0)^*}{\zeta_F(1)^*}=\pm \frac{d^{1/2}}{2^{r_1}(2\pi)^{r_2}} \]
where $d$ is the absolute value of the discriminant of $F$.
On the other hand, we consider again the natural lattice $T_B$ in $V_B$ 
and the
exact sequence
\[ 0\to 
((2\pi i)T_B)^+\tensor \R\to
\Oh_F\tensor\R 
\to T_B^+
\tensor \R\to 0 \]
The two lattices $\det \Oh_F$ and $\det T_B^+\tensor \det T_B(1)^+$
on $\det \Oh_F\tensor \R$ differ by 
the factor $d^{1/2}/2^{r_1}(2\pi)^{r_2}$. Hence, (up to powers of $2$) the element
$\epsilon(\Q/\Q,V(1))$ in conjecture \ref{localconj} is just
\[ \epsilon(\Q/\Q,V(1)) = \det^{-1} T_B(1)^+\tensor \det\Oh_F\tensor \det T_B^+ \]
It remains to show that the image of $\det \Oh_F$ is a generator of 
$\det R\Gamma(F\tensor\Q_p,\Z_p(F)(1))$ under the identification \ref{localident}.

We have defined before \ref{flocalization} integral structures such that
\[ \det R\Gamma(F\tensor\Q_p,\Z_p(1))=
\det R\Gamma_f(F\tensor\Q_p,\Z_p(1))\tensor \det R\Gamma_{/f}(F\tensor\Q_p,\Z_p(1))\]
and \ref{flocalization} also shows that 
\[ \det^{-1} R\Gamma_{/f}(F\tensor\Q_p,\Z_p(1)) \]
maps to $\Z_p$ under the identification in \ref{localident}.

The map
\[ \exp_p:\Oh_F\tensor \Q_p\to H^1_f(F\tensor \Q_p,\Q_p(1)) \]
is an isomorphism. By \ref{localfactor} we have
to check that
\[ \det \Oh_{F_v}\tensor \Z_p = (1-1/q)\det^{-1} R\Gamma_f(F_v,\Z_p(1)) \]
where $F_v$ is a completion of $F$ at a $p$-adic place and 
$q$ is the order of the residue class field of $F_v$.
Indeed, we have $H^1_f(F_v,\Z_p(1))=\Oh_{F_v}^*\tensor\Z_p$. The index of
$\exp_p(\Oh_{F_v})\subset \Oh_{F_v}^*\tensor\Z_p$ is $(q-1)/q=(1-1/q)$
because $\exp_p$ is a bijection of 
$\pi_v^n\Oh_{F_v}$ (with $\pi_v$ the uniformizer of $F_v$) towards the group $1+\pi_v^n\Oh_{F_v}$
for some $n\geq 1$.
This finishes the proof also in the case $r=1$.
\bewende

\section{Cyclotomic elements}

Cyclotomic elements play a decisive role in our paper. First they provide
explicit elements in Galois cohomology groups and second they form 
an Euler system, which will be used to prove the main conjecture.
The importance of the cyclotomic elements lies in the 
close connection with $L$-values via the explicit reciprocity law. Our approach
to the main conjecture through the Bloch-Kato conjecture explains
precisely, which Euler system is needed to give the equality in the 
main conjecture.

\subsection{The Euler system of cyclotomic elements}\label{Eulersystem}
Recall that $E$ is a finite extension of $\Q$ containing all values of
the Dirichlet character $\chi$ considered below. $\Oh$ is its ring of integers.
We fix a prime $p$ and 
put $E_p=E\tensor \Q_p$, $\Oh_p=\Oh\tensor\Z_p$.

Fix a collection $(\zeta_m)_{m\geq 1}$ of primitive $m$-th roots of 
unity in $\C$ satisfying $\zeta_{mn}^n=\zeta_m$. 
To be explicit we can take $\zeta_m=\exp\frac{2\pi i}{m}$.
Let $F=\Q(\zeta_N)$. We obtain a fixed embedding $\sigma_0:F\to \C$.

Let $\chi$ be a character of conductor $N$.
We need to fix a lattice in $V_B(\chi)$.
Recall that 
$h^0(F)_B$ are  the maps from the set of embeddings $\sigma:F\to \C$ to $\Q$ with
action of $g\in G:=\Gal(F/\Q)$ given by $(gf)(\sigma):=f(\sigma g)$. We write 
$\delta_{\sigma}$ for the delta function at $\sigma$. Note that we
have a distinguished element
$\delta_{\sigma_0}\in h^0(F)_B$.
\begin{defn}\label{generatordefn} Define the element $t_B(\chi)\in V_B(\chi)$ by 
\[
t_B(\chi):=p_{\chi^{-1}}\delta_{\sigma_0}
\]
and let $T_B(\chi):=\Oh t_B(\chi)$. We use this to define also
$T_p(\chi):=T_B(\chi)\otimes \Z_p$ and $t_p(\chi):=t_B(\chi)\otimes 1\in T_p(\chi)$.
\end{defn}

\rem Recall that we have by definition $V(\chi)=p_{\chi^{-1}} h^0(F)$ and
with this normalization $V_p(\chi)$ is the Galois module with operation
of $G_\Q$ via $\chi$.\\

Following Soul{\'e} we define cyclotomic elements 
\[
c_r(\zeta_m)\in H^1(\Z[\zeta_m][1/p], \Z_p(r)),
\]
for all  $r\in\Z$ and $m\geq 1$.
 For every integer 
$m\geq 1$ and prime $l$ the  
following norm compatibility is satisfied:
\begin{equation*}
\Norm_{\Q(\mu_{ml})/\Q(\mu_m)}(1-\zeta_{ml}) =\left\{\begin{array}{ll}
1-\zeta_{m}& \mbox{if } l|m\\
(1-\zeta_{m})^{1-\Fr_l}& \mbox{if } l\nmid m\mbox{ and } m>1\\
l&\mbox{if }m=1
\end{array}\right.
\end{equation*}
Here $\Fr_l$ is the (geometric) Frobenius at $l$ in $\Gal(\Q(\zeta_m)/\Q)$.
For $n\geq 1$ and $m\geq 1$ the elements $1-\zeta_{p^nm}$ are units
in $\Z[\zeta_{p^nm}][1/p]$ 
and we can consider for every $n\geq 1$
\begin{align*}
(1-\zeta_{p^nm})\otimes(\zeta_{p^nm}^m)^{\otimes r-1}&\in 
H^1(\Z[\zeta_{p^nm}][1/p], \Z/p^n\Z(1))\otimes \Z/p^n\Z(r-1)\\
&=H^1(\Z[\zeta_{p^nm}][1/p], \Z/p^n\Z(r)).
\end{align*}
We let 
\[
c_r(\zeta_m)_n:=\cores_{\Q(\zeta_{p^nm})/\Q(\zeta_m)}
(1-\zeta_{p^nm})\otimes(\zeta_{p^nm}^m)^{\otimes r-1}
\in H^1(\Z[\zeta_m][1/p], \Z/p^n\Z(r)).
\]
The $c_r(\zeta_m)_n$ are compatible with respect to the maps 
$\Z/p^n\Z\to \Z/p^{n-1}\Z$. 
\begin{defn} For $r\in \Z$ and $m\geq 1 $,  
the Soul{\'e}-Deligne {\em cyclotomic elements } $c_r(\zeta_m)$ are defined
as
\[
c_r(\zeta_m):=\prolim_n c_r(\zeta_m)_n \in 
H^1(\Z[\zeta_m][1/p], \Z_p(r)).
\]
\end{defn}

We are now going to define elements 
$c_r(\zeta_m)(\chi)\in H^1(\Z[\zeta_{m}][1/p], T_p(\chi)(r))$. 
Let $N$ be the conductor of $\chi$ and $K$ the least common multiple
of $N$ and $m$.
Let $F:=\Q(\zeta_N)$, then $F(\zeta_m)=\Q(\zeta_K)$. Let $t_p(\chi)$ be the
$\Oh_p$ generator of $T_p(\chi)$ fixed in \ref{generatordefn}. 
\begin{defn}\label{cdefn}Define $c_r(\zeta_m)(\chi)$ to be
the image of $c_r(\zeta_K)\tensor t_p(\chi) $ under the composition
\begin{align*}
H^1(\Z[\zeta_K][1/p], \Z_p(r))\otimes T_p(\chi)
&\isom H^1(\Z[\zeta_K][1/p],T_p(\chi)(r))\\
&\xrightarrow{\cores}H^1(\Z[\zeta_m][1/p],T_p(\chi)(r)).
\end{align*}
\end{defn}
Note that the elements $c_r(\zeta_m)$ are modulo torsion 
in the $(-1)^{r-1}$-eigenspace of complex
conjugation, so that $c_r(\zeta_m)(\chi)$ can be non torsion only for 
$\chi(-1)=(-1)^{r-1}$. 
If $m=1$ we have $c_r(1)(\chi) \in H^1(\Z[1/p],T_p(\chi)(r))$.

\rem Rationally, we have the isomorphism 
$H^1(\Z[\zeta_N][1/p], \Q_p(r))=\bigoplus_{\chi} H^1(\Z[1/p],V_p(\chi))$ where
the sum is taken over all characters of conductor dividing $N$. 
Consider a character $\chi$ such that the primes dividing its conductor
are the same as the primes dividing $N$ with the possible exception of $p$.
Then 
$c_{r}(1)(\chi)$ is the $\chi$-component of
$c_r(\zeta_N)$, i.e., $c_{r}(1)(\chi)=p_{\chi^{-1}}c_r(\zeta_N)$ up to
torsion, by the norm compatibility of cyclotomic elements.


Our aim is to show that the elements $c_r(\zeta_m)(\chi)$
 form an Euler system for $(T_p(\chi), pN)$. 
Euler systems were invented by Kolyvagin. A general theory of 
Euler systems was developed by Kato \cite{KatoEuler}, 
Perrin-Riou \cite{PerrinEuler} and Rubin \cite{Ru}. We follow Rubin,
because his approach is closest to our setting.

Let us recall the definition of an Euler system:
\begin{defn} An {\em Euler system} for $(T_p(\chi), pN)$
is a  collection of elements
\[
e_r(m)\in H^1(\Q(\zeta_m), T_p(\chi)(r))
\]
for all $m\geq 1$, such that for all primes $l$
\[
\cores_{\Q(\zeta_{ml})/\Q(\zeta_m)}(e_r(ml)) =
\left\{ \begin{array}{ll}
e_r(m)&\mbox{ if } l|mpN\\
(1-\chi^{-1}(l)l^{r-1})e_r(m)&\mbox{ if } l\nmid mpN.
\end{array}\right.
\]
\end{defn} 
\begin{lemma}(Soul{\'e} \cite{Soulepadic}) 
The elements $c_r(\zeta_m)(\chi)$   
form an Euler system for $(T_p(\chi)(r),pN)$.
\end{lemma}
\bew This follows from  proposition 2.4.2. in \cite{Ru} 
and the norm compatibility of the cyclotomic units.
\bewende

We need also the following variant of  $c_r(\zeta_N)$ and $c_r(1)(\chi)$, see \cite{Kato Main conj} {\S} 5.

\begin{defn}\label{ctildedefn} Define 
\[\ctilde_r(\zeta_N)=\begin{cases}
c_r(\zeta_N) &\text{if $p|N$}\\[1ex]
\sum_{i\geq 0}(p^{r-1})^ic_r(\zeta_N^{p^{-i}})
    &\text{if $p\nmid N$, $r\geq 2$}\\[1ex]
-\sum_{i\geq 1}(p^{1-r})^ic_r(\zeta_N^{p^{i}}) &\text{if $p\nmid N$, $r\leq 0$}
\end{cases}
\]
Here  $\zeta_N^{p^{-i}}$ denotes the unique $N$-th root of unity whose $p^i$-th
power is $\zeta_N$. Similarly we let
\[
\ctilde_r(1)(\chi)=\begin{cases}
c_r(1)(\chi) &\text{if $p|N$}\\[1ex]
\sum_{i\geq 0}(p^{r-1})^i\chi^{-1}(p)^ic_r(1)(\chi)&\text{if $p\nmid N$, $r\geq 2$}\\[1ex]
-\sum_{i\geq 1}(p^{1-r})^i\chi(p)^ic_r(1)(\chi) &\text{if $p\nmid N$, $r\leq 0$}
\end{cases}
\]
\end{defn}
Note that the sums converge in $H^1(\Z[\zeta_N][1/p],\Z_p(r))$ and  
$H^1(\Z[1/p],T_p(\chi)(r))$ respectively. 
If $p\nmid N$,
\[
(1-p^{r-1}\Fr_p)\ctilde_r(\zeta_N)=c_r(\zeta_N)
\]
and 
\[
(1-p^{r-1}\chi^{-1}(p))\ctilde_r(1)(\chi)=c_r(1)(\chi).
\]
The following lemma allows the comparison with the $p$-adic regulator from
$K$-theory.
\begin{lemma}\label{Kvergleich} The following identity holds 
in $H^1(\Z[\zeta_{p^nN}][1/p],\Z/\Z p^n(r))$
if $p\nmid N$ and $r>1$:
\[
\sum_{\beta^{p^n}=\zeta_N}(1-\beta)\otimes(\beta^N)^{\otimes r-1}=
\sum_{i=0}^n(p^{r-1})^ic_r(\zeta_N^{p^{-i}})_{n}.
\]
For $p|N$ 
\[
\sum_{\beta^{p^n}=\zeta_N}(1-\beta)\otimes(\beta^N)^{\otimes r-1}=c_r(\zeta_N)_n.
\]
In particular the elements $\ctilde_r(\zeta_N)$ and 
\[
\left(\sum_{\beta^{p^n}=\zeta_N}(1-\beta)\otimes(\beta^N)^{\otimes r-1}
\right)_n
\]
are the same in $H^1(\Z[\zeta_{p^\infty N}][1/p],\Z_p(r))$.
\end{lemma} 
\bew In the case $p|N$ there is nothing  to show, so assume $p\nmid N$.
Write 
$H_i(\zeta_N):=\{\beta|\beta^{p^{n-i}N}=1\mbox{ and } \beta^{p^n}=\zeta_N\}$. 
Then 
\[
\sum_{H_0(\zeta_N)}(1-\beta)\otimes(\beta^N)^{\otimes r-1}=
\sum_{i=0}^n\sum_{ H_i(\zeta_N)\ohne H_{i+1}(\zeta_N)}
(1-\beta)\otimes(\beta^N)^{\otimes r-1}.
\]
We have 
\[
\sum_{ H_i(\zeta_N)\ohne H_{i+1}(\zeta_N)}
(1-\beta)\otimes(\beta^N)^{\otimes r-1}=
(p^{r-1})^i\sum_{H_0(\zeta_N^{p^{-i}})\ohne H_1(\zeta_N^{p^{-i}})}
(1-\beta)\otimes(\beta^N)^{\otimes r-1}
\]
and the last sum is nothing but $c_r(\zeta_N^{p^{-i}})_{n}$.
\bewende

\rem  Note that
the $c_r(\zeta_N)$ are the norm compatible elements, which define Euler
systems and hence are related to $p$-adic $L$-functions. They miss the
Euler factor at $p$. On the other hand, the elements 
$\ctilde_r(\zeta_N)$  are known to
come from elements in integral $K$-theory (see \ref{pregulator}). These
are the elements which appear in \cite{Hu-Wi}, corollary 9.7 or \cite{Bloch-Kato} 6.2.

\subsection{Relation to $L$-values and non-vanishing of the Euler system}
We want to relate the elements $\ctilde_r(\zeta_N)$ defined in
\ref{ctildedefn} 
to $L$-values via the exponential map. This will also show that they are non
torsion.

\begin{defn} For $\alpha$ an $N$-th root of unity let 
\[
Li_s(\alpha):= \sum_{n\geq 1} \frac{\alpha^n}{n^s}
\]
be the {\em polylogarithm function}. The sum converges for $Re\; s> 1$.
\end{defn}

On the other hand we have the Hurwitz zeta function, defined as follows:
\begin{defn}
Let for $c\in \Z$
\[
\zeta_N(s,c):=\sum_{n\geq 1,\; n\equiv c\  \mathrm{mod}\ N}\frac{1}{n^s}
\]
be the {\em Hurwitz zeta function} of modulus $N$ with respect to $c$.
The sum  converges for $Re \; s>1$.
\end{defn}
\begin{thm}\label{hurwitz} The functions $Li_s(\alpha)$ and $\zeta_N(s, c)$
have meromorphic continuations to $\C$ and for $r\in \Z$ and
$\alpha=\exp(\frac{2\pi ic}{N})$ a primitive $N$-th root of unity, 
the equality 
\[
\frac{1}{2} \left( Li_{1-r}(\alpha) + (-1)^rLi_{1-r}((\alpha^{-1}))\right)=
\left(\frac{-2\pi i}{N}\right)^{-r}\lim_{s\to r}\Gamma(s)\left(
\zeta_N(s, c)+(-1)^r\zeta_N(s,-c)\right)
\]
holds.
\end{thm} 
\bew This is the well known functional equation of the Hurwitz zeta function,
see e.g. \cite{Apostol} thm. 12.6.
\bewende
Let $L(\chi,s)$ be the $E\tensor_\Q\C$-valued $L$-function of the primitive character $\chi$. The above theorem
immediately implies with our normalization $\chi(\Fr_p)=\chi(p)$: 
\begin{cor} \label{Li-ident}
Let $\chi$ be a character with $\chi(-1)=(-1)^r$ and 
$r\geq 1$,  then
\[
\sum_{\tau}\chi^{-1}(\tau)\otimes Li_{1-r}(\tau \zeta_N)=2\left(\frac{-2\pi i}{N}\right)^{-r}
(r-1)!L(\chi, r)
\]
as elements in $E\tensor_\Q\C$.
\end{cor}

Consider for $r\geq 1$ the composition
\begin{equation*}
H^1(\Z[\zeta_{N}][1/p], {\Oh_p}(1-r))\to 
H^1(\Q_p\otimes\Q(\zeta_{N}), E_p(1-r))
\xrightarrow{\exp^*_p}V_{\DR}(\Q(\zeta_{N}))^{\lor}\otimes_E E_p.
\end{equation*}
where $\exp_p^*$ is the dual of $\exp_p$ via local duality. 
Let us define an element in 
$V_{\DR}(\Q(\zeta_N))^{\lor}$.

\begin{defn}
Let $d_{1-r}(\zeta_N)^{\lor}\in V_{\DR}(\Q(\zeta_N))^{\lor}$ be the $E$-linear
 map
\[
x\mapsto \Tr_{\Q(\zeta_N)/\Q}(x Li_{1-r}(\zeta_N)).
\]
\end{defn}
The next result is crucial for our approach to the Bloch-Kato conjecture.
\begin{thm}\label{expcomputation}(Kato \cite{Kato Main conj} theorem 5.12) For 
$r\geq 1$ the image of $ \ctilde_{1-r}(\zeta_N)$
under the  map $\exp_p^*$ is 
$(-1)^{r-1}N^{-r}(r-1)!^{-1}d_{1-r}(\zeta_N)^{\lor}$.
In particular the element $\ctilde_{1-r}(\zeta_N)$ is not torsion in 
$H^1(\Z[1/p],\Z_p(1-r))$.
\end{thm}
\bew This is a standard consequence of the explicit reciprocity
law of Kato, \cite{galaxy} II, Theorem 2.1.7. For the calculation see e.g. Kato \cite{galaxy} chapter
III, 1.1.7. 

As Benois pointed out to us, Theorem 2.1.7 in loc. cit. is proved only for norm compatible systems of invertible elements of 
$u_n\in \Oh_{K_n}$ (notation of loc.cit). However, as we learned from Colmez,
the argument also works for a norm compatible system of elements $u_n\in \Oh_{K_n}\smallsetminus\{0\}$.
Then the case $N=p^n$ (where $1-\zeta_{p^n}$ is not a unit) is also covered by the theorem.
\bewende
Let $p_{\chi}:=\frac{1}{\# G}\sum_{\tau}\chi^{-1}(\tau)\tau$
 be the
projector of $V_{\DR}(\Q(\zeta_N))$ onto $V_{\DR}(\chi^{-1})$. 
We consider
\[
p_{\chi}\left(d_{1-r}(\zeta_N)^{\lor}(\chi^{-1})\right)\in 
V_{\DR}(\chi)^{\lor}.
\]
From the above, we have:
\begin{cor}\label{recchi}
For $r\geq 1$ and $\chi(-1)=(-1)^r$,  the image of $\ctilde_{1-r}(1)(\chi^{-1})$ under
\[ \exp_p^*: H^1(\Z[1/p],V_p(\chi^{-1})(1-r))\to V_\DR(\chi)^\lor \]
maps to the linear map which maps $x\in V_\DR(\chi)$ to
\[ 
- \Tr_{\Q(\zeta_N)/\Q}\left(
  \frac{2}{\# G}(2\pi i)^{-r} L(\chi, r)x
\right).
\]
In particular, as $L(\chi, r)\neq 0$ the element $\ctilde_{1-r}(1)(\chi^{-1})$ is not torsion in 
\[ H^1(\Z[1/p],T_p(\chi^{-1})(1-r)).\]
\end{cor}
\bew 
$\ctilde_{1-r}(1)(\chi^{-1})$ maps to 
$(-1)^{r-1}N^{-r}(r-1)!^{-1} p_\chi(d_{1-r}(\zeta_N)^{\lor})\in V_{\DR}(\chi)^{\lor}$ by the theorem. By \ref{Li-ident}
this is the element which maps $x$ to 
\[
- \Tr_{\Q(\zeta_N)/\Q}\left( \frac{2}{\# G}(2\pi i)^{-r}
L(\chi, r)x\right).
\]
\bewende

\subsection{Reformulation of the Bloch-Kato conjecture}
We reformulate the Bloch-Kato conjecture for 
characters $\chi$ with $\chi(-1)=(-1)^r$  and $r\geq 1$
using the elements $\ctilde_{1-r}(1)(\chi^{-1})$ in 
$H^1(\Z[1/p],T_p(\chi^{-1})(1-r))$ defined in \ref{ctildedefn}.

\begin{lemma}\label{finitecoh}Let $r\geq 1$ and $\chi(-1)=(-1)^r$, then
$H^0(\Z[1/p],T_p(\chi^{-1})(1-r))$ and $H^2(\Z[1/p],T_p(\chi^{-1})(1-r))$ are
finite.
\end{lemma}
\bew According to \ref{recchi} the elements
\[
(1-\chi(p)p^{-r})\ctilde_{1-r}(1)(\chi^{-1})=c_{1-r}(1)(\chi^{-1})
\]
are non torsion  and they are the first layer of an Euler system in $H^1(\Z[1/p],T_p(\chi^{-1})(1-r))$.
Theorem 2.2.3 in \cite{Ru} implies then that 
\[
\ker\left(H^2(\Z[1/p],T_p(\chi^{-1})(1-r))\to H^2(\Q_p,T_p(\chi^{-1})(1-r))\right)
\]
is finite. By local duality and our conditions on $\chi$, the
group $H^2(\Q_p,T_p(\chi^{-1})(1-r))$ is finite as well. The statement
about $H^0(\Z[1/p],T_p(\chi^{-1})(1-r))$ is clear.
\bewende

Thus, we have an isomorphism
\[
{\det}_{E_p}H^1(\Z[1/p],V_p(\chi^{-1})(1-r))\isom
{\det}^{-1}_{\Oh_p}R\Gamma(\Z[1/p],T_p(\chi^{-1})(1-r))\otimes\Q_p
\]
and we can consider $c_{1-r}(1)(\chi^{-1})\in {\det}^{-1}_{\Oh_p}R\Gamma(\Z[1/p],T_p(\chi^{-1})(1-r))\otimes\Q_p$.

Our aim is to prove the following theorem:
\begin{thm}\label{BlochKatochar} Let $\chi$ be a character with conductor $N$
and with $\chi(-1)=(-1)^{r}$ and assume 
that $r\geq 1$. Then the Bloch-Kato
conjecture for $V(\chi)$ and $r$ is true up to powers of $2$,
 if and only if for all $p\neq 2$, the cyclotomic element
$c_{1-r}(1)(\chi^{-1})$ is a generator of ${\det}^{-1}_{\Oh_p}R\Gamma(\Z[1/p],T_p(\chi)(1-r))$.
\end{thm}
\bew
We are in the case $\hm^1(\Z,V(\chi)(r))=0$. Hence,
the fundamental line for $V(\chi)(r)$ reduces to 
\[
\Delta_f(V(\chi)(r))={\det}_{E}V_{\DR}(\chi)\otimes
{\det}^{-1}_{E}V_B(\chi)(r)
\]
Our first aim is to describe the element $\delta$ of conjecture \ref{conj}
explicitly. 
Let $t_B(\chi)$ be  the 
$\Oh$-generator of $T_B(\chi)$ fixed in \ref{generatordefn}.
Recall that we have a distinguished embedding $\sigma_0:\Q(\zeta_N)\to \C$ by
our choice of root of unity.
Let $I_\infty:V_\DR(\chi)_\R\to V_B(\chi)(r)_\R$ be the comparison isomorphism and 
$f\in p_{\chi^{-1}}(\Q(\zeta_N)\tensor E)$. Then, for all $\sigma\in G=\Gal(\Q(\zeta_N)/ \Q)$, $\sigma(f)=\chi^{-1}(\sigma)f$ and $f$ is 
mapped to
\begin{align*}
 \sum_\sigma\sigma_0\sigma(f)\delta_{\sigma_0\sigma} &=
\sum_{\sigma}\chi^{-1}(\sigma)\sigma_0(f)\delta_{\sigma_0\sigma}\\
&= (2\pi i)^{-r}(\# G)\sigma_0(f)
 \left(\frac{1}{\# G}\sum_\sigma \chi^{-1}(\sigma) \sigma^{-1}\delta_{\sigma_0}
(2\pi i)^r\right)\\
&=(2\pi i)^{-r}(\# G)\sigma_0(f)t_B(\chi)(r).
\end{align*}
Conversely, $I_\infty^{-1}(t_B(\chi)(r))=(2\pi i)^{r}(\# G)^{-1}$ where
we use $\sigma_0$ to interpret elements of $V_\DR(\chi)_\R$ as complex
numbers. By definition
\[ \delta= \frac{1}{L(\chi,r)} I_\infty^{-1}(t_B(\chi)(r)) \tensor (t_B(\chi)(r))^{-1}
= \frac{(2\pi i)^{r}}{ (\# G)L(\chi,r) }\tensor (t_B(\chi)(r))^{-1}
\]
in $\Delta_f(V(\chi)(r))_\R$.
We identify $\det V_\DR(\chi)=\det^{-1} V_\DR(\chi)^\lor$. In this description,
$\delta= ( v^\lor)^{-1} \tensor (t_B(\chi)(r))^{-1}$ where
$v^\lor$ is multiplication by $(2\pi i)^{-r}(\# G)L(\chi,r)$.
By \ref{recchi}, the element $\exp_p^*(\ctilde_{1-r}(1)(\chi^{-1}))\in V_\DR(\chi)^\lor$ maps $x$ to 
\[  - \Tr_{\Q(\zeta_N)/\Q}\left(
  \frac{2}{\# G}(2\pi i)^{-r} L(\chi, r)x \right) = - 2 (2\pi i)^{-r} L(\chi, r)x \]
Hence we get
\[ \delta= \left( -\frac{\# G}{2} \exp_p^*(\ctilde_{1-r}(1)(\chi^{-1}) )\right)^{-1}
\tensor(t_B(\chi)(r))^{-1} \]
Duality on $h^0(\Q(\zeta_N))$ induces an
 isomorphism $V_B(\chi^{-1})\isom V_B(\chi)^\lor$.  
Under this duality the standard lattice $T_B\subset h^0(\Q(\zeta_N))_B$ is self-dual. This implies that $t_B(\chi)^\lor=(\# G) t_B(\chi^{-1})$, i.e,
\[ \delta= \left( -\frac{1}{2} \exp_p^*(\ctilde_{1-r}(1)(\chi^{-1}) )\right)^{-1}
\tensor t_B(\chi^{-1})(-r) \]
Consider the Bloch-Kato conjecture for the lattice $T_p(\chi^{-1})^\lor(r)$ in $V_p(\chi)(r)$.
It holds if and only if $\delta$ is a generator
of ${\det}_{\Oh_p}R\Gamma_c(\Z[1/p],T_p(\chi^{-1})^\lor(r))$, 
 or equivalently
by  proposition \ref{ohnec} a generator of
\[
{\det}_{\Oh_p}R\Gamma(\Z[1/p],T_p(\chi^{-1})(1-r))\otimes{\det}_{\Oh_p}
T_p(\chi^{-1})(-r).
\]
This holds if and only if $-\frac{1}{2} \exp_p^*(\ctilde_{1-r}(1)(\chi^{-1}))$
corresponds to a generator of ${\det}^{-1}_{\Oh_p}R\Gamma(\Z[1/p],T_p(\chi^{-1})(1-r))$ via $\exp_p$ and using the map on local factors in \ref{localident}. By \ref{localfactor} this
is equivalent to 
\[ -\frac{1}{2}
c_{1-r}(1)(\chi^{-1})=-\frac{1}{2}(1-\chi(p)p^{-r})\ctilde_{1-r}(\zeta_N)(\chi^{-1})
\]
 being a generator of ${\det}^{-1}_{\Oh_p}R\Gamma(\Z[1/p],T_p(\chi^{-1})(1-r))$
using the natural map.
Hence the Bloch-Kato conjecture is reduced to $-\frac{1}{2}
c_{1-r}(1)(\chi^{-1})$ being a generator of
${\det}^{-1}_{\Oh_p}R\Gamma(\Z[1/p],T_p(\chi^{-1})(1-r))$.
\bewende

Let $\Q_n$ be the cyclotomic $\Z/p^n$-extension. Let $\Z_n$ be its ring of integers. We denote by
$c_{1-r}(\Q_n/\Q)(\chi^{-1})$
the image of $c_{1-r}(\zeta_{p^{n+1}})(\chi)$ under the corestriction
\[
H^1(\Z[\zeta_{p^{n+1}}][1/p],T_p(\chi^{-1})(1-r))\to H^1(\Z_n[1/p],T_p(\chi^{-1})(1-r)).
\]
Note that these elements are norm-compatible for varying $n$.

\begin{cor}\label{equivreform}
Let $\chi$ and $r$ be as in the theorem, $p\neq 2$ a prime. Then the $p$-part of the equivariant Bloch-Kato conjecture \ref{equivconj} is
true for the motive $V(\chi)$ and $r$ and $\Q_n/\Q$ if and only if $c_{1-r}(\Q_n/\Q)(\chi^{-1})$ is a generator
of ${\det}^{-1}_{\Oh_p[\Gal(\Q_n/\Q)]}R\Gamma(\Z_n[1/p],T_p(\chi^{-1})(1-r))$.
\end{cor}
\bew 
Consider the characters $\omega\chi$ where $\omega$ runs through the characters of $G=\Gal(\Q_n/\Q)$. As $\Q_n$ is totally real,
the parity condition is unchanged. 
Let $N$ be the conductor of $\chi$ and
$N_{\omega\chi}$ the conductor of $\omega\chi$. They differ by a $p$-power.
 We can apply 
the arguments in the proof of the theorem to the $V(\omega\chi)$'s. Hence
the image of $\delta(\Q_n/\Q)$ under $r_p$ with the local Euler factors from \ref{localident} taken into account has $\omega$-component
\[
\left(-\frac{1}{2}c_{1-r}(1)(\omega^{-1} \chi^{-1})\right)^{-1}
      \tensor \; t_p(\omega^{-1}\chi^{-1})(-r).
\]
This is nothing but the $\omega$-component of 
\[ \left( -\frac{1}{2}c_{1-r}(\Q_n/\Q)(\chi^{-1}) \right)^{-1}
\tensor t_p(\Q_n/\Q)(\chi^{-1})(-r) \]
with 
$t_p(\Q_n/\Q)(\chi^{-1})=t_p(\chi^{-1})\tensor 1\in T_p(\chi^{-1})\otimes\Oh_p[G]$ .
The equivariant Bloch-Kato conjecture
is equivalent to $\delta(\Q_n/\Q)$ being a generator of
\[ \det_{\Oh_p[G]}R\Gamma(\Z_n[1/p],T_p(\chi)^\lor(1-r)\tensor {\det}_{\Oh_p[G]}
(T_p(\chi^{-1})\otimes\Oh_p[G])(-r) \]
As in the absolute case this is equivalent to the assertion.
\bewende

\begin{cor}\label{BlochKato} Let $X$ be the set of characters $\chi$ with conductor 
dividing a fixed $N$ and $\chi(-1)=-1$. Then
the element $\prod_{\chi\in X}c_{0}(1)(\chi)$ is 
a generator of 
\[
{\det}_{\Oh_p}R\Gamma(\Z[1/p],\bigoplus_{\chi\in X}T_p(\chi)).
\]
\end{cor}
\bew By \ref{classnumber} we know that the Bloch-Kato conjecture is true for $V(\Q(\zeta_N))$ and $r=1$ and for the motive of
the real subfield $V(\Q(\zeta_N)^+)$.
Thus it is true 
for 
\[
V(\Q(\zeta_N))^-\isom\bigoplus_{\chi\in X}V(\chi)
\]
and $r=1$
with an arbitrary 
choice of lattice in $V(\Q(\zeta_N))^-$.
Repeating the proof of the last theorem for the direct sum of the $V(\chi)$ 
implies the claim.
\bewende

\rem This simple argument is one of the key insights of this
paper. We originally have proved the Bloch-Kato conjecture for the motive
of a number field using certain choices of lattices. We now apply it
with a completely different choice of lattice which is compatible with
direct sum decomposition. The relation of the two lattices is non-trivial --
but we do not need to know anything about the comparison factors. The
extra freedom in the Bloch-Kato conjecture allows to avoid these
computations.

\section{The main conjecture} \label{mc}
In this section we formulate and prove the main conjecture in Iwasawa theory
for all characters. It is essential for our proof of the Bloch-Kato
conjecture to have the main conjecture for all characters. 
\subsection{Iwasawa modules}
The theory of Euler systems gives relations between certain
Iwasawa modules. These modules will be defined in this section.
We start by defining the Iwasawa algebra:

Recall
that $E$ is a finite extension of $\Q$ containing all values of the
Dirichlet character $\chi$ of conductor $N$.
Let ${\Oh_p}:=\Oh\otimes_{\Z}\Z_p$ be the ring of integers in  $E_p$. 
This is a product of 
discrete valuation rings. 
\begin{defn}
Denote by $\Q_{\infty}$ the maximal $\Z_p$-extension
inside $\Q(\zeta_{p^{\infty}})$ and by $\Q_n$ the finite extensions of $\Q$ inside
$\Q_{\infty}$ with $\Gal(\Q_n/\Q)=\Z/p^n\Z$. Let $\Gamma:=\Gal(\Q_{\infty}/\Q)$ and  
$\Lambda:=\prolim_n{\Oh_p}[\Gal(\Q_{n}/\Q)]$. 
The algebra $\Lambda$ is the {\em Iwasawa algebra}.
\end{defn}

 Then
we have the standard identification 
\[
\Lambda\isom{\Oh_p}[[t]]
\]
with a power-series ring in one variable over ${\Oh_p}$. 
In particular $\Lambda$
is a regular ring. 

The following $\Lambda$ modules are our main object of study: 

Let $\chi$ be a Dirichlet character and $T_p(\chi)$ an
${\Oh_p}$-lattice in $V_p(\chi)$.
Define:
\begin{equation*} 
T_p(\chi)^*:=\Hom_{{\Oh_p}}(T_p(\chi),E_p/\Oh_p)
\end{equation*}
the ${\Oh_p}$-Pontryagin  dual of $T_p(\chi)$.
\begin{defn}\label{fettH} Let $\Z_n$ be the ring of integers in the $\Z/p^n\Z$-extension
$\Q_n$ of $\Q$. Define
\begin{align*}
\Hh^q_{\gl}( T_p(\chi)(k))&:=
\prolim_n 
H^q(\Z_{n}[1/p], T_p(\chi)(k))\\
\Hh^q_{\loc}(T_p(\chi)(k))&:=
\prolim_n H^q(\Q_p\otimes\Q_{n}, T_p(\chi)(k))\\
\Hh^q_{\gl}(T_p(\chi)^*(1-k))^* &:=
H^q(\Z_\infty[1/p],T_p(\chi)^* (1-k))^*,
\end{align*}
where the limit is taken with respect to the corestriction maps.
\end{defn}
In the sequel we collect some facts about these $\Lambda$-modules.
 Let 
\[
\epsilon_{\cycl}:\Gal(\Qbar/\Q)\to \Z_p^{\mal}
\]
be the cyclotomic character (i.e., $\Oh_p(1)=\Oh_p(\epsilon_{\cycl})$)
and write 
$\epsilon_{\cycl}=\epsilon\times\epsilon_{\infty}$ according to the 
decomposition $\Z_p^{\mal}\isom(\Z/p\Z)^{\mal}\times\Z_p$. We write 
$\Oh_p(\epsilon_{\infty})$ for the $\Oh_p$ module of rank $1$ with Galois
action given by $\epsilon_{\infty}$. 

\begin{lemma}\label{twist} There are  isomorphisms of $\Lambda$-modules
\[
\Hh^q_{\gl}( T_p(\chi)(k))\otimes\Oh_p(\epsilon_{\infty})
\isom \Hh^q_{\gl}(T_p(\chi\epsilon^{-1})(k+1))
\]
and
\[
\Hh^q_{\loc}( T_p(\chi)(k))\otimes\Oh_p(\epsilon_{\infty})
\isom \Hh^q_{\loc}(T_p(\chi\epsilon^{-1})(k+1)).
\]
The parity of $\chi\epsilon^{-1}$ is minus the parity of $\chi$.
\end{lemma}
\bew See \cite{Ru} 6.2.1. 
\bewende
\begin{lemma}\label{vanishing} The following $\Lambda$-modules are zero:
\[
\Hh^2_{\gl}( {T_p(\chi)}^*(1-k))^*=0=
\Hh^0_{\loc}( {T_p(\chi)}(k)).
\]
\end{lemma}
\bew The $\Lambda$-module $\prolim_nH^2(\Z[\zeta_{p^nN}][1/p], {E_p/\Oh_p}(1-k))$ 
is independent of $k$ and
is zero by \cite{Schneider} paragraph 4, lemma 7.

 As the functor $\prolim$ is exact in our situation and
$H^2$ is right exact for $p\neq 2$, we get a surjection
$\prolim_nH^2(\Z[\zeta_{p^nN}][1/p], {E_p/\Oh_p}(1-k))\to\Hh^2_{\gl}( {T_p(\chi)}^*(1-k))$.
The group
$\Hh^0_{\loc}( {T_p(\chi)}(k))$ is the inverse
limit of $H^0(\Q_p\otimes\Q_{n}, {T_p(\chi)}(k))$, 
which is zero for
$k\neq 0$ for weight reasons. With lemma \ref{twist} we get the
statement also for $k=0$. 
\bewende
We define, following Kato, a Selmer group:
\begin{defn}\label{0defn} Let
\[
\Hh^2_{\gl,0}({T_p(\chi)}(k)):=
\coker\left(\Hh^1_{\loc}( {T_p(\chi)}(k))\to\Hh^1_{\gl}( {T_p(\chi)}^*(1-k))^*
\right)
\]
and 
\[
\Hh^2_{\loc,0}({T_p(\chi)}(k)):=
\coker\left(\Hh^2_{\gl,0}( {T_p(\chi)}(k))\to
\Hh^2_{\gl}({T_p(\chi)}(k))\right).
\]
\end{defn}
The following sequence gives a connection between the $\Lambda$-modules
defined above. 
\begin{lemma}\label{classfieldseq}The Poitou-Tate localization sequence
induces for all $k\in\Z$ exact sequences
\begin{gather*}
0\to \Hh^1_{\gl}({T_p(\chi)}(k))\to 
\Hh^1_{\loc}({T_p(\chi)}(k))\to
\Hh^1_{\gl}({T_p(\chi)}^*(1-k))^*
\to \Hh^2_{\gl,0}( {T_p(\chi)}(k))\to 0,\\
0\to \Hh^2_{\loc,0}({T_p(\chi)}(k))\to \Hh^2_{\loc}({T_p(\chi)}(k))
\to \Hh^0_{\gl}({T_p(\chi)}^*(1-k))^*\to 0.
\end{gather*}
\end{lemma}
\bew By biduality the functor $\prolim$ is exact on the Poitou-Tate 
localization sequence. It suffices to see that 
\[
\Hh^2_{\gl}( T_p(\chi)^*(1-k))^*=0,
\]
which is the content of lemma \ref{vanishing}.
\bewende
 The elements 
\[
c_k(\Q_n/\Q)(\chi)=\cores_{\Q(\zeta_{p^{n+1}})/\Q_n}c_k(\zeta_{p^{n+1}})(\chi)\in 
H^1(\Z_{n}[1/p], {T_p(\chi)}(k))
\]
are compatible with corestriction.
\begin{defn}\label{iwasawac}The {\em cyclotomic element } in 
$\Hh^1_{\gl}({T_p(\chi)}(k))$ is defined to be
\[
c_k(\chi):=\prolim_n c_k(\Q_n/\Q)(\chi)\in
\Hh^1_{\gl}({T_p(\chi)}(k)).
\]
\end{defn}
\begin{lemma}\label{twist2} Under the isomorphism 
$\Hh^1_{\gl}({T_p(\chi)}(k))\isom\Hh^1_{\gl}({T_p(\chi\epsilon^{-1})}(k+1))$
 of lemma 
\ref{twist} the element  $c_k(\chi)$ maps to
$c_{k+1}(\chi\epsilon^{-1})$.
\end{lemma}
\bew Clear from the definition of  $c_k(\chi)$.
\bewende
\subsection{The main conjecture of Iwasawa theory} 
In this section we formulate the main conjecture of Iwasawa theory. 

Define 
\[
{\det}_{\Lambda}\Rh\Gamma_{\gl}(T_p(\chi)(k)):=
\bigotimes_{i=0}^2{\det}^{(-1)^i}_{\Lambda}\Hh^i_{\gl}(T_p(\chi)(k)).
\]
Denote by $Q(\Lambda)$ the total quotient ring of $\Lambda$.

\begin{prop}\label{identification}
Let $\chi$ be a Dirichlet character of conductor $N$. \\
a) If $\chi(-1)=(-1)^{k-1}$, then 
$\Hh^i_{\gl}(T_p(\chi)(k))\otimes_{\Lambda}Q(\Lambda)=0$ for $i=0,2$, $\Hh^1_{\gl}(T_p(\chi)(k))$ has $\Lambda$-rank $1$ and 
\[
\left(\Hh^1_{\gl}(T_p(\chi)(k))/\Lambda c_k(\chi)\right)\otimes_{\Lambda}Q(\Lambda)=0.
\]
Hence, 
\[
{\det}_{\Lambda}\Lambda c_k(\chi)\otimes_{\Lambda}Q(\Lambda)\isom
{\det}_{\Lambda}\Hh^1_{\gl}(T_p(\chi)(k))\otimes_{\Lambda}Q(\Lambda).
\]
In particular, this gives an isomorphism
\[
\Psi:{\det}_{\Lambda}\Lambda c_k(\chi)\otimes_{\Lambda}Q(\Lambda)\isom
{\det}^{-1}_{\Lambda}\Rh\Gamma_{\gl}(T_p(\chi)(k))\otimes_{\Lambda}Q(\Lambda).
\]
b) If $\chi(-1)=(-1)^{k}$, then for $i=0,1,2$
\[
\Hh^i_{\gl}(T_p(\chi)(k))\otimes_\Lambda Q(\Lambda)=0.
\]
In particular, one has   an isomorphism
\[
\Psi:Q(\Lambda)\isom
{\det}^{-1}_{\Lambda}\Rh\Gamma_{\gl}(T_p(\chi)(k))\otimes_\Lambda Q(\Lambda).
\]
\end{prop}
\bew 
The module $\Hh^0_{\gl}(T_p(\chi)(k))$ is 
zero. Hence, 
$\Hh^0_{\gl}(T_p(\chi)(k))\otimes_{\Lambda}Q(\Lambda)=0$. The result from \ref{vanishing} 
that $\Hh^2(T_p(\chi )^*(1-k))=0$ implies by proposition 1.3.2 from
\cite{Perrin-Riou-Asterisque} that $\Hh^2_{\gl}(T_p(\chi)(k))$ is a torsion $\Lambda$-module
and that $\Hh^1_{\gl}(T_p(\chi)(k))$ has $\Lambda$-rank $1$ if $\chi(-1)=(-1)^{k-1}$ and $\Lambda$-rank $0$,
if $\chi(-1)=(-1)^{k}$. 
It remains to show that $\Hh^1_{\gl}(T_p(\chi)(k))/\Lambda c_k(\chi)$ is
$\Lambda$-torsion. It suffices to consider $k<0$. As before $\Gamma=\Gal(\Q_\infty/\Q)$. It suffices to
to prove that
$\Hh^1_{\gl}(T_p(\chi)(k))_\Gamma/\Oh_p c_k(1)(\chi)$ is $\Oh_p$-torsion. 
$\Hh^1_{\gl}(T_p(\chi)(k))_\Gamma\subset H^1(\Z[1/p],T_p(\chi))$ has $\Oh_p$-rank $1$ by the formula for the Euler-Poincar\'e characteristic  and because $\Gamma$-invariants and $¸\Gamma$-coinvariants of $\Hh^2_{\gl}(T_p(\chi)(k))$ have the same rank. Finally by corollary \ref{recchi}, $c_k(\chi)\in \Hh^1_{\gl}(T_p(\chi)(k))_\Gamma$ is non-torsion.
\bewende

The main conjecture can now be formulated as follows:
\begin{mainconjecture}{\bf [Theorem \ref{mcthm}] (equal parity case)}\\
\label{mc1} 
Let $\chi$ be a character
with $\chi(-1)=(-1)^{k-1}$. The element $c_k(\chi)$ (\ref{iwasawac}) is 
mapped to a generator of the free $\Lambda$-module
\[
{\det}^{-1}_{\Lambda}\Rh\Gamma_{\gl}(T_p(\chi)(k))
\]
under the isomorphism $\Psi$ in proposition \ref{identification} a).
Equivalently,  $\Psi$ and the Poitou-Tate sequence  
induce an isomorphism 
\[
{\det}_{\Lambda}\left(\Rh\Gamma_{\loc}(T_p(\chi)(k))/\Lambda c_k(\chi)[-1]\right)
\isom
{\det}_{\Lambda}\left(\Rh\Gamma_{\gl}(T_p(\chi)^*(1-k))^*\right).
\]
\end{mainconjecture}

\rem  For $p\nmid \Phi(N)$ it
follows from the main conjecture as shown by Mazur and Wiles \cite{MaWi}. Under this
condition it was also proved directly by Rubin \cite{Ru}.\\

Let $\Gamma=\Gal(\Q_\infty / \Q)$. As before we decompose the cyclotomic character $\epsilon_{\cycl}=\epsilon\times\epsilon_\infty$. 
For every continuous character  $\tau:\Gamma\to \Oh_p^*$ define the twisting map 
\[
Tw(\tau ):\Lambda\to \Lambda
\]
by $\gamma\mapsto \tau (\gamma)\gamma$. 
Denote also by $Tw(\tau)$ the map induced on $Q(\Lambda)$.
The character $\tau$ also induces a map $\tau:\Lambda\to \Oh_p$.
If $\rho:\Lambda\to \Oh_p$ denotes the augmentation map, we have 
\[ \rho(Tw(\tau)(f))=\tau(f) \]
 for
any element $f\in\Lambda$. We extend the map $\tau$ to $Q(\Lambda)$ by $\tau(f/g)=\tau(f)/ \tau (g)$,
whenever $\tau(g)\neq 0$. Let $\chi$ be a Dirichlet character of conductor $N$ and $\chi(-1)=(-1)^k$.
Recall from \cite{Washington} theorem 7.10, that there is for $k>1$ an
element $\Lh_p(\chi, 1-k)\in Q(\Lambda)$ such that 
$\rho(\Lh_p(\chi, 1-k))=(1-\chi(p)p^{k-1})L(\chi,1-k)$. It is characterized by the
fact that for all characters $\tau:\Gamma\to \Oh_p^*$ of finite order
\[ \tau^{-1} (\Lh_p(\chi, 1-k))=(1-\chi\tau (p)p^{k-1})L(\chi\tau ,1-k)\ .\]

\rem The power series $f$ in \cite{Washington} theorem 7.10 is related to 
our $ \Lh_p(\chi, 1-k)$ as follows. Write $\chi\epsilon^{-k}=\theta\tau$ with $\theta$ of the first
and $\tau $ of the second kind. Then $Tw(\tau^{-1} \epsilon_\infty^{k-1})(f(t,\theta))=\Lh_p(\chi, 1-k)$.

\begin{defn}\label{p-Lfct} 
Let $\chi$ be a Dirichlet character with $\chi(-1)=(-1)^k$. For $k>1$ 
we call the above $\Lh_p(\chi, 1-k)\in Q(\Lambda)$ the {\em $p$-adic
$L$-function} at $k-1$. For $k\leq 1$ we define $\Lh_p(\chi, 1-k)\in Q(\Lambda)$  
as $Tw(\epsilon_\infty^{k-k'}) \Lh_p(\chi\epsilon^{k'-k}, 1-k')$ for some 
$k'>1$.
\end{defn}

\begin{mainconjecture}\label{mc2}{\bf [Theorem \ref{mc2true}] (unequal parity case)}\\
Let $p\neq 2$ and  $\chi$ be a Dirichlet character with $\chi(-1)=(-1)^{k}$. Then
the isomorphism $\psi$ from \ref{identification} b)
\[
Q(\Lambda)\isom {\det}^{-1}_{\Lambda}
\Rh\Gamma( T_p(\chi^{-1})(k))\otimes_\Lambda Q(\Lambda)
\]
maps the   $p$-adic $L$-function $\Lh_p(\chi,1-k)$
to a generator of ${\det}^{-1}_{\Lambda}\Rh\Gamma( T_p(\chi^{-1})(k))$.
\end{mainconjecture}

\rem  This result is already due to Mazur and Wiles \cite{MaWi}.

\subsection{Application of the Euler system to Iwasawa modules}
The general theory of Euler systems as developed by Kato \cite{KatoEuler},
Perrin-Riou \cite{PerrinEuler} and Rubin \cite{Ru} allows us to prove a result about
the determinants of certain Iwasawa modules.
We follow the formulation by Rubin.  

\begin{thm}\label{Hdivisibility} Let $\chi$ be a character with $\chi(-1)=(-1)^{k-1}$, then
\[
{\det}_{\Lambda}^{-1}\left(\Hh^1_{\gl}(T_p(\chi)(k))\left/\Lambda c_k(\chi)
\right)\right. \subset {\det}_{\Lambda}^{-1}\Hh^2_{\gl,0}(T_p(\chi)(k)).
\]
\end{thm}

\begin{rem} This should be formulated as follows: the element $c_k$ maps to
a generator of the invertible $\Lambda$-module 
\[
{\det}_{\Lambda}^{-1}\Hh^1_{\gl}(T_p(\chi)(k))\otimes {\det}_{\Lambda}\Hh^2_{\gl,0}(T_p(\chi)(k))
\]
in ${\det}_{\Lambda}^{-1}\Hh^1_{\gl}(T_p(\chi)(k))\otimes_{\Lambda}Q(\Lambda)$. 
Observe that we have $\Hh^2_{\gl,0}(T_p(\chi)(k))\otimes_{\Lambda}Q(\Lambda)=0$.
\end{rem}

\bew  By proposition \ref{identification} a)  
$\Hh^1_{\gl}(T_p(\chi)(k))$ has $\Lambda$-rank $1$.
We check the prerequisites for theorem 2.3.3. in \cite{Ru}.
These are called $\Hyp(K_{\infty}/K)$ and $\Hyp(K_{\infty}, T_p(\chi))$ in loc. cit.
As $K=\Q$ and $K_{\infty}=\Q_{\infty}$,
the hypothesis $\Hyp(K_{\infty}/K)$ is trivially verified. For 
$\Hyp(K_{\infty}, T_p(\chi))$ we take $\tau=\id$. 
Note that although our ${\Oh_p}$ is 
a product of discrete valuation rings, the theory of Rubin still goes through,
because we can apply it to every factor. The module called $X_{\infty}$ 
in theorem 2.3.3. in \cite{Ru}
 is our $\Hh^2_{\gl,0}(T_p(\chi)(k))$. Rubin defines
\[
\ind_{\Lambda}(c_k):=\left\{ \Phi(c_k): \Phi\in \Hom_{\Lambda}\left(\Hh^1_{\gl}(T_p(\chi)(k)),
\Lambda\right)\right\}\subset \Lambda.
\]
By the classification theory of $\Lambda$-modules there is a pseudo-isomorphism
\[
\rho:\Hh^1_{\gl}(T_p(\chi)(k))\to \Lambda\oplus \Hh^1_{\tors} 
\]
where $ \Hh^1_{\tors} $ is $\Lambda$-torsion. 
The index $\ind_{\Lambda}(c_k)=\ind_{\Lambda}(\rho(c_k))$ is given by the ideal $\det^{-1} \Lambda/ pr_1(\rho(c_k))\Lambda$.
There is an exact sequence
\[ 0\to K\to \left(\Lambda\oplus \Hh^1_{\tors}\right)/ \rho(c_k)\to \Lambda/pr_1(\rho(c_k))\to 0.\]
As $\det^{-1}K\subset \Lambda$ this 
implies  
\[
{\det}_{\Lambda}^{-1}\left(\Hh^1_{\gl}(T_p(\chi)(k))\left/\Lambda c_k(\chi)\right)\right.
\subset \ind_{\Lambda}(c_k).
\]
Kato's observation, \cite{Kato Main conj} proposition 6.1.,
shows that $\car(X_{\infty})={\det}_{\Lambda}^{-1}(X_{\infty})$.
Finally theorem 2.3.3 of Rubin in \cite{Ru} tells us that 
\[
\ind_{\Lambda}(c_k)\subset {\det}_{\Lambda}^{-1}(X_{\infty}).
\]
This gives the result stated in the theorem.
\bewende
The exact sequence in lemma \ref{classfieldseq} allows to reinterpret this: 
\begin{cor}\label{Ydivisibility} Let  $\chi$ be as in the theorem, then
\[
{\det}_{\Lambda}^{-1}\left(\Hh^1_{\loc}(T_p(\chi)(k))\left/\Lambda 
c_k(\chi)\right)\right.
\subset{\det}_{\Lambda}^{-1}\Hh^1_{\gl}( T_p(\chi)^*(1-k))^*
\]
holds.
\end{cor}
Our aim is to strengthen corollary \ref{Ydivisibility} to:
\begin{thm}\label{strongdivisibility} Let $\chi$ be a Dirichlet character of conductor $N$
with $\chi(-1)=(-1)^{k-1}$, then
\[
{\det}_{\Lambda}\left(\Rh\Gamma_{\loc}(T_p(\chi)(k))/\Lambda c_k(\chi)[-1]\right)\subset
{\det}_{\Lambda}\Rh\Gamma_{\gl}(T_p(\chi)^*(1-k))^*.
\]
Equivalently,
\[ \det_\Lambda \Lambda c_k(\chi)[-1]\supset \det_\Lambda\left(\Rh\Gamma_{\gl}(T_p(\chi)(k))\right) .
\]
\end{thm}
The proof of this theorem will be given at the end of this section.

\begin{rem} a) This is to be interpreted in the same way as theorem \ref{Hdivisibility}.\\
b) Note that by  
lemma \ref{vanishing} and the second exact sequence in lemma \ref{classfieldseq}, to prove the theorem it is enough to show that
\[
{\det}_{\Lambda}^{-1}\left(\Hh^1_{\loc}(T_p(\chi)(k))\left/\Lambda c_k(\chi)\right)\right.
\otimes{\det}_{\Lambda}\Hh^2_{\loc,0}(T_p(\chi)(k))\subset
{\det}_{\Lambda}^{-1}\Hh^1_{\gl}( T_p(\chi)^*(1-k))^*.
\]
Thus, it is necessary to study the module $\Hh^2_{\loc,0}(T_p(\chi)(k))$.\\
c) As will be shown by the explicit computation, the existence of the error term 
$\Hh^2_{\loc,0}(T_p(\chi)(k))$ is related to trivial zeroes of the $p$-adic
$L$-function. The problem already appears in Rubin's case $p\nmid \Phi(N)$, see
\cite{Ru} thm. 3.2.7. If $\chi$ is purely ramified, e.g., $N=p$ (the very first
case treated by Euler system methods), the error term vanishes because
it cancels against a global term.
\end{rem}\\

For a $\Lambda$-module $M$ denote by $M_{\pf}:=M\otimes_{\Lambda}\Lambda_{\pf}$the localization of $M$ at $\pf$. A Dirichlet character $\chi$ induces
 a finite character $\chi:\Gamma\to \Oh_p^*$, which extends to a map $\chi:\Lambda\to \Oh_p$ via $\gamma\mapsto \chi(\gamma)$.
Let $\af_\chi$ be the kernel of $\chi^{-1}:\Lambda\to\Oh_p$.  Then $\af_\chi$ is
a prime ideal of height $1$ and a principal ideal (generated by $\chi(\gamma_0)\gamma_0 -1$ where $\gamma_0$ is a topological generator of $\Gamma$).

\begin{prop}[see \ref{h2constant}, \ref{detcomputation}]\label{h2loccomputation}
Let $\chi$ be a Dirichlet character and $G_{\infty,p}:=\Gal(\Qbar_p/ \Q_\infty\otimes\Q_p)$,  then
\[
 \Hh^2_{\loc}(T_p(\chi)(1))\isom
\begin{cases}
 \text{finite} &\text{ if $\chi|_{G_{\infty,p}}\neq 1$,} \\
  T_p(\chi)& \text{ if $\chi|_{G_{\infty,p}}= 1$.} 
\end{cases}
\]
In particular, the localization at all prime ideals $\pf\neq \af_\chi$ of height $1$ is 
${\det}_{\Lambda_\pf} \Hh^2_{\loc}(T_p(\chi)(1))_\pf\isom \Lambda_\pf$.
\end{prop}


\begin{prop}[see \ref{deth1trivial}]\label{h1glcomputation} For all 
characters $\chi$ with 
$\chi(-1)=1$,
\[
\Hh^1_{\gl}( T_p(\chi)^*)^*_{\af_\chi}=0.
\]
In particular
${\det}_{\Lambda_{\af_\chi}}^{-1}\Hh^1_{\gl}( T_p(\chi)^*)^*_{\af_\chi}=\Lambda_{\af_\chi}$.
\end{prop}
\begin{prop}[see \ref{detrelation}]\label{h1loccomputation} Suppose that the conductor of $\chi$ is not 
a $p$-power and that $\chi(-1)=1$, then
\[
{\det}_{\Lambda_{\af_\chi}}^{-1}\left(\Hh^1_{\loc}(T_p(\chi)(1))/\Lambda c_1(\chi)\right) _{\af_\chi}
\subset \det_{\Lambda_{\af_\chi}}^{-1}\Hh^2_{\loc}(T_p(\chi)(1))_{\af_\chi}.
\]
\end{prop}
Before we can  prove theorem \ref{strongdivisibility}, we need one more observation:
\begin{lemma}\label{h2loc0} Suppose that the conductor of $\chi$ is not a $p$-power, then 
\[
{\det}_{\Lambda} \Hh^2_{\loc,0}(T_p(\chi)(1))\isom{\det}_{\Lambda} \Hh^2_{\loc}(T_p(\chi)(1)).
\]
If the conductor of $\chi$ is a $p$-power and $\chi(-1)=1$, then
\[
{\det}_\Lambda \Hh^2_{\loc,0}(T_p(\chi)(1))\isom \Lambda.
\]
\end{lemma} 
\bew If the conductor of $\chi$ is not a $p$-power, then $\chi$  is not trivial over $\Q_\infty$ 
and the map 
$\Hh^2_{\loc,0}(T_p(\chi)(1))\to \Hh^2_{\loc}(T_p(\chi)(1))$ is a pseudo-isomorphism,
because the cokernel is the finite group $\Hh^0_{\gl}(T_p(\chi)^*)^*$.
If the conductor of $\chi$ is a $p$-power, then $\chi$ is trivial over $\Q_\infty$ because of our
assumption $\chi(-1)=1$. Hence
$\Hh^0_{\gl}(T_p(\chi)^*)^*\isom T_p(\chi)$
and the map
\[
\Hh^2_{\loc}(T_p(\chi)(1))\to \Hh^0_{\gl}(T_p(\chi)^*)^*
\]
is the identity by 
\ref{h2loccomputation}.
 This implies that $\Hh^2_{\loc,0}(T_p(\chi)(1))=0$. 
\bewende

\bew{\em (of theorem \ref{strongdivisibility})}: First of all
it is enough to consider $k=1$. 
We use the remark after theorem \ref{strongdivisibility}. 
If the conductor of $\chi$ is a $p$-power, there is nothing to
show by lemma \ref{h2loc0}. From now on assume that the conductor of $\chi$
is not a $p$-power. Again by lemma \ref{h2loc0} we have to show that
\[
{\det}_{\Lambda}^{-1}\left(\Hh^1_{\loc}(T_p(\chi)(k))\left/\Lambda c_k(\chi)\right)\right.
\otimes{\det}_{\Lambda}\Hh^2_{\loc}(T_p(\chi)(k))\subset
{\det}_{\Lambda}^{-1}\Hh^1_{\gl}( T_p(\chi)^*(1-k))^*.
\]
As $\Lambda$ is a regular ring, the determinant of a $\Lambda$-module
is determined by its localizations at all primes of height $1$, see
\cite{Kato Main conj} proposition 6.1. Proposition \ref{h2loccomputation}
and corollary \ref{Ydivisibility} imply the theorem after localization at a prime
ideal of height $1$ different from $\af_\chi$. 
After localization at $\af_\chi$, the statement  follows from
\ref{h1glcomputation} and 
\ref{h1loccomputation}. 
The equivalence with the second statement follows from the Poitou-Tate
exact sequence.
\bewende

\subsection{The proof of the main conjecture}\label{equality}
In this section we reduce the main conjecture to the Bloch-Kato conjecture 
for $r=1$, i.e., the class number formula
proved in \ref{classnumber}.

\begin{thm} \label{mcthm}
Let $\chi$ be a character with $\chi(-1)=(-1)^{k-1}$, then
the main conjecture \ref{mc1} is true, i.e., the 
element $c_{k}(\chi)$ defined in \ref{iwasawac} maps to a generator of
\[
{\det}^{-1}_{\Lambda}\left(\Rh\Gamma_{\gl}(T_p(\chi)(k))\right).
\]
\end{thm}
The proof occupies the rest of this section. 
We  use that the main conjecture is invariant under twist. We will
prove it for $k=0$ and not for $k=1$. The reason is that  
we have to avoid zeroes 
 of the local $L$-factors.

Let $X$ be the set of characters $\chi$ with  $\chi(-1)=-1$
and conductor dividing some fixed $N$.
Given the inclusions of theorem \ref{strongdivisibility} for all $\chi\in X$, 
\[
{\det}_{\Lambda}\left(\Lambda c_0(\chi)[-1]\right)
\supset
{\det}_{\Lambda}\left(\Rh\Gamma_{\gl}(T_p(\chi))\right),
\]
it is enough to prove that the inclusion
\begin{equation*}
{\det}_{\Lambda}\left(\bigoplus_{\chi\in X}\Lambda c_0(\chi)[-1]\right)\supset
{\det}_{\Lambda}\left(\Rh\Gamma_{\gl}(\bigoplus_{\chi\in X}T_p(\chi))\right)
\end{equation*}
is an isomorphism.
This is a standard trick in Iwasawa theory. 
\begin{lemma} The above inclusion 
is an isomorphism if and only if it is an isomorphism after tensoring with
$\Lambda/\af$. 
\end{lemma}
\bew As $\Lambda$ is a local ring and $\af$ is contained in its radical, this 
is just Nakayama's lemma. 
\bewende
Note that $\af$ is the augmentation ideal so that the tensor product with 
$\Lambda/\af$ just means to take the coinvariants under $\Gamma=\Gal(\Q_\infty/ \Q)$, thus
\[
{\det}_{\Lambda}\Rh\Gamma_{\gl}\left(\bigoplus_{\chi\in X}T_p(\chi)\right)
\otimes^{\Lbb}_{\Lambda}\Lambda/\af\isom{\det}_{\Oh_p}
R\Gamma\left(\Z[1/p],\bigoplus_{\chi\in X}T_p(\chi)\right).
\]
On the other hand, we know from proposition \ref{identification} that 
$\bigoplus_{\chi\in X}\Lambda c_0(\chi)$ is a free 
$\Lambda$-submodule of $\Hh^1_{\gl}\left(\bigoplus_{\chi\in X}T_p(\chi)\right)$,
so that we have to show that 
\begin{equation*}
{\det}_{\Oh_p}\left(\bigoplus_{\chi\in X}\Oh_p c_0(1)(\chi)[-1]\right)\isom{\det}_{\Oh_p}
R\Gamma\left(\Z[1/p],\bigoplus_{\chi\in X}T_p(\chi)\right).
\end{equation*}
An $\Oh_p$-generator of the left hand side is 
$\prod_{\chi\in X}c_0(1)(\chi)$ and the claim is just the 
statement of corollary \ref{BlochKato}. 
This proves the main conjecture.

\section{Proof of the Bloch-Kato conjecture}
We now want to prove the Bloch-Kato conjecture. Note that there are four cases:
The character $\chi$ can have the same parity as $r$, this is the equal
parity case or the parities are different and we are in the unequal parity
case. Moreover $r$ can be $\geq 1$ or we can have $r\leq 0$. This makes four
cases. 
If we compare with the program laid out in section \ref{secoverview},
we see that  points 1 and 2 have meanwhile been settled.



\subsection{The Bloch-Kato conjecture in the equal parity case and $r\geq 1$}
\begin{thm}\label{Bloch-Katorge1} Let $\chi$ be a character with $\chi(-1)=(-1)^{r}$ with $r\geq 1$, $p\neq 2$, then the $p$-part of the
equivariant Bloch-Kato conjecture is true for $V(\chi)$ and $r$ and
the cyclotomic $\Z/p^n$-extension $\Q_n/\Q$.
\end{thm}
\bew
Let $G_n=\Gal(\Q_n/\Q)$.
We have by the main conjecture \ref{mcthm} for $\chi^{-1}$ with $k=1-r$
\[
{\det}_{\Lambda}\left(\Lambda c_{1-r}(\chi^{-1})[-1]\right)
\isom
{\det}_{\Lambda}\left(\Rh\Gamma_{\gl}(T_p(\chi^{-1})(1-r))\right)
\]
and if we tensor this equality over $\Lambda$ with $\Oh_p[G_n]$ and use lemma \ref{finitecoh}, we get 
that the element $c_{1-r}(\Q_n/\Q)(\chi^{-1})$ as defined before \ref{equivreform} maps to a generator of
\[
{\det}_{\Oh_p[G_n]}\left(R\Gamma(\Z_n[1/p],T_p(\chi^{-1})(1-r))\right)
\]
(compare \cite{Kato Main conj} lemma 6.3.).
By corollary \ref{equivreform} this implies the equivariant Bloch-Kato conjecture.
\bewende

Note that the theorem includes the absolute Bloch-Kato conjecture for $\chi$
up to powers of $2$.
Applying the functional equation (proposition \ref{localconjequiv}), we get from this:
\begin{cor}\label{Bloch-Katoung1}
Let $\chi$ be a character with $\chi(-1)=(-1)^{r-1}$ with $r<0$, $p\neq 2$, then the $p$-part of the
equivariant Bloch-Kato conjecture is true for $V(\chi)$ and $r$ and the
cyclotomic $\Z/p^n$-extension $\Q_n/\Q$.
\end{cor}
\begin{cor}\label{mc2true}
Let $p\neq 2$.
The Main Conjecture of Iwasawa theory holds in its second form
\ref{mc2}, i.e., for $\chi$ with $\chi(-1)=(-1)^{r-1}$ and $r\in \Z$,
the isomorphism \ref{identification} b)
\[
Q(\Lambda)\isom {\det}^{-1}_{\Lambda}
\Rh\Gamma_\gl( T_p(\chi^{-1})(1-r)){\otimes} Q(\Lambda)
\]
maps  $p$-adic $L$-function $\Lh_p(\chi,r)$ (see \ref{p-Lfct}) to a generator 
of 
${\det}^{-1}_{\Lambda}\Rh\Gamma( T_p(\chi^{-1})(1-r))$.
\end{cor}
\bew 
It suffices to consider one $r$. We take $r<0$.
Let $Q$ and $P$ be the characteristic power series of  
$\Hh^1_\gl(T_p(\chi^{-1})(1-r))$ and $\Hh^2_\gl(T_p(\chi^{-1})(1-r))$.
Then the lattice defined by ${\det}^{-1}_{\Lambda}\Rh\Gamma( T_p(\chi^{-1})(1-r))\subset Q(\Lambda)$ is generated by $P/Q$. On the other hand
$\Lh_p(\chi,r)= g/h$ with $g,h\in \Lambda$. We
have to show that 
\[ (Ph)=(Qg) \]
as ideals of $\Lambda$.

As above let $G_n=\Gal(\Q_n/\Q)$. Let $P_n,Q_n,g_n,h_n$ be the images of $P,Q,g,h$
in $\Oh_p[G_n]$. The modules $\Hh^i_\gl(T_p(\chi^{-1})(1-r))$ are
not only $\Lambda$-torsion but in addition 
$\Hh^i_\gl(T_p(\chi^{-1})(1-r)) \tensor \Oh_p[G_n]$ is finite because
all cohomology groups of 
$R\Gamma(\Z_n[1/p],T_p(\chi^{-1})(1-r))$ are finite. Hence 
$P_n,Q_n\in E_p[G_n]^*$ and the
lattice
defined by 
\[ \det^{-1}R\Gamma(\Z_n[1/p],T_p(\chi^{-1})(1-r))\isom 
\det^{-1}\Rh\Gamma_\gl( T_p(\chi^{-1})(1-r))\tensor \Oh_p[G_n] \]
in $E_p[G_n]$ is given by $P_n/Q_n$. By definition of the $p$-adic 
$L$-function \ref{p-Lfct} we have for all characters
$\tau$ of $\Gal(\Q_n/ \Q)$ the equality 
\[
\tau^{-1}(g_n/h_n)=(1-\tau\chi(p)p^r)L(\tau\chi,r)= L_{\{p\}}(\tau \chi,r),
\]
where $\tau:E_p[G_n]\to E_p$.
By definition of the equivariant $L$-function \ref{Lequidefn} with
the Euler factor at $p$ removed, this implies
$g_n/h_n=L_{\{p\}}(\Q_n/\Q,V(\chi),r)$,  Taking the Euler factor
from \ref{localfactor} into account, 
the equivariant Bloch-Kato conjecture in this case (\ref{Bloch-Katoung1})
implies
\[ (P_nh_n)=(g_nQ_n) \]
as ideals in $\Oh_p[G_n]$. Recall from \cite{Washington} proof of Theorem 7.1
that the kernel of $\Oh_p[[T]]\isom \Lambda\to \Oh_p[G_n]$ is contained in
$(p,T)^{n+1}$. Hence the equality of ideals in $\Oh_p[G_n]$ implies
\[ Ph\in \bigcap_{n}\left((gQ)+(p,T)^{n+1}\right)= (gQ)\]
because the $(p,T)$-adic topology on $\Lambda$ is separated. Conversely,
$gQ\in (Ph)$ and we have proved the claim.
\bewende

\subsection{The Bloch-Kato conjecture in the equal parity case and $r<0$}

We start by recalling a theorem of Beilinson \cite{Be} (as corrected by
Neukirch \cite{Ne} and Esnault \cite{Es}) which gives  elements in
motivic cohomology mapping to the polylogarithm in Deligne cohomology.
More precisely:
\begin{thm}\label{inftyregulator}(Beilinson \cite{Be}, Neukirch \cite{Ne}, Esnault \cite{Es}) There is
an element for $k> 1$
\[
b_k(\zeta_N)\in H^1_{\Mh}(\Z, V(\Q(\zeta_N))(k))
\]
such that 
\[
r_{\infty}(b_k(\zeta_N))=(-Li_k(\sigma\zeta_N))_{\sigma\in G}\in 
V_B(\Q(\zeta_N))(k-1)^+.
\]
\end{thm}
The next result computes the regulator of $b_k(\zeta_N)$ in 
$H^1(\Z[\zeta_N][1/p],E_p(k))$. This was first proved by Huber and 
Wildeshaus \cite{Hu-Wi} cor. 9.7 following  Beilinson and Deligne. A different approach to the theorem
is developed in \cite{Hu-Ki}.

\begin{thm}\label{pregulator} (\cite{Hu-Wi} cor. 9.7.) For $k>1$, the image of 
$b_k(\zeta_N)\in H^1_{\Mh}(\Z, V(\Q(\zeta_N))(k))$ under the map $r_p$ is
\[
\frac{1}{N^{k-1}(k-1)!}\ctilde_k(\zeta_N)\in H^1(\Z[\zeta_N][1/p],E_p(k)).
\]
\end{thm}
\bew This is corollary 9.7. in \cite{Hu-Wi} combined with lemma 
\ref{Kvergleich} and \cite{Hu-Wi} lemma B.4.9.
\bewende

\rem For $k=1$ the Beilinson element is $b_1(\zeta_N)=1-\zeta_N$. 
This is not an element of $\hm^1(\Z,V(\Q(\zeta_N)))=\Z[\zeta_N]^*\tensor E$
if $N=l^i$ is a prime power. However, even in this case, its $\chi$-component for $\chi\neq 1$
\[ b_1(\chi)= p_{\chi^{-1}}(1-\zeta_N)=p_{\chi^{-1}}\left( \frac{(1-\zeta_{l^i})^{l^{i-1}(l-1)}}{l} \tensor \frac{1}{l^{i-1}(l-1)}\right) \ .
\]
is in fact an element of $\hm^1(\Z,V(\chi)(1))$. The formula of the last
theorem holds with $\ctilde_1(\zeta_N)$ the class of $1-\zeta_N$ in Galois
cohomology. Again its $\chi$-component is an element of
$H^1(\Z[1/p],V_p(\chi)(1))$ for $\chi\neq 1$. This suffices for the
computations in the sequel.

\begin{thm}\label{Bloch-Katorle0} Let $r<0$ and $\chi(-1)=(-1)^{r}$
(or $r=0$, $\chi(-1)=1$ and $\chi(p)\neq 1$), 
then the 
Bloch-Kato conjecture for $V(\chi)$ and $r$ is true up to powers of $2$.
\end{thm} 
\bew 
We first repeat the computations used in the proof of the Beilinson conjecture
in this case.
The fundamental line for 
$V(\chi)$ and $r=1-k$ (see \ref{fundline}) is 
\[
\Delta_f(V(\chi)(1-k))={\det}_{E}^{-1}H^1_{\Mh}(\Z, V(\chi)^{\lor}(k))\otimes
{\det}_{E}^{-1}V_B(\chi)(1-k).
\]
We identitfy $V(\chi)^{\lor}\isom V(\chi^{-1})$. 
Let $t_B(\chi^{-1})$ be as in \ref{generatordefn}.
 Let $p_{\chi}$ be the projector onto the $\chi$-eigenspace  and let 
$b_k(\chi^{-1}):=p_{\chi}b_k(\zeta_N)$, where  $N$ is the conductor of $\chi$.
We claim that the element 
\[
\delta=\left(-\frac{N^{k-1}(k-1)!}{2}b_k(\chi^{-1})\right)^{-1}\otimes\;
(2\pi i)^{k-1}t_B(\chi^{-1})
\]
maps to $(L(\chi,1-k)^*)^{-1}$ under the isomorphism
\[
\Delta_f(V(\chi)(1-k))\otimes\R\isom E_{\infty}.
\]
Theorem \ref{inftyregulator} implies that $r_\infty$ maps  
$-\frac{N^{k-1}(k-1)!}{2}b_k(\chi^{-1})$ 
to 
\[
\frac{N^{k-1}(k-1)!}{2}
\sum_{\tau\in G}\chi^{-1}(\tau)Li_k(\tau \zeta_N)t_B(\chi^{-1}).
\]


Using the functional equation \ref{hurwitz} and the fact that the $\Gamma$-function has
residue $\frac{(-1)^{k-1} }{(k-1)!}$ at $s=1-k$, this is 
\[
L(\chi,1-k)^*(2\pi i)^{k-1}t_B(\chi^{-1}).
\]
Thus the Bloch-Kato conjecture with the lattice
$T_p(\chi^{-1})^\lor\subset V_p(\chi)$
says that $\delta$ is a generator of
${\det}_{\Oh_p} R\Gamma_c(\Z[1/p],T_p(\chi^{-1})^\lor(1-k))$. 
Using \ref{localfactor} for the local Euler factor 
(it applies as $(1-\chi(p)p^{k-1})\neq 0$ if $k\neq 1$ or $\chi(p)\neq 1$)
and \ref{ohnec}, this
is equivalent to
\[
(1-\chi(p)p^{k-1})r_p\left(-\frac{N^{k-1}(k-1)!}{2}b_k(\chi^{-1})\right)
\]
being a generator of 
${\det}^{-1}_{\Oh_p}R\Gamma(\Z[\zeta_N][1/p],T_p(\chi^{-1})(k))$
for all $p$.  With theorem \ref{pregulator}  and the 
relation $(1-\chi(p)p^{k-1})\ctilde_k(1)(\chi^{-1})=c_k(1)(\chi^{-1})$, this means that  
$\frac{-1}{2}c_k(1)(\chi^{-1})$ has to be a generator of 
${\det}_{\Oh_p}R\Gamma(\Z[1/p],T_p(\chi^{-1})(k))$.
The main conjecture implies that
\[
{\det}_{\Lambda}\left(\Lambda c_{k}(\chi^{-1})[-1]\right)
\isom
{\det}_{\Lambda}\left(\Rh\Gamma_{\gl}(T_p(\chi^{-1})(k))\right)
\]
so that 
\[
{\det}_{\Oh_p}\left(\Oh_p c_{k}(1)(\chi^{-1})[-1]\right)
\isom
{\det}_{\Oh_p}\left(R\Gamma_{\gl}(\Z[1/p],T_p(\chi^{-1})(k))\right),
\]
i.e., $c_{k}(1)(\chi^{-1})$ is a generator of 
${\det}^{-1}_{\Oh_p}R\Gamma(\Z[1/p],T_p(\chi^{-1})(k))$ as claimed.
\bewende

Applying the functional equation (proposition \ref{localconjequiv}), we get from this:
\begin{cor}\label{Bloch-Katoung2}
Let $\chi$ be a character with $\chi(-1)=(-1)^{r-1}$ with $r\geq 2$, then the
Bloch-Kato conjecture is true up to powers of $2$ for $V(\chi)$ and $r$.
\end{cor}
\subsection{The Bloch-Kato conjecture for $r=0,1$}
Using the main conjecture in its second form 
we can now prove the unequal parity case $\chi(-1)=(-1)^{r-1}$ and $r=0,1$.
For $r=0$ we still have  to assume $\chi(p)\neq 1$. 

\begin{thm}\label{r=0} Let $\chi$ be a character  with $\chi(-1)=-1$ and $\chi(p)\neq 1$.
 Then the Bloch-Kato conjecture for $V(\chi)$ and
$r=0$ is true up to powers of $2$. 
\end{thm}
\bew 
We have 
\[
\Delta_f(V(\chi))=E
\]
and we have to show that the element $(1-\chi(p))L(\chi,0)$ maps to a 
generator of 
\[ 
{\det}_{\Oh_p}^{-1}R\Gamma(\Z[1/p], T_p(\chi^{-1})(1)).
\]
Let $\rho$ be the augmentation map on $Q(\Lambda)$ as considered before \ref{p-Lfct}. Then 
$\rho(\Lh_p(\chi,k))=(1-\chi(p)p^{-k})L(\chi,k)$ for all $k\leq 0$.
By the main conjecture in the unequal parity case proved in corollary \ref{mc2true}, we have
that $\rho(\Lh_p(\chi,0))$ maps to a generator
of 
\[
{\det}_{\Lambda}^{-1}\Rh\Gamma( T_p(\chi^{-1})(1))\otimes_{\Lambda}\Oh_p
\isom 
{\det}_{\Oh_p}^{-1}R\Gamma(\Z[1/p], T_p(\chi^{-1})(1)).
\]
This proves the result.
\bewende

\rem If $\chi(p)=1$, then the main conjecture gives an isomorphism
between
$\det _{\Oh_p}0=\Oh_p$ and 
$\det _{\Oh_p}^{-1}R\Gamma(\Z[1/p], T_p(\chi^{-1})(1))$. This 
 does not tell us anything about
$L(\chi,0)^*$. In fact, the leading coefficients of $L(\chi,s)$ and
$(1-\chi(p)p^{-s})L(\chi,s)$ at $s=0$ now differ
 by a transcendental factor.

\begin{thm}\label{r=1} Let $\chi$ be a character  with $\chi(-1)=1$.
 Then the Bloch-Kato conjecture for $V(\chi)$ and
$r=1$ is true up to powers of $2$. 
\end{thm}

\bew If $\chi$ is the trivial character, then $V(\chi)=V(\Q)$ and 
the result follows from the class number formula \ref{classnumber}.
Suppose that $\chi\neq 1$. Let $N$ be the conductor of $\chi$.
We have 
\[
\Delta_f(V(\chi)(1))= {\det}_{E}^{-1}H^1_{\Mh}(\Z, V(\chi)(1))\tensor\det V_\DR(\chi) .
\]
Let
$b_1(\chi)=p_{\chi^{-1}}(1-\zeta_N)\in \hm^1(\Z,V(\chi)(1))$ as in the remark
after theorem \ref{pregulator}.
Similarly,
$\ctilde_1(1)(\chi)=p_{\chi^{-1}}(1-\zeta_N)\in H^1(\Z[1/p],V_p(\chi)(1))$.
Consider
\[
\delta=\left( -\frac{1}{2}b_1(\chi)\right)^{-1}\otimes p_{\chi^{-1}}(\zeta_N)
\in \Delta_f(V(\chi)).
\]
The regulator $r_\infty$ maps $-\frac{1}{2}b_1(\chi)$ to $1/2\sum_\sigma \chi(\sigma)Li_1(\sigma(\zeta_N))t_B(\chi)$, which is
$L(\chi,1)\tau(\chi^{-1})t_B(\chi)$, where $\tau(\chi^{-1}):=\sum_{a=1}^N\chi^{-1}(a)\zeta_N^a=\sum_\sigma\chi(\sigma)\sigma(\zeta_N)$ is the Gauss sum. 
The identification $ V_{\DR}(\chi)_\R\isom V_{B}(\chi)_\R$
maps $p_{\chi^{-1}}(\zeta_N)$ to $\tau(\chi^{-1})t_B({\chi})$, so that $\delta$
is the element which maps to $L(\chi,1)^{-1}$ under the isomorphism
$\Delta_f(V(\chi))\otimes_{\Q}{\R}\isom E\tensor\R$.
We have to compare its image with some natural integral structure of
\[ \det R\Gamma_c(\Z[1/p], V_p(\chi)(1))
=\det R\Gamma_f (\Z[1/p], V_p(\chi)(1))\tensor \det^{-1}R\Gamma_f(\Q_p, V_p(\chi)(1)) \]
under the identification \ref{globalident}.
Taking the Euler factor from \ref{localfactor} into account, the image of $\delta$ is
\[\left(-\frac{1-\chi(p)p^{-1}}{2}\ctilde_1(1)(\chi))\right)^{-1} \tensor \exp_p(p_{\chi^{-1}}(\zeta_N)) \]
where the first element is mapped to $H^1_f(\Z[1/p],V_p(\chi)(1))$ and the
second to $H^1_f(\Q_p,V_p(\chi)(1))$.
Note that we are in the case $H^i_c(\Z[1/p], V_p(\chi)(1))=0$ for all $i$,
hence we can also compare with the lattice
$\Oh_p\subset E_p=\det R\Gamma_c(\Z[1/p], V_p(\chi)(1))$.

The image of $-\frac{1-\chi(p)p^{-1}}{2}\ctilde_1(1)(\chi)$ under the composition
\[ H^1_f(\Z[1/p], V_p(\chi)(1))\to H^1_f(\Q_p, V_p(\chi)(1))\xrightarrow{\log_p} V_{\DR}(\chi)\otimes E_p\]
is 
\[ -\frac{1-\chi(p)p^{-1}}{2\# G} \sum_{\sigma}
\chi(\sigma)\log_p(1-\zeta_N^{\sigma})  
=\left( -\frac{1-\chi(p)p^{-1}}{2\tau(\chi^{-1})} \sum_{\sigma}
\chi(\sigma)\log_p(1-\zeta_N^{\sigma})  \right) p_{\chi^{-1}}(\zeta_N)
\]
The conjecture allows the choice of lattice in $V_p(\chi)$ and we use
$T_p(\chi^{-1})^\lor$.
Hence it remains to show that the element
\[  -\frac{1-\chi(p)p^{-1}}{2\tau(\chi^{-1})} \sum_{\sigma}
\chi(\sigma)\log_p(1-\zeta_N^{\sigma}) \]
of $E_p$ is a generator of 
$\det R\Gamma_c(\Z[1/p],T_p(\chi^{-1})^\lor)(1))\isom\det R\Gamma(\Z[1/p],T_p(\chi^{-1}))$.
By the main conjecture \ref{mc2true}  it is indeed generated by
$\rho(\Lh_p(\chi,1))$ where $\rho$ is the augmentation of $Q(\Lambda)$ as considered before \ref{p-Lfct}. The formula from \cite{Washington} theorem 5.18 together with  our 
normalization of the isomorphism $\Gal(\Q(\zeta_N)/\Q)$ with $(\Z/N)^*$ gives
\[  \rho(\Lh_p(\chi,1))=-(1-\chi(p)p^{-1})\frac{\tau(\chi)}{N}\sum_{\sigma}
\chi(\sigma)\log_p(1-\zeta_N^{\sigma}).
\]
The claim follows using the relation $\tau(\chi)\tau(\chi^{-1})=N$.
\bewende

\subsection{End of proof of the Bloch-Kato conjecture}

\begin{thm}\label{Bloch-Katofcteq} The Bloch-Kato conjecture is true for all
abelian Artin motives. 
\end{thm}
\bew 
Extension of coefficients is faithfully flat, hence it suffices to prove
the conjecture after extension of coefficients.
Every abelian Artin motive is the direct sum of motives of the
form $V(\chi)$ for Dirichlet characters $\chi$ after an appropriate extension, so it suffices to prove the
theorem for these. Collecting our results in
theorem \ref{Bloch-Katorge1}, corollary \ref{Bloch-Katoung1}, theorem \ref{Bloch-Katorle0}, corollary \ref{Bloch-Katoung2}, theorem \ref{r=0} and theorem \ref{r=1}, we see that only the case $r=0$ and $\chi(p)=1$ is still missing.
We get the final case by proposition \ref{B.1.2} via the functional equation.
\bewende

\section{Investigation of some Iwasawa modules}\label{investigation}
In this section we prove some structural results about Iwasawa modules,
which are 
only needed for theorem \ref{strongdivisibility} and 
thus can be skipped at a first reading. For the definitions of the 
$\Lambda$-modules we refer to \ref{fettH}.

\subsection{The Iwasawa module $\Hh^2_{\loc}(T_p(\chi)(1))$}
In this section we investigate 
$\Hh^2_{\loc}(T_p(\chi)(1))$, which is necessary for the proof of 
proposition \ref{h2loccomputation}.

\begin{lemma}\label{h2constant} Let $G_{\infty,p}:=\Gal(\Qbar_p/ \Q_\infty\otimes\Q_p)$, then
\[
\Hh^2_{\loc}(T_p(\chi)(1))=
\begin{cases} 
    \text{finite} & \text{if $ \chi|_{G_{\infty,p}}\neq 1$,}\\
     T_p(\chi) &\text{ if $\chi|_{G_{\infty,p}}= 1$}. 
\end{cases}
\]
Moreover there exists $n$ such that $\Hh^2_{\loc}(T_p(\chi)(1))\isom H^2(\Q_n\otimes\Q_p,T_p(\chi)(1))$.
\end{lemma}
\bew 
As $\Hh^2_{\loc}(T_p(\chi)(1))\isom\Hh^0_{\loc}(T_p(\chi)^*)^*$ we get that 
$\Hh^2_{\loc}(T_p(\chi)(1))$ are the coinvariants $T_p(\chi)_{G_{\infty ,p}}$ and
the statement of the lemma is clear.
\bewende

\begin{cor} \label{detcomputation} 
Let $\chi$ be a Dirichlet character and $G_{\infty,p}$ as above, then
\[
{\det}_\Lambda \Hh^2_{\loc}(T_p(\chi)(1))\isom
\left\{ \begin{array}{ll}\Lambda &\mbox{ if }\chi|_{G_{\infty,p}}\neq 1 
\\
{\af_\chi^{-1}}& \mbox{ if }\chi|_{G_{\infty,p}}= 1, 
\end{array}\right.
\]
where $\af_\chi$ is the kernel of $\chi^{-1}:\Lambda\to \Oh_p$. 
In particular, the localization
${\det}_{\Lambda_\pf} \Hh^2_{\loc}(T_p(\chi)(1))_\pf\isom \Lambda_\pf$ for all prime ideals $\pf\neq \af_\chi$ of height $1$.
\end{cor}
\bew Consider the computation in lemma \ref{h2constant}. In the first case,
the corollary 
follows from the fact that the determinant of a finite $\Lambda$ module is trivial. In the second case, it follows from the fact that $\chi(\gamma_0)\gamma_0-1$ generates
$\af_\chi$ (where $\gamma_0$ generates $\Gamma$) and 
$T_p(\chi)\isom \Lambda/(\chi(\gamma_0)\gamma_0-1)$ as $\Lambda$-modules.
\bewende


\subsection{The Iwasawa module $\Hh^1_{\gl}(T_p(\chi)^*)^*$} 
In this section we prove \ref{h1glcomputation}.
Let $\af_\chi$ be the kernel of $\chi:\Lambda\to \Oh_p$.

\begin{prop}\label{deth1trivial} Let $\chi$ be a character with $\chi(-1)=1$, then
the module 
\[
\left(\Hh^1_{\gl}(T_p(\chi)^*)^*\right)_{\af_\chi}=0
\]
is zero. 
\end{prop}
\bew Let $\chi=\theta\tau$ where $\theta$ is of the first and $\tau$ of the second kind. $\tau$ can be viewed as a finite character  of $\Gamma=\Gal(\Q_\infty/\Q)$. 
The statement remains
invariant under twisting with $\tau^{-1}$ (which is trivial over $\Q_\infty$), hence it suffices to consider
the case $\tau=1$, i.e., $\af_\chi=\af$ is the augmentation ideal.

Let $\kappa(\af) $ be the residue class field of $\af$.
By Nakayama's lemma it is enough to show
that
$\kappa(\af)\otimes_{\Lambda_{\af}}(\Hh^1_{\gl}( T_p(\chi)^*)^*)_{\af}=0$.

As $\Lambda/\af={\Oh_p}$  it suffices to show that 
\[
\Lambda/\af\otimes_{\Lambda}(\Hh^1_{\gl}( T_p(\chi)^*)^*)
\]
is finite.  By the definition of 
$\Lambda={\Oh_p}[[\Gamma]]$ we have to show that the
$\Gamma$-invariants
\[
\left(\left(\Hh^1_{\gl}({T_p}(\chi)^*)\right)^{\Gamma}\right)^*
\]
are finite. The Hochschild-Serre spectral sequence gives
\[
0\to H^1(\Gamma, \Hh^0_{\gl}( {T_p(\chi)}^*))\to 
H^1(\Z[1/p],{T_p(\chi)}^*)\to
\Hh^1_{\gl}( {T_p(\chi)}^*)^{\Gamma}\to 0,
\]
because $\Gamma\isom\Z_p$ has cohomological dimension $1$. The group
in the middle 
is finitely cogenerated. It suffices to compute the corank
of the cokernel.
We compute the rank of $H^1(\Z[1/p],{T_p(\chi)}^\lor)$ via the
Euler-Poincar{\'e} characteristic
\[
\sum_{i=0}^2(-1)^i\rk_{E_p}H^i(\Z[1/p],V_p(\chi)^{\lor})
\]
which is equal to $0$ by \cite{Jannsen}, lemma 2, as $\chi(-1)=1$. 
We know by \cite{Schneider} 
paragraph 7, Satz 2 that $H^2(\Z[\zeta_{N}][1/p],{E_p})^+$ (where $+$ denotes the invariants of 
complex conjugation)
is zero because Leopoldt's conjecture is true for $\Q(\zeta_{N})$.
This implies that $H^2(\Z[1/p],{V_p(\chi)^{\lor}})$ is zero,
because it is a direct summand. 

If  $\chi\neq 1$ is not the trivial character
$H^0(\Z[1/p],V_p(\chi)^{\lor})=0$  
 and hence $\rk_{E_p}H^1(\Z[1/p],V_p(\chi)^{\lor})=0$
in this case. Hence the corank of 
$H^1(\Z[1/p],T_p(\chi)^{*})$ is also $0$.

If $\chi=1$ is the trivial character, the same argument shows that
$H^1(\Z[1/p],T_p(\chi)^{*})$ has corank $1$.
However, in this case
$H^1(\Gamma, \Hh^0_{\gl}( {T_p(\chi)}^*))$ also has corank $1$.
\bewende

\subsection{The Iwasawa module $\Hh^1_{\loc}(T_p(\chi)(1))$}
We  show proposition \ref{h1loccomputation}. As before let $\af_\chi$ be the kernel of $\chi^{-1}:\Lambda\to \Oh_p$.
Recall from section \ref{artinmotives} the definition of 
$H^1_{/f}(\Q_p\otimes\Q_n, V_p(\chi)(1))$
and let
\[
\Hh^1_{/f}(V_p(\chi)(1)):=\prolim_nH^1_{/f}(\Q_p\otimes\Q_n, V_p(\chi)(1)).
\]
We will later  consider the  canonical map
\[
\Hh^1_{\loc}(T_p(\chi)(1))\to \Hh^1_{/f}(V_p(\chi)(1)).
\]
the next lemma shows that the image of this map is a finitely generated
${\Oh_p}$-module and computes its determinant.
\begin{lemma}\label{h/fdet} The $\Lambda\otimes_{{\Oh_p}}E_p$-module 
$\Hh^1_{/f}(V_p(\chi)(1))$ is
a free module of finite rank over $E_p$. Moreover 
\[
{\det}_{\Lambda_{\af_\chi}}\Hh^1_{/f}(V_p(\chi)(1))_{\af_\chi}\isom 
{\det}_{\Lambda_{\af_\chi}}\Hh^2_{\loc}({T_p(\chi)}(1))_{\af_\chi}
\]
\end{lemma}
\bew We show first that the inverse limit defining $\Hh^1_{/f}(V_p(\chi)(1))$
is in fact constant. By duality
\[
H^1_{f}(\Q_p\otimes\Q_n,V_p(\chi^{-1}))\isom
H^1_{/f}(\Q_p\otimes\Q_n,V_p(\chi)(1))^{\lor}
\]
and the exact sequence in section \ref{artinmotives} 
\begin{multline*}
 0\to H^0(\Q_p\otimes\Q_n,V_p(\chi^{-1}))\to
D_\cris(V(\Q_n)\otimes V_p(\chi^{-1}))\\
\to D_\cris(V(\Q_n)\otimes V_p(\chi^{-1}))\to H^1_{f}(\Q_p\otimes\Q_n,V_p(\chi^{-1}))
\to 0
\end{multline*}
implies that 
$\rk_{E_p}H^1_{f}(\Q_p\otimes\Q_n,V_p(\chi^{-1}))=
\rk_{E_p} H^0(\Q_p\otimes\Q_n,V_p(\chi^{-1}))$
and the latter is constant for $n\gg 0$,
 because $\Q_n$ is totally ramified at $p$. This implies that 
$\Hh^1_{/f}(V_p(\chi)(1))$ is a free module of finite rank over $E_p$.
We compute the determinants using the above exact sequence for $n$ big enough: 
\begin{align*}
{\det}_{\Lambda_{\af_\chi}}\Hh^1_{/f}(V_p(\chi)(1))_{\af_\chi}&\isom
{\det}_{\Lambda_{\af_\chi}}H^1_{/f}(\Q_p\otimes\Q_n,V_p(\chi)(1))_{\af_\chi}\\
&\isom {\det}_{\Lambda_{\af_\chi}}^{-1}H^1_{f}(\Q_p\otimes\Q_n,V_p(\chi^{-1}))_{\af_\chi}\\
&\isom {\det}_{\Lambda_{\af_\chi}}^{-1}H^0(\Q_p\otimes\Q_n,V_p(\chi^{-1}))_{\af_\chi}\\
 &\isom {\det}_{\Lambda_{\af_\chi}}H^2(\Q_p\otimes\Q_n,V_p(\chi)(1))_{\af_\chi}.
\end{align*}
Note that $p\not\in \af_\chi$ so that the map $\Lambda\to \Lambda_{\af_\chi}$
factors through $\Lambda[1/p]=\Lambda\otimes_{{\Oh_p}}E_p$. In particular
\[
{\det}_{\Lambda_{\af_\chi}}H^2(\Q_p\otimes\Q_n,V_p(\chi)(1))_{\af_\chi}\isom
{\det}_{\Lambda_{\af_\chi}}H^2(\Q_p\otimes\Q_n,{T_p(\chi)}(1))_{\af_\chi}
\]
and by lemma \ref{h2constant} the last determinant is equal to 
\[
{\det}_{\Lambda_{\af_\chi}}\Hh^2_{\loc}({T_p(\chi)}(1))_{\af_\chi}.
\]
This is the desired result.
\bewende
\begin{lemma} Suppose that $N$ is not a $p$-power, then the element 
$c_1(\chi)\in\Hh^1_{\loc}({T_p(\chi)}(1)) $ maps to zero
in $\Hh^1_{/f}(V_p(\chi)(1))$. 
\end{lemma}
\bew By appendix \ref{applocalcomp} the map at finite level 
\[
H^1(\Z_n[\zeta_{p^nN}][1/p], \Oh_p(1))\to H^1(\Q_p\otimes\Q(\zeta_{p^nN}), \Oh_p(1))
\to H^1_{/f}( \Q_p\otimes\Q(\zeta_{p^nN}), E_p(1))
\]
maps a unit in 
$\Z[\zeta_{p^nN}][1/p]^{\mal}\otimes\Oh_p\isom H^1(\Z[\zeta_{p^nN}][1/p], \Oh_p(1))$
to the system of its valuations in 
$\bigoplus_{v|p}E_p\isom H^1_{/f}( \Q_p\otimes\Q(\zeta_{p^nN}), E_p(1))$.
As $N$ is not a $p$-power the elements $c_1(\zeta_{p^nN})=1-\zeta_{p^nN}$ are units
in $\Z[\zeta_{p^nN}]$, hence have zero valuation. Using the commutative
diagram
\[
\begin{CD} 
\prolim_nH^1(\Q_p\otimes\Q(\zeta_{p^nN}),\Oh_p(1))@>>>\prolim_nH_{/f}^1(\Q_p\otimes\Q(\zeta_{p^nN}),E_p(1))\\
@VVV@VVV\\
\Hh^1_{\loc}({T_p(\chi)}(1))@>>>\Hh^1_{/f}(V_p(\chi)(1))
\end{CD}
\] 
proves the claim.
\bewende
As $p$ is invertible in $\Lambda_{\af_\chi}$, we get a surjection
\[
\left(  \Hh^1_{\loc}(T_p(\chi)(1))/\Lambda c_1(\chi)\right)_{\af_\chi}\to 
\Hh^1_{/f}(V_p(\chi)(1))_{\af_\chi}.
\]
\begin{cor}\label{detrelation} Suppose that $N$ is not a $p$-power and $\chi(-1)=1$,
then this map induces 
\[
{\det}_{\Lambda_{\af_\chi}}^{-1}
\left(\Hh^1_{\loc}(T_p(\chi)(1))/\Lambda c_1(\chi)\right)_{\af_\chi}\subset
{\det}_{\Lambda_{\af_\chi}}^{-1}\Hh^2_{\loc}(T_p(\chi)(1))_{\af_\chi}.
\]
\end{cor} 
\bew Let 
\[
K:=\ker\left(\left(\Hh^1_{\loc}(T_p(\chi)(1))/\Lambda c_1(\chi)\right)_{\af_\chi}\to 
\Hh^1_{/f}(V_p(\chi)(1))_{\af_\chi}\right).
\]
Then 
\[
{\det}_{\Lambda_{\af_\chi}}^{-1}K\otimes{\det}_{\Lambda_{\af_\chi}}^{-1}\Hh_{/f}^1(V_p(\chi)(1))_{\af_\chi}\isom {\det}_{\Lambda_{\af_\chi}}^{-1}
\left(\Hh^1_{\loc}(T_p(\chi)(1))/\Lambda c_1(\chi)\right)_{\af_\chi}
\] 
and ${\det}_{\Lambda_{\af_\chi}}^{-1}K$ is an ideal in $\Lambda_{\af_\chi}$
because $\Lambda_{\af_\chi }$ is a discrete valuation ring. 
This together with lemma \ref{h/fdet} implies the result.
\bewende

\begin{appendix}

\section{A local computation}\label{applocalcomp}
The following lemmas are known to the experts but for lack of
finding a precise reference we repeat the proofs here. They provide
the proof of proposition \ref{flocalization}.

We keep the notations of
the main text.  Let $K$ be a finite extension of $\Q_p$.
 Following Bloch-Kato we define 
$H^1_f(K, \Z_p(1)):=\iota^{-1}(H^1_f(K, \Q_p(1))$, where 
$\iota:H^1(K, \Z_p(1))\to H^1(K, \Q_p(1))$, and $H^1_{/f}:=H^1/H^1_f$. For an
abelian group $A$ we denote by 
$A^\land := \prolim_nA\otimes_\Z\Z/p^n\Z$ the $p$-adic completion of $A$.

\begin{alemma}
Let $K$ be a finite extension of $\Q_p$. Then
\[ H^1_f(K,\Z_p(1))= (\Oh_K^*)^\land  \]
under the identification $H^1(K,\Z_p(1))=(K^*)^\land $.
\end{alemma}
\bew
Bloch-Kato show in \cite{Bloch-Kato} p. 358 that their exponential
agrees with the usual $p$-adic exponential in this case. Hence
$H^1_f(K,\Q_p(1))\subset (K^*)^\land \otimes \Q_p$ agrees with the
image of the exponential which is $(\Oh_K^*)^\land\tensor \Q_p$.
The preimage of this group in $(K^*)^\land $ is
$(\Oh_K^*)^\land$.
\bewende

\begin{acor}
There is a unique isomorphism 
\[ H^1_{/f}(K,\Z_p(1)) \to \Z_p \]
compatible with 
\[ H^1(K,\Z_p(1))=(K^*)^\land \xrightarrow{v} \Z_p \]
where $v$ is the valuation map of $K$ normalized such that $v(\pi)=1$
for a uniformizer.
\end{acor}
\bew
By definition, $H^1_{/f}(K,\Z_p(1))$ is the image of 
$H^1(K,\Z_p(1))=(K^*)^\land $ in 
$H^1(K,\Q_p(1))/H^1_f(K,\Q_p(1))$. Hence it can be identified
with \hbox{$(K^*)^\land /(\Oh_K^*)^\land$} which is isomorphic
to $\Z_p$ via the valuation.
\bewende
\begin{aprop}\label{anhangPT}
Let $F$ be a finite extension of $\Q$.
The Poitou-Tate localization sequence induces an exact sequence
\begin{multline*}
0\to(\Oh_K^*)^\land \to H^1(\Oh_F[1/p],\Z_p(1))\to 
H^1_{/f}(F\otimes\Q_p,\Z_p(1)) 
\to Cl(\Oh_F)^\land \to \\
\to H^2(\Oh_F[1/p],\Z_p(1))\to H^2(F\otimes\Q_p,\Z_p(1))\to
H^0(\Oh_F[1/p],\Q_p/\Z_p)^*\to 0.
\end{multline*}
\end{aprop}
\bew 
We apply the previous corollary to the local fields in $F\tensor\Q_p$. In
particular, we identify $H^1_{/f}(F\tensor\Q_p,\Q_p(1))\isom \bigoplus_{v\mid p} \Q_p$ via the valuation maps. 
Thus the kernel of $H^1(\Oh_F[1/p], \Z_p(1))\to H^1_{/f}(F\otimes\Q_p, \Q_p(1))$ is
$(\Oh_F^*)^\land $. On the other hand it follows from \cite{Schneider} that one has 
an exact sequence
\begin{multline*}
0\to Cl(\Oh_F[1/p])^\land \to H^2(\Oh_F[1/p], \Z_p(1))\to\\ 
\to H^0(F\tensor\Q_p,\Q_p/\Z_p)^*
\to  H^0(\Oh_F[1/p],\Q_p/\Z_p)^*        \to 0.
\end{multline*}
Let $j:\spec \Oh_F[1/p]\to \spec\Oh_F$.
The long exact sequence in Zariski-cohomology for the short exact
sequence of sheaves
\[0\to \Oh^*\to j_*\Oh^*\to\bigoplus_{v\mid p}\Z\to 0 \]
yields
\[ 0\to \Oh_F^*\to \Oh_F[1/p]^*\to 
\bigoplus_{v\mid p}\Z
\to Cl(\Oh_F)\to 
Cl(\Oh_F[1/p])\to 0
\ .\]
If we finally observe that 
$H^1(\Oh_F[1/p], \Z_p(1))\isom (\Oh_F[1/p]^*)^\land $ we get the desired result.
\bewende

Let $K_0$ be the maximal unramified
subfield of $K$. Let $\F_q$ be its residue field, $q=p^f$. 
We have $D_{\cris}(\Ind_{\Q_p}^{K}\Q_p)=B_{\cris}^{G_{K}}=K_0$.
As before $\phi$ is the Frobenius on $B_\cris$ normalized as in \cite{FontBou}
3.2.
The fundamental short exact sequence
\[ 0\to \Q_p\to B_{\cris}\cap B_{dR}^+\xrightarrow{\phi-1} B_{\cris}\to 0,\]
induces
\[ 0\to \Q_p\to K_0\xrightarrow{\Fr-1}K_0\to H^1_f(K,\Q_p)\to 0\]
where $\Fr$ is the geometric Frobenius in $G(K_0/\Q_p)=G(\F_q/\F_p)$
by loc. cit. footnote p.\nobreak\  210.

\begin{alemma}\label{anhangexactseq}
The following sequence is exact:
\[ 0\to \Z_p\to \Oh_{K_0}\xrightarrow{\Fr-1}\Oh_{K_0}\xrightarrow{\tr}\Z_p
\to 0 \]
\end{alemma}
\bew
 We take the point of view of $\Z_p[G]$-modules
with $G=G(K_0/\Q_p)$.
Then $\Oh_{K_0}$ is isomorphic to $\Z_p[G]$ by choice of
a primitive element $\omega$. Note that this is possible because it is possible for finite fields and $K_0$ is unramified. We have a commutative diagram
\[\begin{CD}
0@>>> \Z_p@>\Delta >> \Z_p[G]@>\Fr-1 >>\Z_p[G]@>\Sigma >> \Z_p@>>> 0\\
@III @V\tr(\omega) VV @V\omega VV @V\omega VV @VV\tr(\omega) V \\
0@>>> \Z_p @>>>\Oh_{K_0}@>\Fr-1 >> \Oh_{K_0}@>\tr >>\Z_p@>>> 0
\end{CD}\]
where $\Delta$ is the diagonal embedding, $\Sigma$ is summation,
$\omega$ means muliplication by $\omega$ on the right.
The first line is exact as the
image of $(\Fr-1)$ is nothing but the augmentation
ideal. In the second line, $ \Z_p$ is kernel of $\Fr-1$. Hence
$\tr(\omega)$ is a unit of $\Z_p$. All vertical maps are isomorphisms,
hence the second line is also exact as claimed.
\bewende

In particular:
\begin{acor} The boundary morphism
$\delta:K_0\to H^1_f(K,\Q_p)$ factors via the trace $\delta=\delta'\circ\tr$ and
$\delta'\tr(\Oh_{K_0})$ is a natural lattice in $H^1_f(K,\Q_p)$. 
\end{acor} 

We do not need to know
whether this lattice agrees with $H^1_f(K,\Z_p)$, but rather:

\begin{alemma}\label{anhanglocaldual}
Under the local duality isomorphism 
$H^1_{/f}(K,\Q_p(1))^\lor\isom H^1_f(K,\Q_p)$ the lattice
$H^1_{/f}(K,\Z_p(1))^\lor$ is identified with $\delta'\tr(\Oh_{K_0})$.
\end{alemma}
\bew
Recall that $H^1(K,\Q_p(1))\isom (K^*)^\land \tensor \Q_p$ and
$H^1(K,\Q_p)=\Hom(G_K,\Q_p)$.
We have identified $H^1_{/f}(K,\Z_p(1))$ with $\Z_p$ via
the valuation map, i.e., a generator of $H^1_{/f}(K,\Z_p(1))^{\lor}$
is the valuation map.
 Local duality maps the valuation map
to the element $\phi_{K_0}^\lor\in \Hom(G_K,\Q_p)$
which factors through $G_{K_0}$ and maps the
arithmetic Frobenius of $K_0$  to $1$.  
Hence we have to show that $\phi_{K_0}^\lor$
generates $\tr(\Oh_{K_0})\subset H^1(K,\Q_p)$.
As  
$H^1_f$ and $H^1_{/f}$ are orthogonal under local duality, this is
certainly true rationally.

We first consider the special case $K=K_0=\Q_p$. Then $\Fr$ operates as
the identity.
In order to compute the $\delta(1)$ 
 it suffices to
construct a preimage under $\phi-1$ in the Wittring 
$W(\bar{\F_p})\subset B_{\cris}\cap B_{dR}^+$ and note that 
on this subring the operation of the Frobenius of the Galois group and
of the Frobenius in $B_{\cris}$ agrees. The computation yields 
$\delta(1)=\phi_{\Q_p}^\lor$.

In general, we consider the commutative diagram
\[ \begin{CD}
\Oh_{K_0}@>\tr >>\Z_p@>\delta'>> H^1(G_{K_0},\Q_p)\\
@AAA @AAfA @AAA\\
\Z_p@>id>> \Z_p@>\delta>>H^1(G_{\Q_p},\Q_p)
\end{CD}
\]
and note that the image of $\phi_{\Q_p}^\lor$ in $H^1(G_{K_0},\Q_p)$
is $f\phi_{K_0}^\lor$.
\bewende
\rem This computation is prone to normalization problems between arithmetic
and geometric Frobenius. However, this is irrelevant for the lattice
statements we need.

\begin{aprop}\label{anhanglocallattices}
Let $K$ be a finite extension of $\Q_p$. Under the composition
\[ 
\det^{-1}R\Gamma_{/f}(K,\Q_p(1))
\isom \det R\Gamma_{f}(K,\Q_p)
\xrightarrow{\alpha} E_p
\]
of local duality with the isomorphism of \ref{localident} the lattice $\det^{-1}R\Gamma_{/f}(K,\Z_p(1))$ is identified
with $\Z_p$. 
\end{aprop} 
\bew 
The lattice on the left is
\[ 
\det^{-1}R\Gamma_{/f}(K,\Z_p(1))= 
\det H^1_{/f}(K,\Z_p(1)) \tensor \det^{-1} H^2(K,\Z_p(1))
\]
Local duality identifies $H^2(K,\Z_p(1))\isom H^0(K,\Q_p/\Z_p)^* =
H^0(K,\Z_p)^\lor$. The lattice in $H^1_{/f}$ was identified in lemma
\ref{anhanglocaldual}. Hence
$\det^{-1}R\Gamma_{/f}(K,\Z_p(1))$ is identified with
\[ 
\det^{-1}\delta'\tr(\Oh_{K_0})\tensor \det H^0(K,\Z_p)
\stackrel{\ref{anhangexactseq}}{=}\det{\Oh_{K_0}}\tensor \det^{-1}{\Oh_{K_0}}
\]
where we use $\Oh_{K_0}\xrightarrow{\Fr-1} \Oh_{K_0}$ as integral structure of
$R\Gamma_f(K,\Z_p)$.
By definition this integral structure is mapped to 
$\Z_p$ by $\alpha$.
\bewende

\section{The functional equation}\label{appfcteq}
The aim of this appendix is to prove Conjecture \ref{localconj}, 
(compatibility of the Bloch-Kato conjecture with the functional equation) in
the case of Dirichlet motives and for $r>1$. See the remarks after \ref{localconj} for
the cases which are in the literature. 

\subsection{Global considerations}

Let $\chi$ be a Dirichlet character of conductor $N$ and $r\geq 1$, 
 $\Oh$  a finite extension of $\Z$ containing all values of $\chi$, 
$\zeta_N=\exp(2\pi i/N)$ and $G$ the Galois group of $\Q(\zeta_N)$
over $\Q$. In particular, we have distinguished an embedding $\sigma_0:\Q(\mu_N)\to \C$.
Recall from section \ref{artinmotives} the projector 
\[ p_{\chi^{-1}}: V(\Q(\zeta_N))\to V(\Q(\zeta_N)) \]
whose image is $V(\chi)$.  Recall also that elements of 
$V_B(\Q(\zeta_N))$ are maps from the set of embeddings $\Q(\zeta_N)\to\C$
with values in $E$.
We choose generators
\begin{align*}
 t_{\DR}(\chi)&= p_{\chi^{-1}}\left(\zeta_N\right)\in V_{\DR}(\chi)\subset
 \Q(\zeta_N)\tensor E\\
 t_B(\chi)& =  p_{\chi^{-1}}( \delta_{\sigma_0})\in
p_{\chi^{-1}}(V_B(\chi)))\\
t_p(\chi)&=t_B(\chi)\tensor 1\in V_B(\chi)\tensor E_p=V_p(\chi)
\end{align*}
Note that $t_\DR(\chi)$ is a generator of 
$p_{\chi^{-1}}(\Z[\mu_N]\tensor \Oh)$ as this group is generated by the set $p_{\chi^{-1}}(\zeta_N^i)$ for 
$i\in \Z/N$ and
these elements either vanish or are multiples of $p_{\chi^{-1}}(\zeta_N)$.

\begin{lemma} Let $\chi$ be a Dirichlet character of conductor $N$.
The element 
\[ \epsilon\in \det^{-1} V_B(\chi)(r)^+ \tensor\det V_{\DR}(\chi)(r)
  \tensor \det V_B(\chi)^\lor(1-r)^+ \]
which is mapped to $\frac{L(\chi^{-1},1-r)^*}{L(\chi,r)^*}$ under the  
isomorphism of the determinant with $E_{\infty}$ defined in \ref{localconj} is
\[ 2 {N^{r-1}(r-1)!}  \left(t_B(\chi)(r)^{+}\right)^{-1} \tensor  t_{\DR}(\chi)(r)
  \tensor  t_B(\chi)^\lor(1-r)^+ \]
(up to sign depending only on $r$).
\end{lemma}
\bew
By the functional equation for Dirichlet-$L$-functions
\[ \frac{L(\chi^{-1},1-r)^*}{L(\chi,r)^*}=
\pm 2\frac{(r-1)!N^r}{\tau(\chi) (2\pi i)^{r-\delta}} \]
 where 
\[ \tau(\chi)= \sum_{a=1}^N\chi(a)\zeta_N^a=
  \sum_{\sigma\in G}\chi(\sigma^{-1})\sigma(\zeta_N) \]
 is the Gauss sum and
$\delta=0$ (resp. $\delta=1$) if $\chi$ and $r$ have the same (different)
parity. On the other hand we have to consider
\[ 0\to V_B(\chi)(r)^+\tensor\R\to V_\DR(\chi)\tensor\R\to 
V_B(\chi)(r-1)^+\tensor\R \to 0\]
If $r$ and $\chi$ have the same (different) parity, then only the 
last (first) Betti-term appears. The comparison morphism maps
$t_\DR(\chi)$ to
\begin{multline*}
 \sum_{\sigma}\sigma_0\sigma(p_{\chi^{-1}}\zeta_N) \delta_{\sigma_0\sigma}
= \sum_{\sigma} \chi^{-1}(\sigma)\sigma_0(\frac{ \tau(\chi^{-1})}{\Phi(N)})\sigma^{-1}\delta_{\sigma_0} \\
= \frac{\tau({\chi^{-1}})}{(2\pi i)^{r-\delta}} t_B(\chi)(r-\delta)
= \frac{N}{\tau({\chi})(2\pi i)^{r-\delta}}t_B(\chi)(r-\delta)
\end{multline*}
\bewende

Hence conjecture \ref{localconj} predicts that 
\[ 2(r-1)!N^{r-1}t_\DR(\chi) \]
is a generator of $\det^{-1} R\Gamma(\Q_p,T_p(\chi)(r))$ under 
$\exp_p$ and the local identification \ref{localident}. 

\begin{prop}\label{B.1.2}
Conjecture \ref{localconj} holds in the case $K=\Q$ (i.e., non-equivariantly), $r=1$, $p\neq 2$, $\chi$ restricted to the Galois group of $\Q_p$ is the trivial character.
\end{prop}
\bew 
As $\chi$ restricted to the local Galois group is trival, $p$ splits completely.
Hence
$\Z[\mu_N]\tensor \Z_p\isom \bigoplus_{v\mid p}\Z_p$. 
As part of the class number case  (with $F=\Q$) we have checked in section \ref{secclassnumber} that 
the lemma holds for the trivial character. Taking the sum over all $v$, this implies even
that $ \det_{\Z_p[G]}(\Z[\mu_N]\tensor \Z_p)$ is identified with
$\det^{-1}_{\Z_p[G]}R\Gamma(\Q_p,\bigoplus_{v\mid p}\Z_p(1))$ where $G=\Gal(\Q(\mu_N)/\Q)$. Tensoring both
determinants with $\Oh(\chi)$ over $\Z_p[G]$ implies that a generator of 
$p_{\chi^{-1}}(\Z[\mu_N]\tensor \Z_p)$ is mapped to a generator of
$\det^{-1}_{\Oh_p} R\Gamma(\Q_p,\Oh_p(\chi)(1))$.
As remarked before, $t_\DR(\chi)$ is a generator of the $\Oh$-module $p_{\chi^{-1}}(\Z[\mu_N]\tensor  \Oh)$, hence also of the $\Oh_p$-module 
$p_{\chi^{-1}}(\Z[\mu_N]\tensor  \Oh_p)$.
\bewende

If $K/\Q$ is an abelian extension, $X(K/\Q)$ the set
of characters of $\Gal(K/\Q)$, then the equivariant version of conjecture \ref{localconj}
predicts that
\[ 2(r-1)!\left(N_{\omega\chi}^{r-1}t_\DR(\omega\chi)\right)_{\omega\in X(K/\Q)} \]
is a generator of $\det^{-1}_{\Oh_p[\Gal(K/\Q)]} R\Gamma(K_p,T_p(\chi)(r))$ under
$\exp_p$ and the local identification \ref{localident}. 

\begin{prop}\label{localconjequiv}
Conjecture \ref{localconj} holds for all Dirichlet characters, $r>1$ and $p\neq 2$
equivariantly for the cyclotomic extension $\Q(\mu_{p^n})/\Q$.
\end{prop}

The proof will cover the rest of this appendix. 

\begin{lemma}\label{formula}
If $r>1$, then conjecture \ref{localconj} holds equivariantly if and only if
\[ \left(\frac{2(r-1)!N_{\omega\chi}^{r-1}(1-\omega^{-1}\chi^{-1}(p)p^{r-1})}{1-\omega\chi(p)p^{-r}}\exp_p(t_{\DR}(\omega\chi))\right)_{\omega\in X(K/\Q)} \]
is a generator of $\det_{\Oh_p[\Gal(K/\Q)]}^{-1}R\Gamma(K_p,T_p(\chi)(r))$.
\end{lemma}

\bew It suffices to compute the image of our generator characterwise. We
write $\chi$ for $\omega\chi$.
Recall that $\exp_p: V_\DR(\chi)(r)\to R\Gamma(\Q_p,V_p(\chi)(r))[1]$
is a quasi-isomorphism.
The whole point of the lemma is that the identification of determinants in \ref{localident}
is {\em not} the one induced by the above quasi-isomorphism. 
It uses
\begin{align*}
 \det^{-1} R\Gamma(\Q_p,V_p(\chi)(r))=&\det^{-1}R\Gamma_f(\Q_p,V_p(\chi)(r))\tensor
\det^{-1} R\Gamma_{/f}(\Q_p,V_p(\chi)(r)) \\
= & \det^{-1}R\Gamma_f(\Q_p,V_p(\chi)(r))\tensor
\det R\Gamma_f(\Q_p,V_p(\chi^{-1})(r))
\end{align*}
and  \ref{localident} on both factors.
 The first factor is quasi-isomophic to $R\Gamma(\Q_p,V_p(\chi)(r))$, the second to zero.
By proposition \ref{localfactor} the lattice of \ref{localident} differs from
the one given by 
\[  \det^{-1} R\Gamma(\Q_p,T_p(\chi)(r))\tensor \det 0 \]
by the local Euler factors $1-\chi(p)p^{-r}$ and     $(1-\chi^{-1}(p)p^{r-1})^{-1}$.
\bewende


\subsection{Homological algebra}
From now on, $\chi$ is considered as a local character of conductor $N=N'p^m$
where $N'$ is prime to $p$. 
Let $K=\Q_p(\mu_N)$. Let $K_0$ be its totally unramified subfield and
$K_n=K_0(\mu_{p^n})$ (this differs from the convention in section \ref{mc}).
 In particular, $K_0=\Q_p(\mu_{N'})$.

Let
$\Gamma=\Gal(K_{\infty}/K_0)$, $\Gamma_{n}=\Gal(K_\infty/K_n)$ and 
$H=\Gal(K_0/\Q_p)$. We put $G_n=\Gal(K_n/K_0)=\Gamma/\Gamma_n$.
In this section we let $\Lambda:=\prolim_n\Oh_p[ [G_n] ]$, which differs from the $\Lambda$ in the 
main text.
Our first step follows \cite{Bloch-Kato} 4.2 p. 367. We put
$P=\Gal(K_\infty^{ab}/K_\infty)$, $U=\prolim \Oh_{K_0}[\zeta_{p^n}]^*$
(or rather their pro-$p$-parts). They
are natural $\Lambda$-modules. There are exact sequences 
\[ 0\to U\to P\to \Z_p \to 0\]
and
\[ 0\to \Z_p(1)\to U\xrightarrow{\Col} \Oh_{K_0}[[\Gamma]]\to \Z_p(1)\to 0 \]
where $\Col$ is the Coleman map, see \cite{Bloch-Kato} 4.2 for more details.

\begin{lemma}
The following $\Lambda$-modules are pseudo-isomorphic:
\[ 
H^1(K_\infty,T_p(\chi)(r))\to \Hom(P,T_p(\chi)(r))^H 
\]
There are natural isomorphism of $\Lambda$-determinants
\begin{align*}
\det R\Gamma(\Q_{p,\infty},T_p(\chi)(r)) [1]&=
\det^{-1} H^0(\Q_{p,\infty},T_p(\chi)(r))\tensor \det \Hom(P,T_p(\chi)(r))^H 
\\
 &= \det \Hom(U,T_p(\chi)(r))^H \\
&= \det \Hom_H(\Oh_{K_0}[[\Gamma]],T_p(\chi)(r))
\end{align*}
\end{lemma}
\bew If $H$ has no $p$-torsion, then this is precisely \cite{Bloch-Kato} 4.2,
p. 367 and follows immediately from the short exact sequences. Hence
we consider only the case where $H$ has $p$-torsion. 
This means that $\chi(r)$ is non-trivial over $\Q_{p,\infty}$ and
that $H^0(\Q_{p,\infty},T_p(\chi)(r))=0$.

The following arguments are very ugly, we apologize for this. We
abbreviate $T=T_p(\chi)(r)$.
The crucial observation is finiteness of
\[  H^i(H,\Hom(P,T)) \]
for $i>0$.
We have the short exact sequence
\[ 0\to \Hom(\Z_p,T)\to \Hom(P,T)\to \Hom(U,T)\to 0 \]
Hence
\[ H^1(H,T)\to  H^1(H,\Hom(P,T))\to H^1(H,\Hom(U,T))\to H^2(H,T) \]
The first and last term are finite. We replace $P$ by $U$. We consider the triangle
\[ \Z_p(1)\to U\to  C\]
with $C=[\Oh_{K_0}[[\Gamma]]\to \Z_p(1)]$
Hence
\[ 0\to \Hom(C,T)\to \Hom(U,T)\to \Hom(\Z_p(1),T) \]
This implies that the quotient is at most a finitely generated $\Z_p$-module without $H$-invariants and
hence
\[ H^1(H,\Hom(C,T))\to H^1(H,\Hom(U,T)) \]
is again a pseudo-isomorphism.
Finally, we have
\[ 0=\Hom(\Z_p(1)[-1],T)\to \Hom(C,T)\to
\Hom(\Oh_{K_0}[[\Gamma]],T)\to \Hom(\Z_p(1)[-1],T[1])=0
\] 
and we have reduced the question to
$H^1(H,\Hom(\Oh_{K_0}[[\Gamma]],T))$. But this is an induced $H$-module as
 $\Oh_{K_0}$ is the free $\Z_p[H]$-module
by existence of a normal basis. 
Hence its $H$-cohomology vanishes.
The same argument works for $H^2$.

Now we can turn to the proof of the lemma. 
Recall that
$H^1(K_\infty,T_p(\chi)(r))=\Hom(P,T)$ as the coefficients become trivial over $K_\infty$. Hence
\[ 0\to H^1(H,\Hom(P,T))\to H^1(\Q_{p,\infty},T)
\to H^0(H,\Hom(P,T))
\to H^2(H,\Hom(P,T))
\]
This is the desired pseudo-isomorphism. Pseudo-isomorphic $\Lambda$-modules have the same determinant, so
this proves the first identification of determinants. By similar arguments
applying $H^i(H,\cdot)$ to the exact sequences and using $T^H=0$, the other identifications
follow.
\bewende

\begin{cor}\label{generatorcomp}
Let $\chi$ be a local character and $r>1$  as before. Let 
$G_n=\Gal(\Q(\mu_{p^n})/\Q)$.
Then there is an isomorphism of determinants
\[ \det_{\Oh_p[G_n]} R\Gamma(\Q(\mu_{p^n})_p,T_p(\chi)(r))[1]\isom 
\det_{\Oh_p[G_n]} \Hom_{\Gamma_n\times H}(\Oh_{K_0}[[\Gamma]],T_p(\chi)(r)) \]
which rationally is induced  by the isomorphism
\[ 
  H^1(\Q(\mu_{p^n})_p,V_p(\chi)(r)) \xrightarrow{s_\chi} \Hom_{\Gamma_n\times H}(K_0[[\Gamma]],V_p(\chi)(r))\]
\end{cor}
\bew
\[ R\Gamma(\Q(\mu_{p^n})_p,T_p(\chi)(r))= R\Gamma(\Gamma_n,R\Gamma(\Q_{p,\infty},T_p(\chi)(r)) \]
Hence by the lemma
\[ \det_{\Oh_p[G_n]} R\Gamma(\Q(\mu_{p^n})_p,T_p(\chi)(r))[1] 
 = \det_{\Oh_p[G_n]} R\Gamma( \Gamma_n,\Hom_H(\Oh_{K_0}[[\Gamma]],T_p(\chi)(r)) \]
In this case there are no $\Gamma_n$-coinvariants, hence 
\[ R\Gamma( \Gamma_n,\Hom_H(\Oh_{K_0}[[\Gamma]],T_p(\chi)(r))=\Hom_H(\Oh_{K_0}[[\Gamma]],T_p(\chi)(r))^{\Gamma_n}  \]

The identification of determinants of Iwasawa modules was induced by pseudo-isomorphisms hence $R\Gamma(\Gamma_n,\cdot)$ of the error terms is torsion. So rationally it is induced by a quasi-isomorphism.
\bewende

\subsection{Reciprocity laws}
Let $N=N'p^m$ with $N'$ prime to $p$.
Let $\zetatilde_{N'}$ be a generator of the free $\Z_p[H]$-module $\Oh_{K_0}$.
We can assume that $\zetatilde_N$ is $\zeta_N$ plus a linear combination of
roots of unity of order less than $N$. In particular
$p_{\chi}(\zetatilde_{N'})=p_{\chi}(\zeta_{N'})$. 
Evaluation in $\zetatilde_{N'}$ 
fixes an isomorphism
$\Hom_{\Gamma\times H}( K_0[[\Gamma]],V_p(\chi)(r))\isom V_p(\chi)(r)$.

\begin{lemma}\label{indexcomp}
\[ (s_\chi\exp_p)(t_\DR(\chi))=  \frac{{N'}^r(1-\chi(p) p^{-r})}{N^{r-1}(r-1)!(1-p^{r-1}\chi^{-1}(p))}t_p(\chi)(r).
\]
\end{lemma}
\bew Recall that $\Gamma/\Gamma_m=G_m$.
We have $\Z_p[H]\zetatilde_{N'}\isom \Oh_{K_0}$. Hence $\zetatilde_{N'}$ also defines 
isomorphisms
\[ \Hom_{\Gamma_m}(\Oh_{K_0}[[\Gamma]],E_p(r))\isom 
 \Hom(\Oh_{K_0}[G_m],E_p(r))\isom
\Hom(\Z_p[H\times G_m],E_p(r))\]
The following diagram commutes
\[\begin{CD}
 V_\DR(\chi)@>s_\chi\exp_p>> \Hom_{\Gamma\times H}( \Oh_{K_0}[[\Gamma]],V_p(\chi)(r))@>>> V_p(\chi)(r)\\
@VVV @VV\iota V @VV\iota' V\\
K_m @>s_{K_m}\exp_p>> \Hom_{\Gamma_m}(\Oh_{K_0}[[\Gamma]],E_p(r))@>>>
  \Hom(\Z_p[H\times G_m],E_p(r))
\end{CD}\]
where $\iota'$ is defined via the inclusion as $p_{\chi^{-1}}$-eigenpart.
Under this inclusion $\iota'(t_p(\chi)(r))= p_{\chi^{-1}}\delta$ where 
$\delta(e)=t_p(r)$ (the standard generator of $E_p(r)$) and $\delta(g)=0$ for $e\neq g\in H\times G_m$.
Hence
\[ \iota( t_p(\chi)(r) )( \zetatilde_{N'} )= \iota'( t_p(\chi)(r)) (e)= \frac{1}{\Phi(N)} \; t_p(r).
\]
This implies that for 
$\alpha\in V_{\DR}(\chi)$ we have
\[
 (s_\chi\exp_p)(\alpha)(\zetatilde_{N'}) = 
 (s_{K_m}\exp_p)(\alpha)(\zetatilde_{N'})\cdot \frac{\Phi(N)}{t_p(r)}\cdot t_p(\chi)(r) 
\]

The map $s_{K_m}:K_m \to \Hom_{\Gamma_m}(K_0[[\Gamma]],E_p(r))$ is 
computed by the explicit reciprocity law for $K=K_0(\mu_{p^{m}})$.
In the unramified case, $s_{K_0}$ reduces to a map
$K_0\to \Hom(K_0,E_p(r))$ . By
\cite{Bloch-Kato} Claim 4.8 (p. 368) it is given by
\[\alpha\mapsto\left( 
   \beta\mapsto
 \frac{1}{(r-1)!} \Tr_{K_0/\Q_p}( (1-p^{r-1}\Fr_p^{-1})^{-1})(\beta)\cdot (1-\Fr_p^{-1} p^{-r})(\alpha) \; t_p(r)\right)
\]
(Recall that $ \Fr_p$ is geometric Frobenius whereas $f$ in loc. cit is
arithmetic Frobenius).

Evaluating the formula with $\alpha=p_{\chi^{-1}}(\zeta_N)$ and $\beta=\zetatilde_N=\zetatilde_{N'}$, we
get 
\[ \frac{N(1-\chi(p) p^{-r})}{(r-1)!(1-p^{r-1}\chi^{-1}(p))\Phi(N)} \; t_p(r)\]
The factor $\Phi(N)/t_p(r)$ cancels. This yields the desired formula in the unramified case.

In the ramified case $m\geq 1$, we follow an approach also used 
 by Benois and Nguyen in an earlier version of
\cite{BeNg}.
We use Kato's higher explicit reciprocity law, \cite{galaxy} 2.1.7.
The map $s_{K_m}\exp_p$ is given by
\[
\alpha \mapsto \left( 
   \beta\mapsto
 \frac{1}{(r-1)!} p^{-rm}\Tr_{K_m/\Q_p}(\alpha\cdot (1-p^{r-1}\phi)^{-1}D^r(
\beta(1+T))(\zeta_{p^m}-1) t_p(r) \right)
\]
where we consider $\beta\cdot(1+T)\in \Oh_{K_0}[[T]]$ and
$D=(1+T)d/dT$, $\phi$ Frobenius on $\Oh_{K_0}[[T]]$ (see also also
\cite{PerrinReciprocity} for facts about Coleman power series.)
We  evaluate with  $\alpha= p_{\chi^{-1}}(\zeta_N) $ and 
$\beta=\zetatilde_{N'}$ 
and get
\[ \frac{1}{(r-1)! p^{rm}} \Tr_{K_m/\Q_p}(p_{\chi^{-1}}(\zeta_N)
 p_{\chi}y_m)\; t_p(r)
\]
with $y_m=[(1-\phi p^{r-1})^{-1}\zetatilde_{N'}(1+T)](\zeta_{p^m}-1)$. 
Using the geometric series for the inverse of $1-p^{r-1}\phi$, we get
\[ y_m= \zetatilde_{N'}\zeta_{p^m} + p^{r-1} \Fr_p^{-1}(\zetatilde_{N'})\zeta_{p^m}^p+
p^{2(r-1)} \Fr_p^{-2}(\zetatilde_{N'})\zeta_{p^m}^{2p}+
\dots \]
Only the first summand contributes to the $\chi$-part as we have
assumened that $\chi$ is primitive of level $N$. Hence
\[ p_{\chi}(y_m)=p_{\chi}(\zetatilde_{N'}\zeta_{p^m})
=p_{\chi}(\zeta_{N'}\zeta_{p^m})=
p_{\chi}(\zeta_{N})
 \]
As $p_{\chi^{-1}}(\zeta_N)p_{\chi}(\zeta_N)=N/\Phi(N)^2$,
this yields
\[s_{K_m}\exp_p(p_{\chi^{-1}}(\zeta_N))(\zetatilde_{N'}) =\frac{N}{(r-1)! p^{rm}\Phi(N)}t_p(r) \]
The factor $\Phi(N)/t_p(r)$ again cancels and $p^{rm}$ is the $p$-part of $N^r$. 
\bewende

\noindent{\em Proof of proposition \ref{localconj}:}
As before let
 $\chi$ be a character of conductor $N=p^mN'$ with $N'$ prime to $p$,
$G_n=\Gal(\Q(\zeta_{p^n})/\Q)$ and $K_0=\Q_p(\zeta_{N'})$. 
We have to check that the element of lemma \ref{formula} (call it $\epsilon'$) is a generator of
$\det_{\Oh_p[G_n]}^{-1}R\Gamma(\Q(\zeta_{p^n})_p, T_p(\chi)(r)$. By 
 lemma \ref{generatorcomp} this is equivalent to showing that the image
of $\epsilon'$
in $\Hom_{\Gamma_n\times H}(\Oh_{K_0}[[\Gamma]],T_p(\chi)(r))
=\Hom_H(\Oh_{K_0}[G_n],T_p(\chi)(r))$
is an $\Oh_p[G_n]$-generator of the latter. As
\[ \Hom_{H}(\Oh_{K_0}[G_n],T_p(\chi)(r))\subset
   \bigoplus_\omega V_p(\chi\omega)(r)\]
where $\omega$ runs through all characters of $G_n$, the image of $\epsilon'$ is
uniquely determined by its $\omega$-components. By lemma \ref{indexcomp}
applied to the character $\chi\omega$ and the definition of $\epsilon'$
\[ (s_{K_n}\epsilon')_{\omega}= 2(N')^r t_p(\chi\omega)(r) .\]
(Note that $N'$, $H$ and $\zetatilde_{N'}$ are the same for all these $\chi\omega$).
On the other hand, $\Hom_{H}(K_0[G_n],V_p(\chi)(r))$ has
a standard generator given by the function $\delta$ with
$\delta(\zetatilde_{N'}e)=t_p(\chi)(r)$, $\delta(\zetatilde_{N'}g)=0$ for 
$g\neq e\in G_n$. It has $\omega$-component $t_p(\chi\omega)(r)$. Hence
\[ s_{K_n}\epsilon'= 2(N')^r\delta .\]
As $N'$ is prime to $p$ and $p\neq 2$
this finishes the proof of proposition \ref{localconjequiv} for Dirichlet 
characters.

\end{appendix}

\end{document}